\magnification=\magstep1
\input amstex
\documentstyle{amsppt}
\NoBlackBoxes \vsize=7in \hsize=5.5in

\catcode`\@=11 \loadmathfont{rsfs}
\def\mycal{\mathfont@\rsfs}
\csname rsfs \endcsname \topmatter
\title AMALGAMATED FREE PRODUCTS OF $w$-RIGID FACTORS \\
AND CALCULATION OF THEIR SYMMETRY GROUPS\endtitle \vskip -.3in
\author A. IOANA, J. PETERSON and S.
POPA \endauthor \rightheadtext{Amalgamated free products} \vskip
-.3in \affil University of California, Los Angeles\endaffil \vskip
-.3in
\address Math.Dept., UCLA, LA, CA 90095-155505\endaddress
\email adiioana@math.ucla.edu, \, jpete@math.ucla.edu, \,
popa@math.ucla.edu\endemail

\thanks J.P. supported
in part by NSF VIGRE Grant DMS 9983726 at UCLA; S.P. supported in
part by NSF-Grant 0100883.\endthanks \vskip -.3in \abstract We
consider amalgamated free product II$_1$ factors $M = M_1 *_B M_2
*_B ...$ and use ``deformation/rigidity'' and ``intertwining''
techniques to prove that any relatively rigid von Neumann subalgebra
$Q\subset M$  can be unitary conjugated into one of the $M_i$'s. We
apply this to the case $M_i$ are w-rigid II$_1$ factors, with $B$
equal to either $\Bbb C$, to a Cartan subalgebra $A$ in $M_i$, or to
a regular hyperfinite II$_1$ subfactor $R$ in $M_i$, to obtain the
following type of unique decomposition results, \`a la Bass-Serre:
If $M = (N_1 *_C N_2 *_C ...)^t$, for some $t>0$ and some other
similar inclusions of algebras $C\subset N_j$ then, after a
permutation of indices, $(B\subset M_i)$ is inner conjugate to
$(C\subset N_i)^t$, $\forall i$. Taking $B=\Bbb C$ and $M_i =
(L(\Bbb Z^2 \rtimes \Bbb F_{2}))^{t_i}$, with $\{t_i\}_{i\geq 1}=S$
a given countable subgroup of $\Bbb R_+^*$, we obtain continuously
many non stably isomorphic factors $M$ with fundamental group
$\mycal F(M)$ equal to $S$. For $B=A$, we obtain a new class of
factors $M$ with unique Cartan subalgebra decomposition, with a
large subclass satisfying $\mycal F(M)=\{1\}$ and Out$(M)$ abelian
and calculable. Taking $B=R$, we get examples of factors with
$\mycal F(M)=\{1\}$, Out$(M)=K$, for any given separable compact
abelian group $K$.
\endabstract
\endtopmatter
\document
\vskip -.2in \centerline{\bf Content}

\settabs\+\indent\indent & 10. \ & 1.1 \quad & More on rigidity
\cr \tenpoint{\+& 0. & Introduction \cr \+& 1. & Conjugating
subalgebras in AFP factors  \cr \+& 2. & Deformations of AFP
factors \cr \+& 3. & Deformation/rigidity arguments \cr \+& 4. &
Existence of intertwining bimodules \cr \+& 5. & Rigid subalgebras
in AFP factors: General Bass-Serre type results \cr \+& 6. &
Amalgamation over $\Bbb C$: Free product factors with prescribed
$\mycal F(M)$ \cr \+& 7. & Amalgamation over Cartan subalgebras:
vNE/OE rigidity results \cr \+& 8. & Amalgamation over $R$:
Factors with no outer automorphisms \cr \+& Appendix: Constructing
freely independent actions \cr}

\newpage

\heading 0. Introduction. \endheading

We prove in this paper a series of rigidity results for
amalgamated free product (hereafter abbreviated AFP) II$_1$
factors $M=M_1 *_B M_2$ which can be viewed as von Neumann algebra
versions of the ``subgroup theorems'' and ``isomorphism theorems''
for AFP groups in Bass-Serre theory. Our main ``subalgebra
theorem'' shows that, under rather general conditions, any von
Neumann subalgebra $Q \subset M$ with the {\it relative property}
(T) in the sense of ([P5]) (also called a {\it rigid inclusion}),
can be conjugated by an inner automorphism of $M$ into either
$M_1$ or $M_2$. We derive several ``isomorphism theorems'' in the
case the amalgamation is over the scalars, $B=\Bbb C$, over a
common Cartan subalgebra, $B=A$, or over a regular hyperfinite
subfactor, $B=R$. The typical such statement shows that if
$\theta: M \simeq N^t$ is an isomorphism from an AFP factor $M=M_1
*_B M_2 *_B ...*_B M_n$ onto the amplification  by some $t > 0$ of
an AFP factor $N=N_1 *_C N_2 *_C ...*_C N_n$, $1 \leq m,n \leq
\infty$, with each $M_i, N_j$ containing a ``large'' subalgebra
with the relative property (T), then $m=n$ and $\theta(B\subset
M_i)$ is unitarily conjugate to $(C\subset N_i)^t$, $\forall i$,
after some permutation of indices.

When applied to the case $B=R$, these results allow us to obtain
the first explicit calculations of outer automorphism groups of
II$_1$ factors and answer in the affirmative a problem posed by A.
Connes in 1973, on whether there exist II$_1$ factors $M$ with no
outer automorphism, i.e. with Out$(M)\overset \text{\rm def} \to =
\text{\rm Aut}(M)/\text{\rm Int}(M)=\{1\}$. More precisely, we
show that if a group $\Gamma$ is the free product of two infinite
property (T) groups with no non-trivial characters, for example
$\Gamma=SL(n_0, \Bbb Z) * SL(n_1, \Bbb Z)$, $n_0, n_1 \geq 3$,
then there exist actions of $\Gamma$ on the hyperfinite II$_1$
factor $R$ such that the corresponding crossed product factors
$M=R \rtimes \Gamma$ have both trivial fundamental group, $\mycal
F(M)=\{1\}$, and trivial outer automorphism group, Out$(M)=\{1\}$.
In fact, the general result shows that for any separable compact
abelian group $K$ there exist factors $M$ with $\mycal F(M)=\{1\}$
and Out$(M)=K$.

In turn, when applied to the case of amalgamated free products
over a common Cartan subalgebra, our ``isomorphism theorem''
provides a Bass-Serre type result for {\it orbit equivalence} (OE)
of actions of free product groups $\Gamma=\Gamma_1 * ...
*\Gamma_n$, $\Lambda = \Lambda_1 * ... * \Lambda_m$ on the
probability space. Thus, we show that if each $\Gamma_i,
\Lambda_j$ has an infinite normal subgroup with the relative
property (T) of Kazhdan-Margulis (for instance $\Gamma_i,
\Lambda_j$ Kazhdan groups, $\forall i,j$) and $(\sigma, \Gamma)$,
$(\theta,\Lambda)$ are free, probability measure preserving (m.p.)
actions with $\sigma_{|\Gamma_i}, \theta_{|\Lambda_j}$ ergodic
$\forall i,j$, then $\sigma\sim_{OE} \theta$ implies $m=n$ and
$\sigma_{|\Gamma_i}\sim_{OE} \theta_{|\Lambda_i}$, $\forall i$,
after a permutation of the indices $i$. Note that the opposite
implication holds true for arbitrary groups $\Gamma_i, \Lambda_j$,
as shown by D. Gaboriau in ([G2]). In fact, we derive the
component by component OE of actions under the weaker assumption
that the group measure space factors associated with $(\sigma,
\Gamma), (\theta, \Lambda)$ are stably isomorphic, i.e. when
$\sigma,\theta$ are {\it von Neumann equivalent} (vNE). We use
this vNE Bass-Serre rigidity and ([MS], [Ge1,2], [Fu1,2], [P6,8])
to give examples of group measure space factors $M$ from free
ergodic m.p. actions $\sigma$ of free product groups
$\Gamma=\Gamma_1 * \Gamma_2 *...$ such that $\mycal F(M) = \{1\}$
and Out$(M)= \text{\rm H}^1(\sigma,\Gamma)$, with explicit
calculation of the abelian group $\text{\rm H}^1(\sigma,\Gamma)$.

Finally, when applied to the case $B = \Bbb C$ our results become
von Neumann algebra analogue of Kurosh's classical theorems for
free products of groups, similar to N. Ozawa's recent results of
this type in ([O2]), but covering a different class of factors
than ([O2]) and allowing amplifications. For instance, we show
that if $N_i, 2\leq i \leq n$, $M_j, 2\leq j \leq m,$ are property
(T) II$_1$ factors in the sense of Connes-Jones (e.g. $N_i, M_j$
group factors associated with Kazhdan groups, [CJ]) then $M_1 *
M_2 * ... *M_m \overset \theta \to \simeq (N_1*N_2*... *N_n)^t$
implies $m=n$ and $\theta(M_i)$ inner conjugate to $N_i^t$,
$\forall i$, after some permutation of indices. In fact, in its
most general form our result only requires $M_i, N_j$ to be {\it
weakly rigid} ({\it w-rigid}), i.e. to have diffuse, regular
subalgebras with the relative property (T). Taking $M=N$ and $M_j
= P^{s_j}$, with $\{s_j\}_j = S$ a multiplicative subgroup of
$\Bbb R_+^*$ and $P$ a w-rigid II$_1$ factor with trivial
fundamental group (for instance the group factor $L(G)$ associated
with $G=\Bbb Z^2 \rtimes SL(2,\Bbb Z)$, cf. [P5]) and using a
result of Dykema-Radulescu ([DyR]), we get $\mycal F(M)=S$ for
$M=*_{s\in S}P^s$. This provides a completely new class of factors
with arbitrary given $S \subset \Bbb R_+^*$ as fundamental group
than the ones in ([P6]). Indeed, the examples constructed in
([P6]) are group measure space factors, while the free group
factors $*_{s\in S}P^s$ have no Cartan subalgebra, by results of
Voiculescu ([V2]; see [Sh] and Remark 6.6 in this paper).

The key technical result behind all these applications is the above
mentioned ``subalgebra theorem'', of Bass-Serre type. We state it in
details below, together with other main results in the paper, and
also explain some of the ideas behind the proofs. An inclusion of
finite von Neumann algebras $B \subset P$ will be called {\it
homogeneous} if there exist $\{y_j\}_j \subset P$ with
$E_B(y_i^*y_j)=\delta_{ij}$, $\forall i,j$ and $\Sigma_i y_j B$
dense in $P$. This technical assumption is satisfied by all
inclusions coming from (cocycle) crossed products and (generalized)
group measure constructions, or Cartan inclusions. It is also
satisfied when $P$ is an arbitrary finite von Neumann algebra and
$B=\Bbb C$. Following ([P5]), a von Neumann subalgebra $Q\subset P$
has the {\it relative property} (T) (or $Q \subset P$ is a {\it
rigid inclusion}) if any ``deformation'' of $id_P$ by completely
positive, sub-unital, sub-tracial maps, $\phi_n \rightarrow id_P$,
is uniform on the unit ball of $Q$ (see also [PeP]).

\proclaim{0.1. Theorem} Let $(M_i, \tau_i)$, $i=1,2$, be finite
factors with a common von Neumann subalgebra $B \subset M_i$, such
that $\tau_{1|B}=\tau_{2|B}$ and such that $B \subset M_i$ are
homogeneous, $i=1,2$. Let $Q\subset M=M_1 *_B M_2$ be a diffuse
von Neumann subalgebra with the relative property $(\text{\rm T})$
such that no corner $qQq$ of $Q$ can be embedded into $B$. Then
there exists a unique partition of $1$ with projections $q'_1,
q'_2$ in the commutant of $Q$ in $M$ such that $u_i(Qq_i')u_i^*
\subset M_i$, $i=1,2$, for some unitaries $u_1, u_2$ in  $M$.
Moreover, if the normalizer of $Q$ in $M$ generates a factor $N$,
then there exists a unique $i \in \{1,2\}$ such that $uQu^*
\subset M_i$ for some $u \in \Cal U(M)$, which also satisfies
$uNu^* \subset M_i$.
\endproclaim

The proof of this result takes Sections 2 through 5 of the paper.
It uses ``deformation/rigidity''  and ``intertwining'' techniques
from ([P3,5,6]). Thus, we embed $M=M_1 *_B M_2$ into the larger
algebra $\tilde{M}=M *_B (B\overline{\otimes} L(\Bbb F_2))$, whose
aboundance of deformations is used to show that ``rigid parts'' of
$M$ have to concentrate on certain subspaces with ``bounded
word-length''. This initial information is then used as a starting
point in a word-reduction argument to obtain a Hilbert bimodule
intertwining $Q$ into one of the $M_i$'s. The homogeneity
condition is needed in order to measure the ``size'' of letters in
$M_i$'s. To get from this a unitary element conjugating $Q$ into
$M_i$, we prove in Section 1 a series of results on the relative
commutants and normalizers of subalgebras in AFP factors, using
(2.1, 2.3 in [P6]).

If we take $B=\Bbb C$ in the Theorem 0.1 and use the fact that
finite von Neumann algebras with Haagerup property have no diffuse
subalgebras with the relative property (T), then we get an
analogue of Kurosh's isomorphism theorem for free products of
groups:

\proclaim{0.2. Theorem} Let $(M_0, \tau_{M_0}), (N_0, \tau_{N_0})$
be finite von Neumann algebras with Haagerup's compact
approximation property. Let $M_i, 1\leq i \leq m,$ and $N_j, 1\leq
j \leq n,$ be w-rigid $\text{\rm II}_1$ factors, where $ m,n \geq
1$ are some cardinals $($finite or infinite$)$. If $\theta$ is an
isomorphism of $M=*_{i=0}^{m}M_i$ onto $N^t$, where $N=*_{j=0}^{n}
N_j$ and $t
> 0$, then $m=n$ and, after some permutation of indices,
$\theta(M_i)$ and $N^t_i$ are unitary conjugate in $N^t$, $\forall
i\geq 1$.
\endproclaim

Ozawa's pioneering result of this type in ([O2]) concerns free
products of group factors $M_i=L(\Gamma_i), N_i=L(\Lambda_i)$ with
each $\Gamma_i, \Lambda_i$ a product of two or more ICC groups,
either word hyperbolic (at least one of them) or amenable, typical
examples being the groups $\Bbb F_{n_i} \times S_{\infty}$, not
covered by 0.2 above. In turn, our typical $M_i, N_i$ are factors
from property (T) (more generally w-rigid) groups.

Letting $M_i=N_i, \forall i,$ and $m < \infty$ in Theorem 0.2, it
follows that if $\mycal F(M_i)=\{1\}$ for some $1\leq i\leq m$
(for example if $M_i=L(\Bbb Z^2 \rtimes \Bbb F_k)$, with $\infty >
k \geq 2$, cf. [P5]), then $\mycal F(M)=\{1\}$. Moreover, taking
$m=\infty$ in 0.2 and using the ``compression formula'' for free
products of infinitely many II$_1$ factors $(*_iM_i)^t\simeq *_i
M_i^t$ in ([DyR]), we can include specific numbers into the
fundamental group. Thus we get:

\proclaim{0.3. Corollary} $1^\circ$. Let $m \in \Bbb N$ and let
$M_1, \ldots, M_m$ be w-rigid $\text{\rm II}_1$ factors. Let
$(M_0, \tau_{M_0})$ be a finite von Neuman algebra with Haagerup's
compact approximation property. If one of the factors $M_i, m\geq
i \geq 1,$ has trivial fundamental group then so does $M=*_{i =
0}^m M_i$.

$2^\circ$. If $S \subset \Bbb R_+^*$ is an arbitrary infinite
$($possibly uncountable$)$ subgroup and $P$ is a w-rigid
$\text{\rm II}_1$
factor with trivial fundamental group $($e.g. $P= L(\Bbb Z^2
\rtimes SL(2,\Bbb Z)))$ then the $\text{\rm II}_1$ factor
$*_{s\in S} P^s$ has fundamental group equal to $S$.
\endproclaim

Since a group measure space factor $M=L^\infty(X,\mu)
\rtimes_\sigma (\Gamma_1 * \Gamma_2)$ associated with a free
ergodic m.p. action $(\sigma, \Gamma_1 * \Gamma_2)$ on a
probability space $(X,\mu)$ can alternatively be viewed as an AFP
factor $M=M_1 *_A M_2$, where $A = L^{\infty}(X,\mu)$, $M_i = A
\rtimes_{\sigma_{|\Gamma_i}} \Gamma_i$, Theorem 0.1 allows us to
obtain Bass-Serre type vNE and OE rigidity results for actions of
free products of groups, as follows:

\proclaim{0.4. Theorem (vNE Bass-Serre rigidity)} Let $\Gamma_0,
\Lambda_0$ be groups with Haagerup property and $\Gamma_i,
\Lambda_j$, $1\leq i \leq n \leq \infty$, $1\leq j \leq m \leq
\infty$, be ICC groups having normal, non virtually abelian
subgroups with the relative property $(\text{\rm T})$.
Assume either $\Gamma_0$ is infinite or $n\geq 2$. Let
$\sigma$ $($resp. $\theta)$ be a free ergodic m.p. action of
$\Gamma=\Gamma_0 * \Gamma_1 * ...$ $($resp. $\Lambda=\Lambda_0 *
\Lambda_1 *...)$ on the probability space $(X, \mu)$ (resp.
$(Y,\nu)$) such that $\sigma_j=\sigma_{|\Gamma_j}$ $($resp.
$\theta_j=\theta_{|\Lambda_j})$ is ergodic $\forall j \geq 1$.
Denote $M=L^\infty(X,\mu)\rtimes_\sigma \Gamma$,
$N=L^\infty(Y,\nu)\rtimes_\theta \Lambda$, $M_i=L^\infty(X,\mu)
\rtimes_{\sigma_i} \Gamma_i \subset M$, $N_j =
L^\infty(Y,\nu)\rtimes_{\theta_j} \Lambda_j\subset N$ the
corresponding group measure space factors. If $\alpha: M\simeq
N^t$ is an isomorphism, for some $t > 0$, then $m=n$ and there
exists a permutation $\pi$ of indices $j \geq 1$ and unitaries
$u_j \in N^t$ such that $\text{\rm
Ad}(u_j)(\alpha(M_{j}))=N^t_{\pi(j)}$, $\text{\rm
Ad}(u_j)(\alpha(L^\infty(X,\mu)))=(L^\infty(Y,\nu))^t$, $\forall
j\geq 1$. In particular, $\Cal R_{\sigma} \simeq \Cal
R_{\theta}^t$ and $\Cal R_{\sigma_{j}} \simeq \Cal
R^t_{\theta_{\pi(j)}}$, $\forall j\geq 1$.
\endproclaim

In particular, taking the isomorphism $\alpha$ between the group
measure space factors in 0.4 to come from an orbit equivalence of
the actions, one gets:

\proclaim{0.5. Corollary (OE Bass-Serre rigidity)} Let $\Gamma_i,
\Lambda_j$, $0\leq i \leq n \leq \infty$, $1\leq j \leq m \leq
\infty$,  $\sigma$, $\theta$ be as in $0.4$. If $\Cal R_{\sigma,
\Gamma} \simeq \Cal R_{\theta, \Lambda}^t$ then $n=m$ and there
exists a permutation $\pi$ of the set of indices $i \geq 1$ such
that $\Cal R_{\sigma_i,\Gamma_i} \simeq \Cal R_{\theta_{\pi(i)},
\Lambda_{\pi(i)}}^t$, $\forall i\geq 1$.
\endproclaim

Like in ([P6II]), the terminology ``vNE rigidity'' is used here in
a broad sense, in the same spirit the terminology ``OE rigidity''
is being used in orbit equivalence ergodic theory ([Z1], [Fu1], [MS],
[S]). It can designate results which from an isomorphism of group
measure space factors derives orbit equivalence of the actions
involved (``vNE/OE rigidity'', like 6.2 in [P5]), or even
conjugacy of the actions (``vNE strong rigidity'', e.g. 7.1 in
[P6II]). Theorem 0.4 brings out a new type of vNE rigidity, which
we have labeled ``Bass-Serre'' because of its analogy to group
theory results. It is a ``vNE/OE''-type result but stronger, as it
derives not only the orbit equivalence of the ``main actions''
$(\sigma, \Gamma)$, $(\theta, \Lambda)$, but also the component by
component orbit equivalence of their restrictions $(\sigma_i,
\Gamma_i)$, $(\theta_i, \Lambda_i)$.

The ``vNE Bass-Serre rigidity'' can be used in combination with OE
rigidity results in orbit equivalence ergodic theory to get more
insight on the group measure space factors involved. Thus, taking
$\Gamma_0=\Lambda_0=\{1\}$ and $2 \leq n,m < \infty$ in 0.4, by
Gaboriau's results in ([G1]) it follows that the $\ell^2$-Betti
numbers of $\Gamma_i, \Lambda_j$ must satisfy
$\beta_k(\Gamma_i)=\beta_k(\Lambda_i)/t$, $\forall 1 \leq i \leq
n=m$, and $\Sigma_i \beta_k(\Lambda_i) + (n-1)= (\Sigma_i
\beta_k(\Lambda_i) +(n-1))/t$, forcing $t=1$.
Also, if we take $\Gamma=\Gamma_0 * \Gamma_1 *...$ and $\sigma$ as
in 0.4 and add the conditions Out$(\Cal R_{\sigma_1})=\{1\}$ and
$(\sigma_1, \Gamma_1)$ not OE to $(\sigma_i, \Gamma_i), \forall i
\neq 1,$ then Out$(\Cal R_\sigma)=\{1\}$ and Out$(M)=\text{\rm
H}^1(\sigma, \Gamma)$. Examples of actions $(\sigma_1, \Gamma_1)$
with the associated orbit equivalence relation $\Cal R_{\sigma_1,
\Gamma_1}$ having trivial outer automorphism group are constructed
in ([Ge2], [Fu2], [MS]) and we construct some more, using the
Monod-Shalom rigidity theorem ([MS]). The group H$^1(\sigma,
\Gamma)$ can in turn be calculated by using ([P8]), thus getting
explicit computations of Out$(M)$ for the group measure space
factors $M$. The fact that one can choose the action $(\sigma,
\Gamma)$ to be free yet have restrictions
$\sigma_{|\Gamma_i}$ isomorphic to specific
$\Gamma_i$-actions, $\forall i$,
is a consequence of ([To]) but we include
a proof for the reader's convenience (see 7.3 and A.1).

Theorem 0.1 is in fact used to obtain another (genuine) ``vNE/OE
rigidity'' result in this paper, for free ergodic m.p. actions
$(\sigma,\Gamma)$ with $\Gamma$ a free product of infinite groups,
$\Gamma=\Gamma_0
* \Gamma_1$, and $\sigma$ satisfying the relative property (T) of
(5.10 in [P5]), i.e. such that $L^\infty(X,\mu) \subset
L^\infty(X,\mu) \rtimes_\sigma \Gamma$ is a rigid inclusion. We
recover this way the uniqueness of the HT Cartan subalgebra in the
group-factors $L(\Bbb Z^2 \rtimes \Bbb F_n)$ and their
amplifications, one of the main results in ([P5]).

Similarly, we obtain rigidity results for crossed product factors
$M=R \rtimes_\sigma (\Gamma_0
* \Gamma_1)$ corresponding to actions $(\sigma, \Gamma_0 *
\Gamma_1)$ on the hyperfinite II$_1$ factor $R$, by regarding M as
an AFP factor $M=(R \rtimes \Gamma_0) *_R (R\rtimes \Gamma_1)$. In
fact, in this case we can control even better the groups of
symmetries $\mycal F(M)$, Out$(M)$, with complete calculations. To
state this result, let $f\Cal T_R$ denote the class of actions
$(\sigma, \Gamma_0*\Gamma_1)$ of free product groups $\Gamma_0 *
\Gamma_1$ on the hyperfinite factor $R$, satisfying the properties:
$(a)$. $\Gamma_0$ is free indecomposable; $(b)$. $\Gamma_1$ is
w-rigid (e.g. $\Gamma_1$ infinite Kazhdan group); $(c)$. $R \subset
R \rtimes_{\sigma} \Gamma_0$ is a rigid inclusion; $(d)$.
$\sigma_{|\Gamma_1}$ is a non-commutative Bernoulli $\Gamma_1$-action,
i.e. $R$ can be represented in the form $R =
\overline{\otimes}_{g\in \Gamma_1} (M_{n\times n} (\Bbb C), tr)_g$,
$n \geq 2$, with $\sigma_{|\Gamma_1}$ acting on it by left Bernoulli
shifts; $(e)$. $\sigma_{|\Gamma_1}$ is freely independent with
respect to the normalizer $\Cal N_0$ of $\sigma(\Gamma_0)$ in
Out$(R)$.

To show that such actions exist, we first prove that for any two
countable sets of automorphisms $S_1, S_2$ of $R$ there exist
$\theta \in \text{\rm Aut}(R)$ such that $S_1$ and $\theta S_2
\theta^{-1}$ are ``freely independent'' (see 8.2 and A.2). Combining
this with results from ([Ch], [P5], [NPS], [Va]) we deduce that for
many arithmetic groups $\Gamma_0$ (in particular for $\Gamma_0 =
SL(n,\Bbb Z), \forall n \geq 2$) and any w-rigid group $\Gamma_1$,
there exist actions $(\sigma,\Gamma_0 * \Gamma_1)$ on $R$ in the
class $f\Cal T_R$. Using Theorem 0.1, properties $(a)-(e)$ above and
[P3], we get:

\proclaim{0.6. Theorem} For any $\Gamma_0 = SL(n_0, \Bbb Z), n_0
\geq 2$ and any w-rigid group $\Gamma_1$ there exist actions
$\sigma$ of $\Gamma_0*\Gamma_1$ on $R$ in the class $f\Cal T_R$.
If $(\sigma, \Gamma_0 * \Gamma_1)$ is a $f\Cal T_R$ action and we
denote $M=R \rtimes_\sigma (\Gamma_0 * \Gamma_1)$, then $\mycal
F(M)=\{1\}$ and $\text{\rm Out}(M) = \text{\rm Char}(\Gamma_0)
\times \text{\rm Char}(\Gamma_1)$.
\endproclaim

\proclaim{0.7. Corollary} Given any compact abelian group $K$,
there exist separable $\text{\rm II}_1$ factors $M$ with $\mycal
F(M)=\{1\}$ and $\text{\rm Out}(M)=K$. For instance, if $(\sigma,
\Gamma_0*\Gamma_1)$ is a $f\Cal T_R$ action and $M = R
\rtimes_\sigma (\Gamma_0 * \Gamma_1)$ is the associated  crossed
product factor, with $\Gamma_0=SL(n,\Bbb Z)$, $\Gamma_1 =
SL(m,\Bbb Z) \times \hat{K}$ for some $n,m \geq 3$, then $\mycal
F(M)=\{1\}$ and $\text{\rm Out}(M) = K$. Moreover, denoting
$M^\infty = M \overline{\otimes} (\Cal B(\ell^2\Bbb N))$ the
associated $\text{\rm II}_\infty$ factor, we have $\text{\rm
Out}(M^\infty)=K$.
\endproclaim

The study of outer automorphisms of type II von Neumann factors
was at the core of Connes decomposition theory for factors of type
III and his classification of amenable factors, in the early 70's
([C1,2]). Two subsequent seminal papers ([C3,4]) gave the first
indications that the outer symmetry groups Out$(M)$ and $\mycal
F(M)$ can reflect rigidity properties of non-amenable factors. In
particular, it was shown in ([C3]) that $\mycal F(M)$ and Out$(M)$
are countable for group factors associated with ICC groups with
the property (T). The recent rigidity results in ([P5,6,9])
provide explicit calculations of $\mycal F(M)$ for large families
of group measure space factors $M$, and reduce the calculation of
Out$(M)$ to the computation of the commutants of the corresponding
group actions. However, such commutants are difficult to compute,
being left as an open problem even in the case of Bernoulli
actions (see [P6II]). The calculation of Out$(M)$ that we obtain
in this paper for crossed product factors arising from actions of
free products of w-rigid groups on the probability space and on
the hyperfinite factor thus give the first such explicit
computations.

This work initiated while the authors were visiting the
Laboratoire d'Alg\'ebres  d'Op\'erateurs at the Institut de
Math\'emathiques of the University Paris 7 during the Fall of
2004. It is a pleasure for us to thank the CNRS and the members of
the Lab for their support and kind hospitality. We are very
grateful to Damien Gaboriau for illuminating discussions and
comments on our work and Bass-Serre theory. We thank Dima
Shlyakhtenko for kindly pointing out to us Remark 6.6.

\heading 1. Conjugating subalgebras in AFP factors\endheading

\vskip .05in \noindent {\bf 1.1. AFP algebras}. Let $(M_1,
\tau_1), (M_2, \tau_2)$ be finite von Neumann algebras with a
common von Neumann subalgebra $B\subset M_i, i=1,2,$ such that
$\tau_{1|B}=\tau_{2|B}$. We denote by $(M_1 *_B M_2, \tau_1 *
\tau_2)$ the finite von Neumann algebra {\it free product with
amalgamation} (AFP) of $(M_1, \tau_1; B)$ and $(M_2, \tau_2; B)$,
as defined in ([V1] and pp. 384-385 of [P4]). Thus, $M_1 *_B M_2$
has a dense $*$-subalgebra
$$
B \oplus \bigoplus_{n\geq1} \bigoplus_{ {i_j \in \{1, 2\} } \atop
{i_1\neq i_2,\cdots,i_{n-1}\neq i_n}} \text{\rm sp}(M_{i_1}\ominus
B) (M_{i_2}\ominus B) \cdots (M_{i_n}\ominus B) \tag 1.1.1
$$
with the trace $\tau = \tau_1 * \tau_2$ defined on {\it reduced
words} by $\tau(x)= \tau_1(x)=\tau_2(x)$ for $x\in B$ and
$\tau(x)=0$ for $x=x_{i_1}x_{i_2} \ldots x_{i_n}$, with
$x_{i_k}\in M_{i_k} \ominus B$, $i_k \in \{1,2\}$, $i_1 \neq i_2,
\ldots, i_{n-1} \neq i_n$. Thus, the vector subspaces $B$ and
$\text{\rm sp}(M_{i_1}\ominus B) (M_{i_2}\ominus B) \cdots
(M_{i_n}\ominus B)\subset M$ in the above sum are all mutually
orthogonal with respect to the scalar product given by the trace
$\tau$. Also, their closure in $L^2(M,\tau)$ gives mutually
orthogonal Hilbert $B$-bimodules,
$$
L^2((M_{i_1}\ominus B) (M_{i_2}\ominus B) \cdots (M_{i_n}\ominus
B))\simeq \Cal H^0_{i_1} \otimes_B \Cal H^0_{i_2} \otimes_B ...
\otimes_B \Cal H_{i_n}^0,
$$
summing up to $L^2(M,\tau)$, where $\Cal H_i^0=L^2(M_i)\ominus
L^2(B)$.

\vskip .05in \noindent {\bf 1.2. Controlling intertwiners and
relative commutants}. In this sub-section we prove a very useful
``dichotomy-type'' result for subalgebras $Q$ of AFP factors
$M=M_1 *_B M_2$. It shows that if $Q$ sits in one of the factors,
say $M_1$, then it can either be conjugated into the ``core'' $B$
of the AFP algebra $M$, or else all its normalizer lies in $M_1$,
and even all ``intertwining'' Hilbert $Q-M_1$ bimodules $\Cal
H\subset L^2(M)$ with dim$(\Cal H_{M_1}) < \infty$ must be
entirely contained in $L^2(M_1)$ ! Also, in this second
alternative, any bimodule intertwining $Q$ into the
other factor, $M_2$, vanishes.

The first results of this type have been obtained in ([P1]), in
the case of ``plain'' free product factors, $M=M_1 * M_2$. The
next theorem provides a sharp generalization of those results. The
proof uses the basic ``intertwining criteria'' (2.1 and 2.3 in
[P6]), following arguments similar to (3.1 in [P6]).

\proclaim{1.2.1. Theorem} Let $(M_1, \tau_1)$ and $(M_2, \tau_2)$
be finite von Neumann algebras and $B$ a common von Neumann
subalgebra such that $\tau_{1|B} = \tau_{2|B}$.  Let $M = M_1 *_B
M_2$, $0\neq q \in \Cal P(M_1)$ and $Q \subset qM_1q$ a von
Neumann subalgebra. Assume no corner of $Q$ can be embedded into
$B$ inside $M_1$, i.e. $Q'\cap q\langle M_1, B \rangle q$ contains
no non-zero finite projections. If $0 \neq \xi \in L^2(qM)$
satisfies $Q\xi \subset L^2(\Sigma_{i}\xi_i M_k)$ for some $k \in
\{1,2\}$ and some $\xi_1,\dots ,\xi_n\in L^2(M)$, then $k=1$ and
$\xi \in L^2(M_1)$. In particular, $Q'\cap qMq\subset M_1$, the
normalizer $\Cal N_{qMq}(Q)$ of $Q$ in $qMq$ is contained in
$qM_1q$ and if $x\in M$ satisfies $Q x \subset x M_2$ then $x=0$.
\endproclaim
\vskip .05in \noindent {\it Proof.} Let $p$ denote the orthogonal
projection of $L^2(M)$ onto the Hilbert subspace $\overline{Q\xi
M_k}^{ \|.\|_2}\subset L^2(M)$. Note that $p\in Q'\cap q\langle M
,e_{M_k} \rangle q$ and $0\not= Tr(p)<\infty$, where
$Tr=Tr_{\langle M, e_{M_k} \rangle}$ denotes the canonical trace
on $\langle M, e_{M_k} \rangle$. To prove that $k=1$ and $\xi \in
M_k$ it is sufficient to show that $p \leq e_{M_k}$, equivalently
$(1-e_{M_k}) p (1-e_{M_k})=0$. Indeed, because then $Q\xi M_k
\subset L^2(M_k)$, so in particular $\xi\in L^2(M_k)$ and $u\xi
\in L^2(M_k)$, $\forall u\in \Cal U(Q)$. But since no corner of
$Q$ can be embedded into $B$ inside $M_1$, by (2.3 in [P6]) it
follows that for any $\varepsilon >0$ there exists $u\in \Cal
U(Q)$ such that $\|E_B(u)\|_2 \leq \varepsilon $. Thus, if $k=2$
then $\xi, u\xi \in L^2(M_2)$ so that
$u\xi=E_{M_2}(u\xi)=E_B(u)\xi$ and by the Cauchy-Schwartz
inequality we have
$$
\|\xi\|_1=\|u\xi\|_1 = \|E_B(u)\xi\|_1 \leq \|E_B(u)\|_2 \|\xi\|_2
\leq \varepsilon \|\xi\|_2.
$$
Since $\varepsilon > 0$ is arbitrary, this shows $\xi=0$. Thus,
the only possibility is that $k=1$, i.e. $\xi \in L^2(M_1)$.

By taking spectral projections, to show that $(1-e_{M_k}) p
(1-e_{M_k})=0$ it is in fact sufficient to show that if $f\in
Q'\cap \langle M, e_{M_k} \rangle$ is a projection such that $0
\neq Tr(f)<\infty$ and $f\leq 1-e_{M_k}$, then $f=0$. To this end,
we'll show that $\|f\|_{2,Tr}$ is arbitrarily small.

Thus, let $\eta_0=1,\eta_1, \ldots, \eta_n, \ldots \subset M$ be
an orthonormal basis of $M$  over $M_k$, i.e. $E_{M_k}(\eta_i^*
\eta_j) = \delta_{ij} p_j \in \Cal P(M_k)$ and $\| \eta \| <
\infty$, $\forall i,j$. If we denote $f_n=\Sigma_{i=1}^{n}
\eta_ie_{M_k}{\eta_i}^*$ then, as $f$ has finite trace and $f\leq
1-e_{M_k}=\Sigma_{i=1}^{\infty}\eta_ie_{M_k}{\eta_i}^*$, there
exists $n\in\Bbb N$ such that $\|f_nf-f \|_{2,Tr}<\epsilon \|f
\|_{2,Tr}$. Thus, if $u\in \Cal U(Q)$  then
$$Tr(f_nuf_nu^*)
\geq Tr(ff_nfuf_nu^*) \tag 1.2.1
$$
$$
-|Tr(ff_n(1-f)uf_nu^*)|-|Tr((1-f)f_n uf_nu^*)|.
$$

Using that $f_nf$ is $\epsilon-$close to $f$ in the norm $\| \cdot
\|_{2,Tr}$ and that $f$ commutes with $u\in Q$ we deduce:

$$
Tr(ff_nfuf_nu^*)= Tr(f_nfuf_nfu^*)\geq
(1-2\epsilon-{\epsilon}^2){\|f \|^2_{2,Tr}}. \tag 1.2.2
$$
Similarly we have:
$$ |Tr(ff_n(1-f)uf_nu^*)|+|Tr((1-f)f_nfuf_nu^*)|\leq 2\epsilon(1+\epsilon)
{\|f \|^2_{2,Tr}}. \tag 1.2.3
$$
Combining $(1.2.1)-(1.2.3)$ we get
$$
Tr(f_nuf_nu^*) \geq (1-4\epsilon-3{\epsilon}^2) {\|f \|^2_{2,Tr}},
\forall u\in \Cal U(Q). \tag 1.2.4
$$
On the other hand,
$$
Tr(f_nuf_nu^*)=
Tr(\Sigma_{i,j}\eta_ie_{M_k}{\eta_i}^*u{\eta_j}e_{M_k}{\eta_j}^*u^*)
=\Sigma_{i,j}||E_{M_k}(\eta_iu\eta_j^*)||_{2}^2. \tag 1.2.5
$$

Thus, in order to prove that $\|f\|_{2,Tr}$ is small, it is
sufficient to prove that $\forall \eta_0, ..., \eta_n \in M \ominus
M_k$, $\forall \varepsilon > 0$, $\exists u\in\Cal U(Q)$ such that
$||E_{M_k}(\eta_iu{\eta_j}^*)||_2\leq \epsilon$, $\forall 0\leq i,j
\leq n$. Furthermore, by $(1.1.1)$ and Kaplansky's density theorem
it is enough to prove this in the case $\eta_i$ are reduced words of
the form $\eta_i=\delta_ix_i$ such that one of the following holds
true: $(a)$. $\delta_i$ is a reduced word that ends with a letter in
$M_{2}\ominus B$ and $x_i$ is either equal to 1 or contained in
$M_1\ominus B$; $(b)$. $k=2$, $\delta_i=1$ and $x_i \in M_1\ominus
B$. Since $x_i u x_j^* \in M_1$, if we denote $y=x_i u x_j^*-E_B(x_i
u x_j^*)\in M_1 \ominus B$ then in both cases $(a)$ and $(b)$ the
reduced word $\delta_i y \delta_j^*$ is perpendicular to $M_k$.
Indeed, in case $(a)$, $\delta_i y \delta_j^*$ lies in
$...(M_{k'}\ominus B) (M_k \ominus B)(M_{k'}\ominus B)...$, where
$\{k,k'\}=\{1,2\}$, and thus it has length at least $3$, so
$\delta_i y \delta_j^*\perp M_k$ by $(1.1.1)$. In case $(b)$,
$\delta_i y \delta_j^* \in M_1\ominus B$ so it is perpendicular to
$M_2=M_k$. Using $\delta_i y \delta_j^*\perp M_k$ we then get:
$$
E_{M_k}(\eta_iu{\eta_j}^*)=E_{M_k}(\eta_iy{\eta_j}^*)+
E_{M_k}({\delta_i}E_B(x_iu{x_j}^*){\delta_j}^*)
$$
$$
=E_{M_k}({\delta_i}E_B(x_iu{x_j}^*){\delta_j}^*),
$$
implying that
$$
\|E_{M_k}(\eta_iu{\eta_j}^*)\|_2=
\|E_{M_k}({\delta_i}E_B(x_iu{x_j}^*){\delta_j}^*)\|_2
$$
$$
\leq||\delta_i|| ||\delta_j|| ||E_B(x_iu{x_j}^*)||_2.
$$
But by the hypothesis and (2.3 in [P6]), for any $\epsilon>0$ we
can find $u\in \Cal U(Q)$ such that $||E_B(x_iu{x_j}^*)||_2\leq
\epsilon/(||\delta_i|| ||\delta_j||)$. \hfill $\square$

\vskip .05in

Note that under the conditions of the above theorem, not only the
normalizer $\Cal N_{qMq}(Q)$ of $Q$ in $qMq$ follows contained in
$M_1$, but also the normalizer of the von Neumann algebra generated
by $\Cal N_{qMq}(Q)$, and so on. In fact, even the unitaries $u\in
qMq$ with the property that $uQu^* \cap qM_1q$ is not embeddable
into $B$ follow contained in $M_1$. More generally, if $Q \subset
qM_1q$ is a von Neumann subalgebra such that $Q'\cap q\langle M_1, B
\rangle q$ contains no non-zero finite projections, and if we denote
by $N_1=N(Q,M_1;B)$ the von Neumann subalgebra of $qMq$ generated by
unitaries $u\in qMq$ such that $(uQu^* \cap M_1)'\cap q\langle M_1,
B \rangle q$ contains no non-zero finite projections then $N_1
\subset M_1$. If we then repeat this operation, taking $N_2=N(N_1,
M_1;B)$ be the von Neumann algebra generated by all unitaries $u\in
qMq$ such that $uN_1u^* \cap q\langle M_1, B \rangle q$ contains no
non-zero finite projections, then $N_2 \subset M_1$. We can of
course continue this procedure inductively until it ``stops'', i.e.
until we reach an $N_\imath$ such that $N(N_\imath,
M_1;B)=N_\imath$. More, formally, consider the following:

\vskip .05in \noindent {\it 1.2.2. Definition}. If $q\in \Cal
P(B)$ and $Q\subset qMq$ we consider by (transfinite) induction
the strictly increasing family of von Neumann algebras $Q = N_0
\subset N_1 \subset ... \subset N_{\jmath} \subset ... \subset
N_{\imath}$, indexed by the first $\imath$ ordinals, such that:
$(a)$. For each $\jmath < \imath$, $N_{\jmath+1}=N(N_{\jmath},
M_1;B)$ and $N_\jmath \neq N_{\jmath +1}$; $(b)$. $N(N_\imath,
M_1; B) = N_{\imath}$; $(c)$. If $\jmath \leq \imath$ has no
``predecessor'' then $N_{\jmath} = \cup_{n < \jmath} N_n$. We then
denote $\tilde{N}(Q,M_1;B)=N_\imath$ and call it the {\it weak
quasi normalizer} (wq-normalizer) of $Q$ in $qMq$ relative to
$(M_1;B)$. Note that in fact both the definition of $N(Q, M_1;B)$
and $\tilde{N}(Q,M_1;B)$ make actually sense for any finite von
Neumann algebra $(M,\tau)$ and von Neumann subalgebras $B, M_1
\subset M$, $Q \subset qMq$, with $q\in \Cal P(M_1)$.

This definition is analogous to the definition of wq-normalizer of a
subgroup $H \subset G$ used in ([P5,6,8]). It is easy to see that
$\tilde{N}(Q,M_1;B)$ is the smallest von Neumann subalgebra $P$ of
$qMq$ such that $(uPu^*\cap qM_1q)' \cap q\langle M, e_B \rangle q$
contains non-zero finite projections $\forall u\in \Cal (qMq)
\setminus P$. Theorem 1.1 thus implies:

\proclaim {1.2.3. Corollary} With $(M_i,\tau_i)$, $q\in B \subset
M_i$, $i=1,2$, $M=M_1 *_B M_2$ as in $1.2.1$, let $Q\subset qM_1q$
be a von Neumann subalgebra such that no corner of $Q$ can be
embedded into $B$ inside $M_1$, i.e. $Q'\cap q\langle M_1, B
\rangle q$ contains no non-zero finite projections. Then
$\tilde{N}(Q,M_1;B) \subset M_1$.
\endproclaim

We will make repeated use of the following application of 1.2.1,
which shows that if one of the algebras $M_i$ involved in an
amalgamated free product $M=M_1*_B M_2$ contains a regular
subalgebra $Q$, then $Q$ must necessarily be contained in $B$,
modulo inner conjugacy:

\proclaim {1.2.4. Corollary} Let $(M_i,\tau_i)$, $B \subset M_i$,
$i=1,2$, be as in $1.1$ and denote $M=M_1 *_B M_2$. Let $Q\subset
qM_1q$ be a von Neumann subalgebra, for some $q \in \Cal P(B)$
with $qBq \neq qM_2q$. Assume $Q$ is regular in $qMq$. Then one
can embed a corner of $Q$ into $B$ inside $M_1$, i.e. $Q'\cap
q\langle M_1,e_B \rangle q$ contains non-zero finite-trace
projections.
\endproclaim
\noindent {\it Proof.} If $Q'\cap  q\langle M_1,e_B \rangle q$
contains no non-zero finite-trace projection then by Theorem 1.1
the normalizer $\Cal N_M(Q)$ of $Q$ in $qMq$ is contained in
$M_1$. Since $\Cal N_M(Q)''=qMq$, this implies $qMq=qM_1q$, thus
$qM_2q=qBq$, a contradiction. \hfill $\square$

\vskip .05in
\noindent

\vskip .05in \noindent {\bf 1.3. Locating subalgebras by means of
normalizers.} In this and the next sub-section we prove that if a
subalgebra of $M=M_1 *_B M_2$ is normalized by ``many'' unitaries
in $M_1$ then it must necessarily be contained in $M_1$. This
technical result will in fact not be needed until Section 7, where
it plays a key role in the proof of the Bass-Serre type Theorem
7.4. The proof uses the intertwining criteria in ([P6]) and a
careful asymptotic analysis of elements written in AFP expansion
$(1.1.1)$. We first prove the result assuming the subalgebra we
want to ``locate'' is a unitary conjugate of a subalgebra of $B$.
This assumption will be shown to be redundant in 1.4, in the case
$B=A$ is Cartan in $M$.

\proclaim{1.3.1. Proposition} Let $\Lambda_1, \Lambda_2$ be
discrete groups and $\sigma: \Lambda \rightarrow \text{\rm
Aut}(B,\tau)$ an action of $\Lambda= \Lambda_1 * \Lambda_2$ on the
finite von Neumann algebra $(B,\tau)$. Denote $M=B \rtimes_\sigma
\Lambda = M_1 *_B M_2$, where $M_i=B \rtimes_{\sigma_{|\Lambda_i}}
\Lambda_i$, $i=1,2$. Let $q\in \Cal P(B)$, $B_0 \subset qBq$ a von
Neumann subalgebra, $u\in \Cal U(qMq)$ and denote $\Cal N=\{v \in
\Cal U(qM_1q) \mid v(uqB_0qu^*)v^* =uqB_0qu^*\}$. Assume no corner
of $\Cal N''$ can be embedded into $B$ inside $M_1$. Then
$uqB_0qu^* \subset M_1$.
\endproclaim
\noindent {\it Proof.} Assume there exists $b_0\in qB_0q$,
$\|b_0\|\leq 1$ such that $d_0=ub_0u^*-E_{M_1}(ub_0u^*)$ satisfies
$c=\|d_0\|_2>0$. To get a contradiction, we first show that all
the unit ball of $uqBqu^*$ can be embedded into a set of the form
$\Sigma_{g\in F} (B)_1 u_g$ with $F \subset \Lambda$ finite,
$(B)_1$ denoting the unit ball of $B$. We need the following:

\proclaim{1.3.2. Lemma} Let $(B,\tau)$ be a finite von Neumann
algebra, $\sigma: \Lambda \rightarrow \text{\rm Aut}(B,\tau)$ an
action, $M=B \rtimes_\sigma G$ the corresponding crossed product
finite von Neumann algebra and $\{u_g\}_g \subset M$ the canonical
unitaries. For any finite set in the unit ball of $M$, $S_0
\subset (M)_1$, and any $\varepsilon
> 0$, there exists $F \subset \Lambda$ finite such that $x(B)_1y^*
\subset_\varepsilon \Sigma_{g \in F} (B)_1 u_g$, $\forall x,y\in
S_0$.
\endproclaim
\noindent {\it Proof}. By Kaplansky's density theorem there exists
a finite set $F_0\subset \Lambda$ and elements $\{b^x_g\in B \mid
x\in S_0, g\in F_0\},$ such that $x_0 = \Sigma_{g\in F_0} b^x_g
u_g$ satisfies $\|x_0\| \leq 1$ and $\|x-x_0\|_2 \leq
\varepsilon/2$, $\forall x\in S_0$. If we put $F=F_0F_0^{-1}$ then
we clearly have $x_0 B y_0^* \subset \Sigma_{g\in F} Bu_g$,
$\forall x,y\in S_0$. On the other hand, if $b\in B$ satisfies
$\|b\|\leq 1$ and we let $x_0by^*_0 = \Sigma_g b_gu_g$ then
$b_g=E_B((x_0by_0^*)u_g^* )$ and thus $\|b_g\| \leq
\|(x_0by_0^*)u_g^*\|  \leq 1$. This implies $\|xby^* - x_0
by_0^*\|_2 \leq \varepsilon$ and thus $xby^* \in_\varepsilon
\Sigma_{g \in F} (B)_1 u_g$. \hfill $\square$

\vskip .05in

By Lemma 1.3.2, it follows that there exists $F \subset \Lambda$
finite such that $uq(B)_1qu^* \subset_{\varepsilon/2} \Sigma_{g
\in F} (B)_1u_g$, where $\varepsilon = c^2/4$. Denote $N=\Cal
N''\subset qM_1q$.

For any $v\in \Cal N\subset qM_1q$ we then have $v(ub_0u^*)v^*
\in_{\varepsilon/2} \Sigma_{g \in F} (B)_1u_g$ as well. Since
$v(ub_0u^*)v^*= vd_0v^* + v(E_{M_1}(ub_0u^*))v^*$, with $vd_0v^*
\perp M_1$ and $v(E_{M_1}(ub_0u^*))v^*\in M_1$, by Pythagora it
follows that $vd_0v^* \in_{\varepsilon/2} \Sigma_{g \in F_0}
(B)_1u_g$, where $F_0=F\setminus \Lambda_1$. Letting $d_1 \in
\Sigma_{g \in F_0} (B)_1u_g$ be so that $\|d_0-d_1\|_2 \leq
\varepsilon/2$, we have thus shown:

\vskip .05in \noindent $(1.3.3)$. $\exists F_0\subset \Lambda
\setminus \Lambda_1$ finite and $d_1 \in \Sigma_{g\in F_0}(B)_1
u_g$ with $\|d_1\| \leq |F_0|$, $\|d_1\|_2 \geq c-c^2/8 >0$, such
that $vd_1v^* \in_\varepsilon \Sigma_{g\in F_0} (B)_1 u_g$ ,
$\forall v \in \Cal N$, where $\varepsilon =c^2/4$. \vskip .05in

Now note that by the condition satisfied by the algebra $\Cal
N''=N$, from (2.3 in [P6]) it follows that:

\vskip .05in \noindent $(1.3.4)$ $\forall K \subset \Lambda_1$
finite and $\delta > 0$, $\exists v \in \Cal N$ such that if $\xi$
denotes the projection of $v$ onto the Hilbert space $\oplus_{h
\in \Lambda_1 \setminus K} L^2(B)u_h$ then $\|v-\xi\|_2 \leq
\delta$. \vskip .05in

At this point, we need the following:

\proclaim{1.3.5. Lemma} With $(B, \tau)$, $(\sigma, \Lambda)$, $M,
\{u_g\}_g$ as in $1.3.2$, assume $\Lambda = \Lambda_1
* \Lambda_2$. Let $F_0 \subset \Lambda\setminus \Lambda_1$ be a finite
set. There exists $K = K(F_0)\subset \Lambda_1$ finite such that
any $\xi \in L^2(B \rtimes \Lambda_1)$ supported by $\Lambda_1
\setminus K$ satisfies $\xi (\Sigma_{g \in F_0} Bu_g) M_1 \perp
\Sigma_{g \in F_0} Bu_g $.
\endproclaim
\noindent {\it Proof.} Note that each irreducible alternating word
$g \in F_0$ has at least one letter from $\Lambda_2$. Let $K_0$
denote the set of elements in $\Lambda_1$ that can appear as first
letter in a word $g$ in $F_0$ (including the trivial letter $e$)
and denote $K=K_0K_0^{-1}$. Then $(\Lambda_1\setminus K)K_0 \cap
K_0 = \emptyset$. Now note that if $\xi \in L^2(B \rtimes
\Lambda_1)$ is supported by $\Lambda_1\setminus K$ then any
element $\eta$ in $\xi (\Sigma_{g \in F} Bu_g)$ is supported on
elements $g\in G$ that begin by a letter in $(\Lambda_1 \setminus
K) K_0 \subset \Lambda_1 \setminus K_0$. Moreover, this is still
the case for elements of the form $\eta x$ for $x \in M_1$. In
turn, any $g$ in the support of an element in $(\Sigma_{g\in F_0}
Bu_g)$ begins by a letter in $K_0$. Thus, the two vector spaces
$\xi (\Sigma_{g \in F_0} Bu_g)M_1$, $(\Sigma_{g \in F_0} Bu_g)$
are supported on disjoint subsets of $\Lambda_1
* \Lambda_2$ and are thus perpendicular. \hfill $\square$

\vskip .05in We continue the proof of Proposition 1.3.1 by letting
$K=K(F_0)$ be given by 1.3.5, for the finite set $F_0 \subset
\Lambda\setminus \Lambda_1$ from $(1.3.3)$. Let
$\delta=\varepsilon/|F_0|$ and choose $\xi\in L^2(B \rtimes
\Lambda_1)$ supported on $\Lambda_1 \setminus K$, as given by
$(1.3.4)$. Then
$$
\|\xi d_1 v^*-vd_1v^*\|_2 \leq \|\xi-v\|_2 \|d_1\| \leq \delta
|F_0| \leq \varepsilon,
$$
which together with $(1.3.3)$ implies $\xi d_1 v^*
\in_{2\varepsilon} \Sigma_{g\in F_0} B u_g$. But by Lemma $1.3.5$
we have $\xi d_1 v^* \perp \Sigma_{g\in F_0} B u_g$. Thus $\|\xi
d_1 v^*\|_2 \leq 2\varepsilon$. On the other hand
$$
\|\xi d_1 v^*\|_2 = \|\xi d_1 \|_2 \geq \|d_1\|_2 - \|\xi-v\|_2
\|d_1\| \geq c-c^2/8  -\varepsilon > 2\varepsilon,
$$
a contradiction which ends the proof of Proposition 1.3.1. \hfill
$\square$

\vskip .05in \noindent {\bf 1.4. A Cartan conjugacy result.} We
now prove that if the ``core'' $B$ of an AFP algebra $M=M_1*_B
M_2$ is maximal abelian and regular (thus Cartan) in $M$, then any
other Cartan subalgebra $A_0\subset M$ which is normalized by
``many'' unitaries in $M_1$ is a unitary conjugate of $B=A$. Note
that it strenghtens both Corollary 1.2.4 and Proposition 1.3.1, in
the case the core $B=A$ is abelian and Cartan in $M$.

\proclaim{1.4.1 Theorem} Let $\Lambda_1, \Lambda_2$ be infinite
discrete groups and $\sigma: \Lambda \rightarrow \text{\rm
Aut}(A,\tau)$ a free ergodic action of $\Lambda= \Lambda_1 * \Lambda_2$ on a
diffuse abelian von Neumann algebra $(A,\tau)$. Denote $M=A
\rtimes_\sigma \Lambda = M_1 *_A M_2$, where $M_i=A
\rtimes_{\sigma_{|\Lambda_i}} \Lambda_i$, $i=1,2$. Let $q\in \Cal
P(A)$, $A_0 \subset qMq$ a Cartan subalgebra such that no corner
of $(\Cal N_{qMq}(A_0) \cap qM_1q)''$ can be embedded into $A$
inside $M_1$. Then $A_0 \subset qM_1q$ and there exists $u\in\Cal
U(qM_1q)$ such that $uA_0u^*=Aq$.
\endproclaim

\proclaim{1.4.2. Lemma} Let $\Lambda_1,\Lambda_2$ be discrete
groups and let $\sigma :\Lambda=\Lambda_1*\Lambda_2\rightarrow
\text{\rm Aut}(B,\tau)$ be a trace preserving action  on a finite
von Neumann algebra $(B,\tau)$. Denote $M=B \rtimes_\sigma \Lambda
= M_1 *_B M_2$ and let $M_1 = B \rtimes_{\sigma_{|\Lambda_1}}
\Lambda_1\subset M$. Let $q\in \Cal P(B)$ and assume that
$A_0\subset qMq$ is a diffuse abelian von Neumann
subalgebra such that no corner
of $A_0$ can be embedded into $qM_1q$ inside $M$. Then for any
$\varepsilon>0$ there exist $F\subset \Lambda\setminus \Lambda_1$
finite and $x_1,x_2\in\Sigma_{g\in F} B{u_g}$ such that any $u\in
\Cal N_{qMq}(A_0)$ satisfies
$$
\|ux_1u^*x_2-x_2ux_1u^*\|_2\leq \varepsilon,\sqrt{\tau(q)}
-\varepsilon\leq \|ux_1u^*x_2\|_2\leq \sqrt{\tau(q)}+\varepsilon.
$$
\endproclaim
\noindent{\it Proof.} By the assumption on $A_0$ it follows from
(2.3 in [Po6]) that there exists $a_1\in\Cal U(A_0)$ such that
$\|E_{M_1}(a_1)\|_2<\varepsilon/4$. Thus, we can find $F_1\subset
\Lambda\setminus \Lambda_1$ finite and $x_1\in\Sigma_{g\in F_1}
B{u_g}$ such that $\|a_1-x_1\|_2\leq\varepsilon/4$. Repeating the
above argument we can now find $a_2\in\Cal U(A_0)$, a finite set
$F$ with $F_1\subset F\subset \Lambda \setminus\Lambda_1$ and
$x_2\in\Sigma_{g\in F}B{u_g}$  such that
$\|a_2-x_2\|_2\leq\varepsilon/(4\|x_1\|)$. Using these
inequalities we get for $u\in \Cal N_{qMq}(A_0)$ (in fact for all
$u\in M$ with $\|u\|\leq 1$):
$$
\|ux_1u^*x_2-ua_1u^*a_2\|_2\leq \|ux_1u^*(x_2-a_2)\|_2 +
\|u(x_1-a_1)u^*a_2\|_2
$$
$$
\leq\|x_1\|\|x_2-a_2\|_2+\|x_1-a_1\|_2 \leq
\|x_1\|\varepsilon/(4\|x_1\|)+\varepsilon/4 =\varepsilon/2,\forall
u\in(M)_1.
$$
Similarly, it follows that $\|x_2ux_1u^*-a_2ua_1u^*\|_2\leq
\varepsilon/2,\forall u\in (M)_1$.  Finally, if $u\in\Cal
N_{qMq}(A_0)$ then $ua_1u^*a_2=a_2ua_1u^*$ (because $A_0$ is
abelian) and $\|ua_1u^*a_2\|_2=\sqrt{\tau(q)}$, and thus by
combining this with the above inequalities we get the desired
estimates. \hfill $\square$

\proclaim{1.4.3. Lemma} With the same notations and assumptions as
in $1.4.2$, let $F\subset\Lambda \setminus\Lambda_1$ be a finite
set and let $x_1,x_2\in\Sigma_{g\in F}Bu_g$. Then $\forall
\varepsilon>0,\exists K=K(F,\varepsilon)\subset \Lambda_1$ finite
and $\delta=\delta(F,\varepsilon)>0$ such that if  $u\in\Cal
(M_1)_1$ satisfies $\|E_{B}(uu_g^*)\|_2\leq \delta,\forall g\in
K$, then it must also satisfy
$ux_1u^*x_2\perp_{\varepsilon}x_2ux_1u^*$.
\endproclaim

\noindent{\it Proof.} Since $F\subset \Lambda\setminus\Lambda_1$
is finite, by the free decomposition $\Lambda=\Lambda_1 *
\Lambda_2$ one readily deduces that there exists $K=K^{-1}\subset
\Lambda_1$ finite such that $ (\Lambda_1\setminus
K)F(\Lambda_1\setminus K)F$
has empty intersection with  $F(\Lambda_1\setminus
K)F(\Lambda_1\setminus K).$ Next, let $u\in (M_1)_1$ and
set $u'=\Sigma_{g\in K}E_{B}(u{u_g}^*)u_g$ and $u''=u-u'$. Then
$u''$ is supported on $\Lambda_1\setminus K$ and we have the
decomposition

$$
ux_1u^*x_2=u''x_1{u''}^*x_2+u'x_1u^*x_2+ux_1{u'}^*x_2-u'x_1{u'}^*x_2.
$$

Denote $x_{12}=u''x_1{u''}^*x_2$ . Then, $x_{12}$ is supported on
$(\Lambda_1\setminus K)F(\Lambda_1\setminus K)F$ and we have the
following estimate:
$$\|ux_1u^*x_2-x_{12}\|_2\leq (2\|u\|+\|u'\|)\|x_1\|\|x_2\|\|u'\|_2
$$
$$
\leq (2+|K|)\|x_1\|\|x_2\|\|u'\|_2
$$

Similarly, if we denote $x_{21}=x_2u''x_1{u''}^*$ , then $x_{21}$
is supported on $F(\Lambda_1\setminus K)F(\Lambda_1\setminus K)$
and $\|x_2ux_1u^*-x_{21}\|_2\leq (2+ |K|)\|x_1\| \|x_2\|
\|u'\|_2$.

Next, we show that $K$ and $\delta= \varepsilon (12|K|(\|x_1\|
\|x_2\|+1))^{-3/2}$ satisfy the conclusion. To this end, let $u\in
(M_1)_1$ be such that $\|E_{B}(uu_g^*)\|_2\leq \delta,\forall g\in
K.$ Then
$$
\|u'\|_2=(\Sigma_{g\in K}\|E_{B}(u{u_g}^*)\|^2_2)^{1/2} \leq
\varepsilon |K|^{-1} (12(\|x_1\| \|x_2\|+1))^{-3/2},
$$
hence $\|ux_1u^*x_2-x_{12}\|_2\leq \varepsilon (\|x_1\|
\|x_2\|+1)^{-1}/4$ and $\|x_2ux_1u^*-x_{21}\|_2\leq\varepsilon
(\|x_1\| (\|x_2\|+1)^{-1}/4$. Also, we have:
$$
\|x_{12}\|_2\leq \|x_{12}-ux_1u^*x_2\|_2+\|ux_1u^*x_2\|_2
$$
$$
\leq\varepsilon (\|x_1\| \|x_2\|+1)^{-1}/4+\|x_1\|\|x_2\|\leq
\|x_1\|\|x_2\|+1.
$$
But by the way we have chosen $K$,  $x_{12}$ and $x_{21}$ have
disjoint supports, hence $x_{12}\perp x_{21}$. Thus
$$
|\langle ux_1u^*x_2,x_2ux_1u^* \rangle|
$$
$$ \leq |\langle
ux_1u^*x_2-x_{12},x_2ux_1u^* \rangle|+|\langle
x_{12},x_2ux_1u^*-x_{21} \rangle|
$$
$$
\leq \|ux_1u^*x_2-x_{12}\|_2\|x_2ux_1u^*\|_2+
\|x_{12}\|_2\|x_2ux_1u^*-x_{21}\|_2
$$
$$
\leq \varepsilon (\|x_1\| \|x_2\|+1)^{-1} \|x_1\|\|x_2\|/4+
(\|x_1\|\|x_2\|+1) (\varepsilon (\|x_1\| \|x_2\|+1)^{-1}/4)
<\varepsilon.
$$
\hfill $\square$

\proclaim{ 1.4.4. Proposition} With the same notations and
assumptions as in $1.4.2$ and $1.4.3$, let $q\in \Cal P(B)$ and
let $A_0 \subset qMq$ be a diffuse abelian von Neumann subalgebra.
Assume no corner of $(\Cal N_{qMq}(A_0)\cap qM_1q)'' \subset qMq$
can be embedded into $B$ inside $M_1$. Then a corner of $A_0$ can
be embedded into $qM_1q$ inside $qMq$.
\endproclaim
\noindent{\it Proof.} Assume that no corner of $A_0$ embeds into
$M_1$.  Apply first 1.4.2 for $\varepsilon=\tau(q)^{1/2}/4$ to
deduce that there exists $F\subset \Lambda\setminus \Lambda_1$
finite and $x_1,x_2\in M$ supported on $F$ such that if $u\in\Cal
N_{qMq}(A_0)$ then $\|ux_1u^*x_2-x_2ux_1u^*\|_2\leq \tau(q)^{1/2}/4$
and $3 \tau(q)^{1/2}/4\leq\|ux_1u^*x_2\|_2\leq 5 \tau(q)^{1/2}/4.$
It then follows that $|\langle ux_1u^*x_2,x_2ux_1u^* \rangle|\geq
\tau(q)/4,\forall u\in\Cal N_{qMq}(A).$

By Lemma 1.4.3, there exist $K\subset \Lambda_1$ finite and
$\delta>0$ such that if $u\in (M_1)_1$ satisfies
$\|E_{B}(uu_g^*)\|_2\leq \delta,\forall g\in K$, then
$ux_1u^*x_2\perp_{\tau(q)/4}x_2ux_1u^*$. But this implies that we
cannot find $u\in \Cal N=\Cal N_{qMq}(A_0)\cap qM_1q$ such that
$\|E_{B}(uu_g^*)\|_2\leq \delta,\forall g\in K$. By (2.3 in [P6])
this contradicts the fact that no corner of $\Cal N''$ embeds into
$B$ inside $M_1$. \hfill $\square$

The next result provides a useful ``transitivity'' property
for the subordination relation considered in Corollary 1.2.4:

\proclaim{ 1.4.5. Lemma} Let $M$ be a finite von Neumann
algebra, $B_0, M_1\subset M$ von Neumann  subalgebras of $M$
and $Q\subset M_1$ a von Neumann subalgebras of $M_1$. Assume
there exist projections $q_0 \in B_0$, $q_1 \in M_1$,
a unital isomorphism of $q_0B_0q_0$ into
$q_1M_1q_1$ and a partial isometry $v \in M$ such that $v^*v
\in (q_0B_0q_0)'\cap q_0Mq_0$,
$vv^* \in \psi(q_0B_0q_0)'
\cap q_1 M q_1$ and $vb = \psi(b)v, \forall b
\in q_0B_0q_0$.
Denote by $q'$ the
support projection of
$E_{M_1}(vv^*)\in \psi(q_0B_0q_0)'\cap q_1M_1q_1$
and $B_1=\psi(q_0B_0q_0)q'$.
If a corner of $B_1=\psi(q_0B_0q_0)q'$ can
be embedded into $Q$ inside $M_1$, then a corner of $B_0$ can be embedded into $Q$ inside $M$.
\endproclaim
\noindent{\it Proof.} Indeed, if $p_1\in
\Cal P(B_1)$, $v_1\in M_1p_1$ is a
non-zero partial isometry and $\psi_1: p_1B_1p_1\rightarrow
Q$ a (not necessarily unital) isomorphism such
that $v_1b = \psi_1(b)v_1$, $\forall b\in p_1B_1p_1$,
then $v_1v\neq 0$ and $v_1v b = \psi_1(\psi(b))v_1v$,
$\forall b\in q_0B_0q_0$.
\hfill $\square$

\vskip .05in \noindent {\it Proof of Theorem 1.4.1}. By Proposition
1.4.4 and (2.1 in [P6]), there exist projections $q_0 \in A_0
\subset qMq$, $q_1 \in qM_1q$, a unital isomorphism of $A_0q_0$ into
$q_1M_1q_1$ and a partial isometry $v \in M$ such that $v^*v = q_0,
vv^* \in \psi(A_0q_0)'\cap q_1 M q_1$ and $va = \psi(a)v, \forall a
\in A_0q_0$. Let $q'$ be the
support projection of $E_{M_1}(vv^*)$ and note that if we denote
$A_1 = \psi(A_0q_0) \subset q_1M_1q_1$ then $q'\in A_1'\cap
q_1M_1q_1$. By replacing if necessary $\psi$ by $q'\psi(\cdot)q'$
and shrinking $q_0\in A_0$ accordingly, we may assume $q_1=q'$.

Now, if a corner of
$A_1$ can be embedded into $A$ inside $M_1$,
then by Lemma 1.4.5 a corner of
$A_0$ can embedded into $A$ inside $M$, so by (A.1 in [P5]) the two
Cartan subalgebras $A_0, Aq \subset qMq$ are unitary conjugate. If
in turn no corner of $A_1$ can be embedded into $A$ inside $M_1$
then by Theorem 1.2.1 we have $vv^* \in A_1'\cap q_1Mq_1 \subset
q_1M_1q_1$, implying that $vA_0v^* \subset q_1M_1q_1$. By
spatiality, $vA_0v^*$ is Cartan in $q_1Mq_1$, which by
Corollary 1.2.4 implies
that a corner of $vA_0v^*$ can be embedded into $A$ inside $M_1$. By
(A.1 in [P5]), this implies $A_0$ and $Aq$ are unitary conjugate in
$qMq$. On the other hand, by Proposition 1.3.1 we have $A_0\subset
qM_1q$. But then Theorem 1.2.1 applies to show that $A_0$, $Aq$ are
conjugate in $qM_1q$ as well. \hfill $\square$

\heading 2. Deformation of AFP factors.
\endheading

Throughout this section,  $(M_1, \tau_1)$, $(M_2, \tau_2)$ will be
finite von Neumann algebras with $B\subset M_i$ a common von
Neumann subalgebra such that $\tau_{1|B} = \tau_{2|B}$, as in
Section 1. We describe in this section several useful ways to
deform the identity map on the AFP algebra $M=M_1 *_B M_2$ by
subunital, subtracial, c.p. maps which arise naturally from the
amalgamated free product structure of $M$. By a {\it deformation}
of $id_M$ we mean here a sequence $\phi_n$ of subunital,
subtracial, c.p. maps on $M$ such that $\lim_n \|\phi_n(x)-x\|_2 =
0$, $\forall x\in M$.

\vskip .1in \noindent {\bf 2.1. Amalgamated free product of c.p.
maps}. We first recall the definition of amalgamated product of
c.p. maps from ([Bo]), and establish some basic properties.

\proclaim{2.1.1. Lemma} Let $\phi_i:M_i \rightarrow M_i$ be
subunital, subtracial, c.p. maps with $\phi_1 (b) = \phi_2 (b)$,
$\forall b \in B$, and $\tau \circ \phi_i = \lambda \tau$,
$\forall i = 1, 2$, with $0 < \lambda \leq 1$. Then $\phi_1 *_B
\phi_2: M_1 *_B M_2 \rightarrow M_1 *_B M_2$ is a well-defined
subunital, subtracial, c.p. map. Moreover, the map $(\phi_1,
\phi_2) \mapsto \phi_1 *_B \phi_2$ is continuous with respect to
the topologies given by pointwise $\| \cdot \|_2$-convergence.
\endproclaim

\vskip .05in
\noindent
{\it Proof}. Since $\tau \circ \phi_i = \lambda \tau$, $\forall i = 1,
2$
we have that
$\tau \circ ( \phi_1 *_B^{\text{alg}} \phi_2 ) = \lambda \tau$, and
hence by [Bo] extends uniquely to a c.p. map
$\phi_1 *_B \phi_2$ on $M_1 *_B M_2$ which follows subtracial.

To show that this correspondence is continuous suppose that
$\varepsilon > 0$ and
$x_1, \ldots, x_n \in M_1 *_B M_2$.  As $M_1 *_B^{\text {alg}} M_2$ is
dense in
$M_1 *_B M_2$ let $x_1', \ldots, x_n' \in M_1 *_B^{\text {alg}} M_2$
such that
$\| x_j - x_j' \|_2 \leq \varepsilon / 3$, $\forall j \leq n$.
Hence $\exists m,l \in \Bbb N$, $F_i \subset M_i$ finite, such that
each
$x_i'$ is the sum of at most $m$ products of lenght at most $l$ from
$F_1 \cup F_2$.
Let $N = \max_{x \in F_1 \cup F_2} \| x \|$, then if $\phi_i': M_i
\rightarrow M_i$
are subunital, subtracial, c.p. maps with
$\| \phi_i'(x) - \phi_i(x) \|_2 < \varepsilon / 3nlN^l$,
$\forall x \in F_i$, $i = 1,2$,
then repeated use of the
triangle inequality together with the fact that subunital c.p.
maps are contractions
in the uniform norm
show that
$\| \phi_1' *_B \phi_2' (x_i) - \phi_1 *_B \phi_2 (x_i) \|_2 <
\varepsilon$,
$\forall i \leq n$.
\hfill $\square$

\vskip .1in \noindent {\bf 2.1.2. Remark}. In general, the free
product of two subunital subtracial c.p. maps need not be
subtracial. In fact, given any c.p. map $\phi_1$ which is unital,
subtracial, but not tracial, the c.p. map $\phi = \phi_1 * id$ is
not subtracial. Even more so, the Radon-Nykodim derivative $d\tau
\circ \phi / d\tau$ of any such free product c.p. map $\phi$ is
unbounded. To see this, let $x \in M_1$ such that $\tau(x) = 0$,
but $\tau \circ \phi_1 (x) \not= 0$, let $v \in M_2$ be a partial
isometry with $\tau(v) = 0$ and set $p = vv^*$. Then we have:

$$
\tau \circ \phi ((vxv^*)^*(vxv^*)) = |\tau \circ \phi_1(x)|^2
\tau(p) + \| v (\phi_1(x) - \tau \circ \phi_1(x)) v^* \|_2^2 $$
$$
= |\tau \circ \phi_1(x)|^2 \tau(p) + \tau (p)^2 \| \phi_1(x) -
\tau \circ \phi_1(x) \|_2^2. $$ Thus, if we choose $v$ such that
$\tau(p) \rightarrow 0$ then
$$
\tau \circ \phi ((vxv^*)^*(vxv^*)) / \tau ((vxv^*)^*(vxv^*))
\rightarrow \infty.
$$

We note that in [Bo] it was assumed that a free product of
subtracial c.p. maps is subtracial in order to show that the
Haagerup property is preserved by free products. But although the
above argument shows that this fact doesn't hold true unless the
c.p. maps are actually tracial, the result on the Haagerup
property is still valid, since by (Proposition 2.2. in [Jol]) one
can take the c.p. maps given by the Haagerup property to be unital
and tracial.

\vskip .1in \noindent {\bf 2.2. Deformation by automorphisms}.
Suppose $\alpha_i \in \text{\rm Aut}(M_i, \tau_{M_i})$
such that $\alpha_1 (b)
= \alpha_2 (b)$, $\forall b \in B$. Then since automorphisms are
unital c.p. maps we have that $\alpha = \alpha_1 *_B
\alpha_2$ is a unital, tracial, c.p. map, moreover $\alpha$
restricted to the dense subalgebra $M_1 *_B^{\text {alg}} M_2$ is
an automorphism and so by continuity we have that $\alpha$ is an
automorphism.

Hence if $\alpha_i^t \in \text{\rm Aut}(M_i)$, $t \in \Bbb R$ is a
one parameter group of automorphisms of $M_i$ which is pointwise
$\| \cdot \|_2$-continuous and satisfies $\alpha^t_{i|B} = id_B$,
for $i=1,2$, then $\alpha^t$ gives a deformation of the identity
of $M$ by automorphisms. In particular, we have:

\proclaim{2.2.1. Lemma} Let $v_i \in \Cal U (B' \cap M_i)$,
$i=1,2$. Then there is a pointwise $\| \cdot \|_2$-continuous one
parameter group of automorphisms $\{ \alpha^t \}_{ t \in \Bbb R }
\subset {\text{\rm Aut}}(M)$ such that $\alpha_1 = {\text{\rm
Ad}}(v_1) *_B {\text{\rm Ad}}(v_2)$.
\endproclaim
\noindent {\it Proof}. Let $h_j=h_j^*\in B'\cap M_j$ be such that
$\text{\rm exp}(\pi i h_j)=v_j$ and define $\alpha^t_j = \text{\rm
Ad}(\text{\rm exp}(\pi t i h_j))$, $t \in \Bbb R$, $j=1,2$, and
the above observation applies. \hfill $\square$

\vskip .1in For the next lemma we denote $\tilde M = M *_B
(B\overline{\otimes} L(\Bbb F_2))$. Note that if we let $L(\Bbb
F_2) = L(\Bbb Z * \Bbb Z)=L(\Bbb Z) * L(\Bbb Z)$ and denote
$\tilde M_i =M_i *_B (B\overline{\otimes} L(\Bbb Z))$, $i=1,2$,
then $\tilde M= \tilde M_1 *_B \tilde M_2$. Also, if $u_1 \in
L(\Bbb Z
* 1) \subset L(\Bbb F_2)$, $u_2 \in L(1 * \Bbb Z) \subset L(\Bbb F_2)$
are the canonical generating unitaries then $u_i \in B'\cap \tilde
M_i$, $i=1,2$. We will use the algebra $\tilde{M}$ as framework
for the main deformation of $M$, 2.2.2 below. The
action of $M$ on the $(1.1.1)$-decomposition of the AFP algebra
$\tilde{M}=\tilde{M}_1 *_B \tilde{M}_2$ can be viewed as the
analogue of the action of an amalgamated free product
group $\Lambda_1 *_H \Lambda_2$
on the Bass-Serre tree with
vertices $\Lambda_1/H \cup \Lambda_2/H$.

In turn, the ``graded deformation'' below is inspired by the
graded deformations of crossed product algebras involving
Bernoulli actions in ([P3], [P6]).

\proclaim{2.2.2. Lemma} There exists a pointwise $\| \cdot
\|_2$-continuous one parameter group of automorphisms $\{ \theta_t
\}_{t\in \Bbb R}$ and a period 2 automorphism $\beta$ of $\tilde
M$ such that:

\noindent $(a)$. $\theta_0=id, \theta_1 = {\text{\rm Ad}}(u_1) *_B
{\text{\rm Ad}}(u_2)$.

\noindent $(b)$. $\beta\theta_t\beta=\theta_{-t},\forall t\in\Bbb
R$.

\noindent $(c)$. $M\subset{\tilde M}^{\beta}$.
\endproclaim
\noindent {\it Proof}. Let $A_j$ be the von Neumann subalgebra
generated by $u_j$, and let $h_j \in A_j$ be self-adjoint elements
with spectrum in $[-\pi, \pi]$ such that $u_j = \text{exp}(\pi i
h_j)$. Set $u_j^t = \text{exp}(\pi i t h_j)$, $j = 1,2$, $t \in \Bbb
R$. Then $\theta_t = {\text{\rm Ad}}(u_1^t) *_B {\text{\rm
Ad}}(u_2^t) \in \text{\rm Aut}(\tilde M)$, $\forall t \in \Bbb R$,
defines a pointwise $\| \cdot \|_2$-continuous one parameter group
of automorphisms which satisfies $(a)$.

Let $\beta$ be the unique automorphism of $\tilde M$ satisfying
$\beta_{\mid M}=id_M$ and $\beta(u_j)=u_j^*,j=1,2$. Then $\beta$
is clearly a period 2 automorphism and it satisfies $(c)$ by
definition. Also, for $x \in M=M_1 *_B M_2$ we have

$$
\beta\theta_t\beta(x)=\beta\theta_t(x)=\beta({\text{\rm Ad}}(\exp(
\pi ith_1))* {\text{\rm Ad}}(\exp (\pi ith_2 )))(x)
$$
$$
=({\text{\rm Ad}}(\exp (-\pi ith_1))*{\text{\rm Ad}}(\exp (-\pi
ith_2)))(x)=\theta_{-t}(x),\forall x\in M.
$$
Similarly, for $u_1, u_2$  we have
$$
\beta\theta_t\beta(u_j)=\beta\theta_t(u_j^*)=
\beta(u_j^*)=u_j=\theta_{-t}(u_j).
$$
Since $u_1$, $u_2$, and $M$ generate $\tilde M$
as a von Neumann algebra it follows that
$\beta\theta_t\beta=\theta_{-t},\forall t$. \hfill $\square$

\vskip .1in \noindent {\bf 2.3. Deformation by free products of
multiples of the identity}. Recall that if $\Cal H_i^0 = L^2(M_i)
\ominus L^2(B)$ then we may decompose $L^2(M_1 *_B M_2)$ in the
usual way as
$$
L^2(M_1 *_B M_2) = L^2(B) \oplus
\bigoplus_{n\geq1} \bigoplus_{
{i_j \in \{1, 2\} } \atop
{i_1\neq i_2,\cdots,i_{n-1}\neq i_n}}
\Cal H_{i_1}^0 \otimes_B \Cal H_{i_2}^0 \otimes_B \cdots \otimes_B \Cal
H_{i_n}^0.
$$

For each $L \in \Bbb N$ we let $\widehat E_L$ be
the projection onto the subspace
$$
\bigoplus_{n\geq L} \bigoplus_{
{i_j \in \{1, 2\} } \atop
{i_1\neq i_2,\cdots,i_{n-1}\neq i_n}}
\Cal H_{i_1}^0 \otimes_B \Cal H_{i_2}^0 \otimes_B \cdots \otimes_B \Cal
H_{i_n}^0.
$$

Let $\{ c_n \}_n \subset [0,1)$ such that $c_n \nearrow 1$, then
by Lemma 1.1.1 we have:

\proclaim{2.3.1. Lemma}  $\phi_n = (c_n $id$) *_B (c_n $id$)$
gives a deformation of the identity of $M$.  Moreover, $\phi_n$
commutes with $\widehat E_L$ as operators on $L^2(M)$ and $\|
\phi_n \circ \widehat E_L (x) \|_2 \leq c_n^L \| x \|_2$,
$\forall n, L \in \Bbb N$, $x \in M$.
\endproclaim
\noindent {\it Proof}. Trivial by the definitions. \hfill

\vskip .1in \noindent {\bf 2.4. Deformation by subalgebras}. For
each $i \in \{ 1,2 \}$ let $N_i^j \subset M_i$ be an increasing
sequence of von Neumann subalgebras such that $B \subset N_i^1$,
and $\overline { \cup_{j \geq 1} } N_i^j = M_i$.  Let $E_i^j: M_i
\rightarrow M_i$ be the conditional expectation onto $N_i^j$. Then
$E^j = E_1^j *_B E_2^j$ gives a sequence of conditional
expectations of $M$ onto $N_1^j *_B N_2^j$, which by 2.1.1
converges to the identity pointwise, i.e. $\overline { \cup_{j
\geq 1} } N_1^j *_B N_2^j = M$.

A particular case of such a deformation, which works whenever $B'
\cap M_i$ is diffuse, is given by $N_i^j = p_i^j M_i p_i^j \oplus
B(1-p_i^j)$, with $p_i^j \in \Cal P(B' \cap M_i)$ satisfying
$p_i^j \nearrow 1_{M_i}$, $i = 1,2$. Indeed, this clearly implies
$\overline { \cup_{j \geq 1} } ( p_1^j M_1 p_1^j \oplus B(1-p_1^j)
) *_B ( p_2^j M_2 p_2^j \oplus B(1-p_2^j) ) = M$.

\heading 3. Deformation/rigidity arguments.  \endheading

In this section we investigate the effect that the deformations
considered in Section 2 have on the relatively rigid subalgebras
of $M_1 *_B M_2$ . To this end, first recall from (Sec. 4 in [P5])
that if $Q \subset M$ is a von Neumann subalgebra of the finite
von Neumann algebra $(M,\tau)$, then $Q \subset M$ is called a
{\it rigid inclusion} (or $Q$ is a {\it relatively rigid
subalgebra} of $M$) if any deformation of $id_M$ by subunital
subtracial c.p. maps $\{\phi_n\}_n$ tends uniformly to $id_Q$ on
the unit ball of $Q$, i.e. $\lim_n {\text{\rm sup}}
\{\|\phi_n(y)-y\|_2 \mid y \in Q, \|y\| \leq 1\} =0$. This
property doesn't in fact depend on the choice of the trace $\tau$
on $M$ and can be given several other equivalent characterizations
(see [P5], [PeP]). The following provides yet another
characterization of relative rigidity, by showing that it is
enough to consider deformations by unital tracial c.p. maps:

\proclaim {3.1. Theorem}  Let $N$ be a finite von Neumann algebra
with countable decomposable center.  Let $Q$ be a von Neumann
subalgebra, then the following are equivalent:

\vskip .05in \noindent $(i)$.  The inclusion $Q \subset N$ is
rigid.

\vskip .05in \noindent $(ii)$.  There exists a normal faithful
tracial state $\tau$ on $N$ such that: $\forall \varepsilon > 0$,
$\exists F = F(\varepsilon) \subset N$ finite and $\delta =
\delta(\varepsilon) > 0$ such that if $\phi : N \rightarrow N$ is
a normal, c.p. map with $\tau \circ \phi = \tau$, $\phi(1) = 1$,
and $\| \phi(x) - x \|_2 \leq \delta$, $\forall x \in F$, then $\|
\phi(b) - b \|_2 \leq \varepsilon$, $\forall b \in Q$, $\| b \|
\leq 1$.

\vskip .05in \noindent $(iii)$.  Condition $(ii)$ is satisfied for
any normal faithful tracial state $\tau$ on $N$.
\endproclaim

\noindent {\it Proof}.  $(i) \implies (iii)$ and $(iii) \implies
(ii)$ are both trivial and so it is enough to show $(ii) \implies
(i)$. That is, assuming $(ii)$ holds we must show that the
following condition holds:

\vskip .05in \noindent $(3.1.1)$ $\forall \varepsilon > 0$,
$\exists F' = F'(\varepsilon) \subset N$ finite and $\delta' =
\delta'(\varepsilon) > 0$ such that if $\phi : N \rightarrow N$ is
a normal, c.p. map with $\tau \circ \phi \leq \tau$, $\phi(1) \leq
1$, and $\| \phi(x) - x \|_2 \leq \delta'$, $\forall x \in F'$,
then $\| \phi(b) - b \|_2 \leq \varepsilon$, $\forall b \in Q$,
$\| b \| \leq 1$.

\vskip .05in

By (Lemma 3 in [PeP]) we may also assume that in the above
condition $\phi$ is symmetric, i.e. $\tau(\phi(x)y) =
\tau(x\phi(y))$, $\forall x,y \in N$. Let $F = F(\varepsilon/2)$,
and $\delta = \delta(\varepsilon/2)$ be given from $(ii)$. Let $F'
= F \cup \{ 1 \}$, and $\delta' = \min_{x \in F} \{ \delta/2,
\delta/(8\|x\|^2+1), \varepsilon^2/8 \}$, suppose that $\phi : N
\rightarrow N$ is a normal, symmetric c.p. map with $\tau \circ
\phi \leq \tau$, $\phi(1) \leq 1$, and $\| \phi(x) - x \|_2 \leq
\delta'$, $\forall x \in F'$. Let $a = \phi(1) = d\tau \circ \phi
/ d\tau$ and define $\phi'$ by $\phi'(x) = \phi(x) + (1-a)^{1/2} x
(1-a)^{1/2}$. Then $\phi'$ is a normal, c.p. map with $\phi'(1) =
1$. Moreover, as $\phi$ is symmetric, so is $\phi'$ and hence
$\tau \circ \phi' = \tau$.

Also, it follows that for each $x\in F$ we have
$$
\| \phi'(x) - x \|_2
      \leq \| \phi(x) - x \|_2 + \|(1-a)^{1/2} x (1-a)^{1/2}\|_2
$$
$$
      \leq \| \phi(x) - x \|_2 + \|(1-a) x (1-a)\|_2^{1/2} \| x \|^{1/2}
$$
$$
      \leq \| \phi(x) - x \|_2 + \|1 - a \|_2 \| x \|
      \leq \delta.
$$
Hence, by $(ii)$ we have $\| \phi'(b) - b \|_2 \leq
\varepsilon/2$, $\forall b \in Q$, $\|b\| \leq 1$. Thus
$$\| \phi(b) - b \|_2
              \leq \| \phi'(b) - b \|_2 + \|(1-a)^{1/2} x (1-a)^{1/2}\|_2
$$
$$
              \leq \| \phi'(b) - b \|_2 + \|1 - a \|_2
              \leq \varepsilon,
$$
for all $b \in Q$ with $\| b \| \leq 1$. \hfill $\square$

\vskip .1in \proclaim {3.2. Corollary}  Let $(M_1, \tau_1)$ and
$(M_2, \tau_2)$ be finite von Neumann algebras and $Q \subset M_1$
a von Neumann subalgebra such that the inclusion $Q \subset M_1 *
M_2$ is rigid.  Then the inclusion $Q \subset M_1$ is rigid.
\endproclaim

\noindent {\it Proof}.  Let $\varepsilon >0$ be given.  Since $Q
\subset M_1 * M_2$ is rigid we can find $F = F(\varepsilon)\subset
M$ finite and $\delta = \delta(\varepsilon)>0$ to satisfy
condition $(ii)$ of Theorem 3.1.  By Lemma 2.1.1 there exists
$F' \subset M_1$ finite, and $\delta' > 0$ such that if $\phi :
M_1 \rightarrow M_1$ is a normal unital, tracial c.p. map such
that $\| \phi(x) - x \|_2 \leq \delta'$, $\forall x \in F'$ then
$\| \phi * id_{M_2} (x) - x \|_2 \leq \delta$, $\forall x \in F$.
Hence by our choice of $F$ and $\delta$ we have that $\| \phi (b)
- b \|_2 = \| \phi * id_{M_2} (b) - b \|_2 \leq \varepsilon$,
$\forall b \in Q$, $\| b \| \leq 1$.  Thus, by 3.1 $Q \subset M_1$
follows a rigid inclusion. \hfill $\square$

We'll now use the deformation 2.2.2 to exploit the relative
rigidity of subalgebras $Q \subset M_1 *_B M_2 \subset \tilde{M}_1
*_B \tilde{M}_2$. This ``deformation/rigidity'' argument is
inspired from (4.3-4.8 in [P3] and Sec. 4 in [P6I])

\proclaim{3.3. Proposition} Let $(M_1, \tau_1), (M_2,\tau_2)$ be
finite von Neumann algebras with a common von Neumann subalgebra
$B\subset M_i, i=1,2,$ such that $\tau_{1|B}=\tau_{2|B}$. Denote
$M=M_1 *_B M_2$, $\tilde M_i = M_i *_B (B \overline{\otimes}
L(\Bbb Z))$, $i=1,2,$ and $\tilde M= \tilde M_1 *_B \tilde M_2 = M
*_B (B \overline{\otimes} L(\Bbb F_2))$, as in $2.2.2$. Denote
$\theta = {\text{\rm Ad}}(u_1)* {\text{\rm Ad}}(u_2)\in \text{\rm
Aut}(\tilde M)$, where $u_1, u_2 \in L(\Bbb F_2)$ are the
canonical generators of $L(\Bbb F_2)$, as in $2.2.2$. Let
$Q\subset M$ be a von Neumann subalgebra such that $Q\subset
\tilde M$ is a rigid inclusion and assume no corner of $Q$ can be
embedded into $B$ inside $M$, i.e. $Q'\cap \langle M, B \rangle$
contains no non-zero finite projections. Then there exists a
non-zero partial isometry $v\in \tilde M$ such that
$vy=\theta(y)v,\forall y\in Q$.
\endproclaim
\vskip .05in
\noindent
{\it Proof}. Since $Q\subset \tilde M$ is rigid, there
exist $F\subset \tilde M$ finite and $\delta>0$ such
that if
 $\phi:\tilde M\rightarrow \tilde M$
  is a subunital, subtracial, c.p. map with
$\|\phi(x)-x\|_2\leq\delta,\forall x\in F$, then
$\|\phi(u)-u\|_2\leq 1/2,\forall u\in\Cal U(Q)$. Using the
continuity of $t\mapsto \theta_t$, we can find $n\geq 1$ such that
$\|\theta_{1/{2^n}}(x)-x\|\leq\delta,\forall x\in F$, which
entails $ \|\theta_{1/{2^n}}(u)-u\|\leq 1/2,\forall u\in\Cal
U(Q)$.

Now, let $a$ be the unique element of minimal norm-$\|\cdot \|_2$
in $K=\overline{co}^w\{\theta_{1/{2^n}}(u)u^*\mid u\in\Cal
U(Q)\}$. From $\|\theta_{1/{2^n}}(u)u^*-1\|_2\leq 1/2,\forall
u\in\Cal U(Q)$, we get that $\| a-1\|_2\leq 1/2$, so $a\not= 0$.
Since $\theta_{1/2^n}(u)Ku^*=K$ and $
\|\theta_{1/2^n}(u)au^*\|_2=\|a\|_2,\forall u\in \Cal U(Q)$ we
deduce, using the uniqueness of $a$, that
$au=\theta_{1/{2^n}}(u)a,\forall u\in\Cal U(Q)$ . Using standard
arguments we can replace $a$ by the partial isometry in its polar
decomposition, thus getting a non-zero partial isometry $v\in M$
such that $vu=\theta_{1/{2^n}}(u)v,\forall u\in\Cal U(Q).$

In what follows we show by induction that for any
$k\geq 0$ there exist a non-zero partial isometry
$v_k\in\tilde M$ such that
$$
v{_k}u=\theta_{1/{2^{n-k}}}(u)v_k,\forall u\in\Cal U(Q), \tag
3.3.1
$$
which for $k=n$ gives the conclusion. Since we have already
constructed $v_0$ we only need to construct $v_{k+1}$, given
$v_k$. Note first that if $v_k$ satisfies $(3.3.1)$ then
$v{_k}^*v{_k}\in Q'\cap \tilde M$ and $v{_k}v{_k}
^*\in\theta_{1/{2^{n-k}}}(Q)'\cap\tilde M. $

By applying the automorphism $\beta$ to the equality $(3.3.1)$ and
using the properties of $\beta$   we get:
$$\beta(v_k)u=\beta(v_k)\beta(u)=
\beta(\theta_{1/{2^{n-k}}}(u))\beta(v_k) \tag 3.3.2
$$
$$
=\theta_{-1/{2^{n-k}}}(\beta(u))\beta(v_k)=
\theta_{-1/{2^{n-k}}}(u)\beta(v_k),\forall
u\in\Cal U(Q).
$$

By replacing $u$ by $u^*$ and taking conjugates in $(3.3.2)$, we
obtain
$$
u\beta(v{_k}^*)=\beta(v{_k}^*)\theta_{-1/{2^{n-k}}}(u),\forall
u\in\Cal U(Q),
$$
which combined with $(3.1.1)$ gives
$$
v_k\beta({v_k}^*)\theta_{-1/{2^{n-k}}}(u)= v_ku\beta(v{_k}^*) \tag
3.3.3
$$
$$
=\theta_{1/{2^{n-k}}}(u)v_k\beta({v_k}^*),\forall u\in\Cal U(Q).
$$

By applying $\theta_{1/{2^{n-k}}}$ to $(3.3.3)$, we further get
for $u \in \Cal U(Q)$ the identity
$$
\theta_{1/{2^{n-k}}}(v_k\beta({v_k}^*))u=\theta_{1/{2^{n-k-1}}}(u)
 \theta_{1/{2^{n-k}}}(v_k\beta({v_k}^*)). \tag 3.3.4
$$

Since no corner of $Q$ can be embedded into $B$ inside $M$ we can
apply Theorem 1.1 to conclude that $Q'\cap \tilde M\subset M$.
Thus, since $v{_k}^*v{_k}\in Q'\cap \tilde M$ and $M\subset{\tilde
M}^{\beta}$, we get $\beta(v{_k}^*v{_k})=v{_k}^*v{_k}$, implying
that $v_k\beta({v_k}^*)$ is a partial isometry of same left
support as $v_k$. Thus, by $(3.3.4)$, $w=\theta_{1/{2^{n-k}}}(v_k
\beta(v_k^*))$ is a non-zero partial isometry satisfying
$wu=\theta_{1/{2^{n-k-1}}}(u) w,\forall u\in\Cal U(Q)$, and the
inductive step follows. \hfill $\square$

\proclaim{3.4. Proposition} As in $2.3$, denote by $\widehat{E}_L$
the orthogonal projection of $L^2(M)$ onto the Hilbert space
spanned by reduced words of length $\geq L$. If $Q \subset M_1 *_B
M_2$ is a rigid inclusion, then for any $\varepsilon > 0$, there
exists $L \in \Bbb N$ such that $\| \widehat {E}_L(x) \|_2 <
\varepsilon$, $\forall x \in Q, \| x \| \leq 1$.
\endproclaim
\noindent {\it Proof}. Let $\phi_n: M \rightarrow M$ be as in
2.3.1, for some $c_n \nearrow 1$, then by 2.1 we have $\lim_n \|
\phi_n(x) - x \|_2 = 0$, $\forall x \in M$. Thus, since $Q \subset
M$ is rigid, there exists $l \in \Bbb N$ such that $\| \phi_l(x) -
x \|_2 < \varepsilon /2$, $\forall x \in Q, \| x \| \leq 1$.  Let
$L \in \Bbb N$ such that $c_l^L < \varepsilon /2$. Then
$$
\| \widehat {E}_L (x) \|_2
  \leq \| \widehat {E}_L (x - \phi_l(x)) \|_2 + \| \phi_l
  \circ \widehat {E}_L (x) \|_2
$$
$$
\leq \| x - \phi_l(x) \|_2 + c_l^L \| x \|_2
  < \varepsilon, \forall x \in Q, \| x \| \leq 1.
$$
\hfill $\square$

\heading 4. Existence of intertwining bimodules. \endheading

In the previous section we saw that a relatively rigid subalgebra
$Q$ of a finite AFP  von Neumann algebra $M=M_1 *_B M_2$ can be
``located'' by certain c.p. deformations of $id_M$. In this
section we will use this information to prove that $L^2(M)$ must
contain non-trivial Hilbert bimodules intertwining $Q$ into either
$M_1$ or $M_2$. The rather long and technical proof will proceed
by contradiction, assuming $Q$ cannot be intertwined in neither
$M_1$ nor $M_2$, inside $M$. We first show that this implies $Q$
cannot be intertwined into $M_1, M_2$ in $\tilde{M}=(M_1*_B
M_2)*_B (B \overline{\otimes} L(\Bbb F_2))$ either. By Proposition
3.4 and (2.3 in [P6]), this shows that $Q$ must contain ``at
infinity'' elements with uniformly bounded free length and at
least two ``very large letters'' in $M_1, M_2$ or $L(\Bbb F_2)$.
This will be shown to contradict Proposition 3.3.

To ``measure'' the letters in $M_i$ we will need these algebras to
have nice orthonormal basis over $B$, in the following sense:

\vskip .1in \noindent {\it 4.1. Definition}. Let $(M, \tau)$ be
separable finite von Neumann algebra and $B \subset M$ a von
Neumann subalgebra. A sequence of elements $\{\eta_n\}_n \subset
M$ satisfying the conditions $\eta_0 = 1$,
$E_B(\eta_i^*\eta_j)=\delta_{ij}$,
$\forall i,j$ and $\Sigma \eta_n B$ dense in $L^2(M,\tau)$ is
called a {\it bounded, homogeneous, orthonormal basis} (BHOB) of
$M$ over $B$. An inclusion $B \subset M$ having a BHOB is said to
be {\it homogeneous}.

\proclaim{4.2. Lemma} Let $(M,\tau)$ be a separable finite von
Neumann algebra and $B \subset M$ a von Neumann subalgebra. Assume
one of the following conditions holds true:

\vskip .05in \noindent $(a)$. $B = \Bbb C1$.

\vskip .05in \noindent $(b)$. $B=A\subset M$ is Cartan (i.e. $B$
is maximal abelian and regular in $M$) and $M$ is type $\text{\rm
II}_1$.

\vskip .05in \noindent $(c)$. $B=N\subset M$ is an irreducible
inclusion of $\text{\rm II}_1$ factors and $B$ is regular in $M$.

\vskip .05in \noindent Then $B \subset M$ is homogeneous, moreover
in both cases $(b)$ and $(c)$, $M$ has a $\text{\rm BHOB}$ made of
unitary elements in $\Cal N_M(B)$.
\endproclaim
\noindent {\it Proof}. Case $(a)$ is clear by the Gram-Shmidt
algorithm. The case $(c)$ is trivial once we notice that such
$N\subset M$ is a crossed product inclusion $N \subset M=N
\rtimes_{\sigma, v} G$ for some cocycle action $\sigma$ of a
discrete countable group $G$ on $N$, with $2$-cocycle $v$. Indeed,
in this case the canonical unitaries $\{u_g\}_g \subset M$
implementing the action $\sigma$ provide a BHOB of $M$ over $N$.

To prove the case $(b)$, we first show the that given any $n \geq 1$
and any finite set $v_0=1, v_1, ..., v_{n-1} \in \Cal N_M(A)$, with
$E_A(v_i^*v_j)=0$ for $0 \leq i\neq j <n$, and $v'_n \in \Cal G\Cal
N_M(A)$, with ${v_n'}^*v_n' \neq 1$, there exists $v \in \Cal G\Cal
N_M(A)$ non-zero such that $v^*v_n'=0, v_n'v^*=0$ and
$E_A(v^*v_j)=0, \forall 1\leq j \leq n-1$. Indeed, note first that
by ([Dy]) there exists $u\in \Cal N_M(A)$ such that $v_n' = u(1-p)$,
where $p=1-{v_n'}^*v_n'\in \Cal P(A)$. Since $Ap \subset pMp$ is
Cartan with $pMp$ of type II$_1$, it follows that there exists $w
\in \Cal N_{pMp}(Ap)$ such that $w \not\in (\Sigma_{i=0}^{n-1} pu^*
v_i Ap)$ (or else, $pMp$ would follow finite dimensional over $Ap$
thus of type I, implying that $M$ has a type I direct summand, a
contradiction). Then $w_i=pu^*v_ip \in \Cal G\Cal N(Ap)$, $0 \leq i
\leq n-1$, satisfy $0\neq w- \Sigma_{i=0}^{n-1} w_iE_A(w_i^*w) \in
\Cal G\Cal N(A)$ (by [Dy]). But then $v = u(w- \Sigma_{i=0}^{n-1}
w_iE_A(w_i^*w))$ clearly satisfies all required conditions.

Now, to finish the proof of $(b)$ let $\{u_n \}_{n \geq 0} \subset
\Cal N_M(A)$ be a sequence of unitaries normalizing $A$, dense in
$\Cal N_M(A)$ in the norm $\|\cdot \|_2$, with $u_0=1$ and with
each $u_n$ appearing with infinite multiplicity. It is sufficient
to  construct a sequence $\{v_n\}_{n \geq 0} \subset \Cal N_M(A)$
such that for all $m \geq 0$ we have:
$$
E_A(v_i^*v_j)=\delta_{ij}, \forall m \geq i,j \geq 0, u_m \in
\Sigma_{i=0}^m v_i A. \tag 4.2.1
$$
Assume we have constructed $v_0=1, v_1, ..., v_{n-1}$ satisfying
$(4.2.1)$ for $m=n-1$, for some $n \geq 1$. Let $v_n''= u_n -
\Sigma_{i=0}^{n-1} v_i E_A(v_i^*u_n)\in \Cal G\Cal N_M(A)$. Let
$v_n' \in \Cal G\Cal N_M(A)$ be maximal with the properties:
$v_n'(v_n'')^*v_n''=v_n''$ and   $v_n'\perp \Sigma_{i=0}^{n-1} v_i
A$. By applying the first part of the proof to $v_0, v_1, ...,
v_{n-1}, v_n'$  and using the maximality, it follows that $v_n'$ is
a unitary element. But then $v_n=v_n'$ clearly satisfies $(4.2.1)$.
\hfill $\square$

\proclaim{4.3. Theorem} Let $(M_1, \tau_1)$ and $(M_2, \tau_2)$ be
finite von Neumann algebras and $B\subset M_i$ a common von
Neumann subalgebra such that $\tau_{1|B} = \tau_{2|B}$. Suppose
that the inclusions $B \subset M_1$ and $B \subset M_2$ are
homogeneous. Let $Q$ be a von Neumann subalgebra of $M=M_1 *_B
M_2$ such that the inclusion $Q \subset M$ is rigid. Then for
either $i=1$ or $i=2$ there exists a non-zero projection $f$ in
$Q' \cap \langle M, e_{M_i} \rangle$ of finite trace
$Tr=Tr_{\langle M, e_{M_i} \rangle}$.
\endproclaim

\vskip .05in \noindent {\it Proof}. By taking spectral
projections, it is sufficient to show there exists $i \in \{1,2\}$
and $a \in Q'\cap \langle M, e_{M_i} \rangle$ with $ 0 \leq a \leq
1$, $0 \neq Tr(a) < \infty$. Assume by contradiction that there
are no such elements in $Q' \cap \langle M, e_{M_i} \rangle$,
$i=1,2$. If we identify $Q$ with the diagonal subalgebra $\{x
\oplus x \mid x\in Q\}$ in $M \oplus M$, then this is equivalent
to saying that $Q$ cannot be intertwined into $M_1 \oplus M_2$
inside $M \oplus M$. By (2.3 in [P6]) this implies the following:

\vskip .05in \noindent {\it Fact} 1:  $\forall \varepsilon > 0$,
$\forall y_1, \ldots, y_n \in M$, $\exists w \in \Cal U(Q)$ such
that $\| E_{M_i}(y_j w y_k^*) \|_2 < \varepsilon$, $\forall i \in
\{ 1,2 \}, \forall j,k \in\{1,\dots n\}$. \vskip .05in

For the next part of the proof we need to introduce some notation.
Thus, for each $i \in \{ 1,2 \}$ let $\{ \xi_i^j \}_j \subset M_i$
be a BHOB of $M_i$ over $B$.
Also take $\{ \xi_3^j \}_j = \{ u_g \}_{g \in \Bbb F_2} \subset L(\Bbb F_2)$
such that $\xi_3^0 = u_e = 1$. Then
using the notation $\tilde M = M *_B (B \otimes L(\Bbb F_2))$ as
in Proposition 3.3 a simple exercise shows that
$$
\beta =
\{ 1 \} \cup
\{
  \xi_{i_1}^{j_1} \cdots  \xi_{i_n}^{j_n} \mid
  n \in \Bbb N,
  i_k \in \{ 1,2,3 \},
  j_k \geq 1,
  i_1 \not= i_2, \cdots, i_{n-1} \not= i_n
\}
$$
is BHOB of $\tilde M$ over $B$.

For each $n_0 \in \Bbb N$ let
$$
S_{n_0} =
\{ 1 \} \cup
\{
  \xi_{i_1}^{j_1} \cdots  \xi_{i_n}^{j_n} \mid
  n \leq n_0,
  i_k \in \{ 1,2 \},
  1 \leq j_k \leq n_0,
  i_1 \not= i_2, \cdots, i_{n-1} \not= i_n
\},
$$
$$
\tilde S_{n_0} =
\{ 1 \} \cup
\{
  \xi_{i_1}^{j_1} \cdots  \xi_{i_n}^{j_n} \mid
  n \leq n_0,
  i_k \in \{ 1,2,3 \},
  1 \leq j_k \leq n_0,
  i_1 \not= i_2, \cdots, i_{n-1} \not= i_n
\}.
$$
Also for $i \in \{ 1,2 \}$ let $\tilde S_{n_0}^{R,i}$
(resp. $\tilde S_{n_0}^{L,i}$) be the subset of
$\tilde S_{n_0}$ consisting of
$1$ and the vectors in $\tilde S_{n_0}$ such that
$i_n \not= i$ (resp. $i_1 \not= i$).
Note that if $\zeta, \zeta' \in \tilde S_{n_0}^{R,i}$
and $x \in M_i$ such that $E_B(x) = 0$, then
$\| \zeta' x \zeta^* \|_2 = \| x \|_2$.  Also if
$\zeta \in \tilde S_{n_0}$ and $b \in B$ then
$\| \zeta b \|_2 = \| b \|_2$.

We now strengthen {\it Fact} 1 so that the elements $y_1, \ldots,
y_n$ may be taken in $\tilde M$:

\vskip .05in \noindent {\it Fact} 2:  $\forall \varepsilon > 0$,
$\forall y_1, \ldots, y_n \in \tilde M$, $\exists w \in \Cal U(Q)$ such
that $\| E_{M_i}(y_j^* w y_k) \|_2 < \varepsilon$, $\forall i \in
\{ 1,2 \}, \forall j,k \in\{1,\dots n\}$. \vskip .05in

As our basis for $\tilde M$ is made up of bounded vectors by first
approximating the $y_k$'s on the right of $w$ and then approximating
the $y_j$'s on the left of $w$ we may assume that all of the $y_j$'s
are basis elements and then use the triangle inequality to deduce
the general case.  Also as $E_{M_i}$ is $B$-bimodular it enough to
suppose that the $y_j$'s all lie in $\beta$. Thus we only need to
show that $\forall \varepsilon > 0$, $\forall n_0 \in \Bbb N$,
$\exists w \in \Cal U(Q)$ such that $\| E_{M_i}(\zeta^* w \zeta')
\|_2 < \varepsilon$, $\forall \zeta, \zeta' \in \tilde S_{n_0}$.

To prove this, we first use {\it Fact} 1 to deduce
that there exists $w \in \Cal U(Q)$ such that
$\| E_{M_i}( \zeta_0^* w \zeta_0' ) \|_2 < \varepsilon$,
$\forall \zeta_0, \zeta_0' \in S_{n_0}$.
Then if $\zeta, \zeta' \in \tilde S_{n_0}$ we may find
$\zeta_1, \zeta_1' \in \tilde S_{n_0}^{L,1} \cap \tilde S_{n_0}^{L,2}$, and
$\zeta_0, \zeta_0' \in S_{n_0}$ such that
$\zeta = \zeta_0 \zeta_1$, and $\zeta' = \zeta_0' \zeta_1'$.
If $\zeta_1 = \zeta_1' = 1$ then from above we have
$\| E_{M_i}( \zeta^* w \zeta' ) \|_2 =
\| E_{M_i}( \zeta_0^* w \zeta_0') \|_2 < \varepsilon$.
Otherwise we have
$\| E_{M_i}( \zeta^* w \zeta' ) \|_2 \leq
\| E_M( \zeta^* w \zeta' ) \|_2
 \leq \| \zeta_1^* E_B(\zeta_0^* w \zeta_0') \zeta_1' \|_2
 \leq \| E_{M_1}(\zeta_0^* w \zeta_0') \|_2 < \varepsilon$.
This proves {\it Fact} 2.

\vskip .05in

We  continue by showing that there are elements of $\Cal U(Q)$
(``at infinity'') which are almost orthogonal to the subspaces
having at most one ``large letter'' from $M_1, M_2$. Specifically,
let $\Cal H_{n_0} = \overline {\text{sp}} (\tilde S_{n_0} M_1
{\tilde S_{n_0}}^* \cup \tilde S_{n_0} M_2 {\tilde S_{n_0}}^*)
\subset L^2(\tilde M)$.

Suppose $\zeta_1, \zeta_1' \in  \tilde S_{n_0}^{R,i}$, $\zeta_2,
\zeta_2' \in \tilde S_{n_0}^{R,j}$, then $\forall b_1, b_2 \in B$,
$\forall K, L > n_0$ we have that $E_B((\zeta_1' \xi_i^K b_1
\zeta_1^*)^* (\zeta_2' \xi_j^L b_2 \zeta_2^*)) = \delta_{\zeta_1'
\zeta_2'} E_B((\xi_i^K b_1 \zeta_1^*)^* (\xi_j^L b_2 \zeta_2^*))$,
and also we have that \newline $E_B((\xi_j^L b_2 \zeta_2^*)
(\xi_i^K b_1 \zeta_1^*)^*) = \delta_{\zeta_1 \zeta_2} E_B((\xi_j^L
b_2) (\xi_i^K b_1)^*)$.  Hence
$$
\langle \zeta_1' \xi_i^K b_1 \zeta_1^*, \zeta_2'
\xi_j^L b_2 \zeta_2^* \rangle
= \tau \circ E_B((\zeta_1' \xi_i^K b_1 \zeta_1^*)^*
(\zeta_2' \xi_j^L b_2 \zeta_2^*))
$$
$$
= \delta_{\zeta_1' \zeta_2'} \tau \circ E_B((\xi_i^K b_1\zeta_1^*)^*
(\xi_j^L b_2 \zeta_2^*)) = \delta_{\zeta_1' \zeta_2'}
\delta_{\zeta_1 \zeta_2} \delta_{ij} \langle \xi_i^K b_1, \xi_j^L
b_2 \rangle.
$$
Also if $\zeta \in \tilde S_{n_0}$, and $b \in B$ then
$\langle \zeta_1' \xi_i^K b_1 \zeta_1^*, \zeta b \rangle =
\tau \circ E_B((\zeta b)^* \zeta_1' \xi_i^K b_1 \zeta_1^*) = 0$.

Hence we have the following direct sum of orthogonal
subspaces of $\Cal H_{n_0}$:
$$
\Cal H_{n_0}' = \bigoplus_{\zeta \in \tilde S_{n_0}} \Cal H_{\zeta}
    \oplus \bigoplus_{i \in \{ 1,2 \}, \zeta, \zeta' \in \tilde S_{n_0}^{R,i}}
            \Cal H_{i, \zeta', \zeta},
$$
where $\Cal H_{\zeta} =\overline {\zeta B}$, and
$\Cal H_{i, \zeta', \zeta} = \overline{\text {sp}}
\zeta' \{ \xi_i^K  \}_{K > n_0} B \zeta^*$.

Since $\{ \xi_3^j \}_j = \{ u_g \}_{g \in \Bbb F_2}$
we may find $m_0 > 2n_0$ such that
$\{ \xi_3^j \}_{j \leq n_0} \{ \xi_3^j \}_{j \leq n_0}^*
\subset \{ \xi_3^j \}_{j \leq m_0}$.
We will show that $\Cal H_{n_0} \subset \Cal H_{m_0}'$.

Let $\Cal K_0$ be the Hilbert space generated by all vectors of the form
$\eta' b \eta^*$ where $b \in B$, $\eta =
\xi_{i_1}^{j_1} \cdots  \xi_{i_n}^{j_n}$,
$\eta' = \xi_{k_1}^{l_1} \cdots  \xi_{k_m}^{l_m}$, such that
$i_n \not= k_m$, $n + m \leq m_0$, and $j_p, l_p \leq m_0$, $\forall p$.
If $\zeta, \zeta' \in \tilde S_{n_0}$, and $x \in M_i \ominus B$,
Then we may find $\zeta_1, \zeta_1' \in \tilde S_{n_0}^{R,i}$, and
$\zeta_0, \zeta_0' \in \tilde S_{n_0} \cap M_i$ such that
$\zeta = \zeta_1 \zeta_0$, and $\zeta' = \zeta_1' \zeta_0'$.
We have that if $P$ is the projection onto the subspace
$\overline {\text {sp}} \{ \xi_i^j \}_{j = 1}^{m_0} B$ then
$\zeta_1' (\zeta_0' x \zeta_0^* - P(\zeta_0' x \zeta_0^*) -
E_B(\zeta_0' x \zeta_0^*)) \zeta_1^* \in \Cal H_{m_0}'$,
$\zeta_1' P(\zeta_0' x \zeta_0^*) \zeta_1^* \in \Cal K_0$, and
$\zeta_1' E_B(\zeta_0' x \zeta_0^*) \zeta_0^* \in
\tilde S_{n_0} B \tilde S_{n_0}^*$.

If $\zeta, \zeta' \in \tilde S_{n_0}$ and $b \in B$ then if $\zeta$
and $\zeta'$ do not end with a letter in the same algebra we have that
$\zeta' b \zeta^* \in \Cal K_0$, also if both $\zeta$ and $\zeta'$ end with
something in $\{ \xi_3^j \}_j$ then as
$L(\Bbb F_2)$ commutes with $B$ and since
$\{ \xi_3^j \}_{j \leq n_0} \{ \xi_3^j \}_{j \leq n_0}^*
\subset \{ \xi_3^j \}_{j \leq m_0}$
we may rewrite $\zeta' b \zeta^*$ to see that it is in $\Cal K_0$,
otherwise as above we may find
$\zeta_1, \zeta_1' \in \tilde S_{n_0}^{R,i}$, and
$\zeta_0, \zeta_0' \in \tilde S_{n_0} \cap M_i$ such that
$\zeta = \zeta_1 \zeta_0$, and $\zeta' = \zeta_1' \zeta_0'$,
and then decompose $\zeta' b \zeta^*$ into parts in
$\Cal H_{m_0}'$, $\Cal K_0$, and something in
$\tilde S_{n_0} B \tilde S_{n_0}$ with shorter words.
Hence by induction to show that $\Cal H_{n_0} \subset
\Cal H_{m_0}'$ it is enough
to show that $\Cal K_0 \subset \Cal H_{m_0}'$.

Let $\eta$ and $\eta'$ be as above, and take $b \in B$, if $n =
0$, i.e. $\eta = 1$ then $\eta' b \eta^* = \eta' b \in \Cal
H_{m_0}$.  Also if $i_n = 3$ then since $L(\Bbb F_2)$ commutes
with $B$ and ${\xi_{i_n}^{j_n}}^* \in \{ \xi_3^j \}_{j \leq m_0}$
we can rewrite $\eta' b \eta^*$ so that $\eta$ and $\eta'$ are
still in $\tilde S_{m_0}$ but such that the length of $\eta$ is
shorter. If $i_n = i \in \{ 1,2 \}$ then since $E_B(b
{\xi_{i_n}^{j_n}}^*) = 0$ as above we may replace $b
{\xi_{i_n}^{j_n}}^*$ with $P(b {\xi_{i_n}^{j_n}}^*)$ and $b
{\xi_{i_n}^{j_n}}^* - P(b {\xi_{i_n}^{j_n}}^*)$ and in so doing
rewrite $\eta' b \eta^*$ as a sum of things where the word on the
right has shorter length plus something in $\Cal H_{m_0}'$.  Thus
by induction we have shown that $\Cal K_0 \subset \Cal H_{m_0}'$
and so $\Cal H_{n_0} \subset \Cal H_{m_0}'$.

Let $P_{n_0}$ be the orthogonal projection of
$L^2(\tilde M)$ onto $\Cal H_{n_0}$.

\vskip .05in \noindent {\it Fact} 3:  $\forall
\varepsilon > 0$, $y_1, \ldots, y_n \in \tilde M$,
$\exists w \in \Cal U(Q)$ such that
$\| P_{n_0}(y_j^* w y_k) \|^2_2 < \varepsilon$, $\forall j,k$. \vskip .05in

Let $m_0$ be as above, then by {\it Fact} 2 $\exists w \in \Cal U(Q)$ such that
$\| E_{M_i}(\zeta^* y_k^* w^* y_j \zeta') \|_2^2
< \varepsilon / 3|\tilde S_{m_0}|^2$,
$\forall j,k \leq n$, $\forall \zeta, \zeta' \in \tilde S_{m_0}$.

Thus $\forall \zeta \in \tilde S_{m_0}$, $b \in B$ we have
$$
| \langle \zeta b, y_j^* w y_k \rangle |^2 = | \tau
(E_B( y_k^* w^* y_j \zeta ) b) |^2
$$
$$
  \leq \| E_{M_1}( y_k^* w^* y_j \zeta ) \|_2^2 \| b \|_2^2
  < (\varepsilon / 3|\tilde S_{m_0}|^2) \| \zeta b \|_2^2,
$$
and so $\| P_{\Cal H_\zeta} (y_j^* w y_k) \|_2^2
< \varepsilon / 3|\tilde S_{m_0}|^2$.
Also $\forall i \in \{ 1,2 \}$, $\zeta, \zeta'
\in \tilde S_{m_0}^{R,i}$, and
$\xi \in \overline {\text {sp}} \{ \xi_i^K  \}_{K > m_0} B$ we have
$$
| \langle \zeta' \xi \zeta^*, y_j^* w y_k \rangle |^2
 = | \tau( E_{M_i}(\zeta^* y_k^* w^* y_j \zeta') \xi ) |^2
$$
$$
 \leq \| E_{M_i} (\zeta^* y_k^* w^* y_j \zeta') \|_2^2 \| \xi \|_2^2
 < (\varepsilon / 3|\tilde S_{m_0}|^2) \| \zeta' \xi \zeta^* \|_2^2,
$$
and so $\| P_{\Cal H_{i, \zeta', \zeta}} (y_j^* w y_k) \|_2^2 <
\varepsilon / 3|\tilde S_{m_0}|^2$.

Therefore
$$
\| P_{n_0} (y_j^* w y_k) \|_2^2 \leq
\| P_{\Cal H_{m_0}'}(y_j^* w y_k) \|_2^2
$$
$$
= (\Sigma_{\zeta \in \tilde S_{m_0}}
\| P_{\Cal H_\zeta} (y_j^* w y_k) \|_2^2 )
    + (\Sigma_{i \in \{ 1,2 \}, \zeta' \in \tilde S_{m_0}^{R,i},
    \zeta \in \tilde S_{m_0}^{L,i}}
         \| P_{\Cal H_{i, \zeta', \zeta}} (y_j^* w y_k) \|_2^2)
  < \varepsilon,
$$
for each $j,k \leq n$, thus we have proved {\it Fact} 3.

\vskip .05in

Next we note that if $Q' \cap \langle M, e_B \rangle$ contains a non-zero,
finite-trace projection then so does $Q' \cap \langle M, e_{M_1}
\rangle$ and so by our assumption we are in the position of
applying Proposition 3.3.

Hence there exists a non-zero
partial isometry $v \in \tilde M$ such that $vy = \theta(y)v$,
$\forall y \in Q$. Let $\varepsilon > 0$ and
take $n_1$ large
enough so that there exists $v_0^* \in {\text {sp}} \tilde S_{n_1}B$
satisfying $\| v - v_0 \|_2 < \varepsilon/ 6$.
Using the notation in Proposition 3.4 let $L \in \Bbb N$, such that
$\| \widehat E_L(x) \|_2 < \varepsilon/ 12 \| v_0 \|^2$,
$\forall x \in Q$, $\| x \| \leq 1$,
also let $n_0 = n_1 + 3L$, and let $m_0$ be as above so
that $\Cal H_{n_0} \subset \Cal H_{m_0}'$.
Then as our basis is bounded we have
that $v_0$ is bounded so by {\it Fact} 3 $\exists w \in
\Cal U(Q)$, such that $\|P_{m_0}(v_0w)\|_2 < \varepsilon/ 4\| v_0 \|$.
Take
$w_0 \in M$ such that
$\widehat E_L(w_0) = 0$, and $\| w - w_0 \|_2 < \varepsilon/6 \| v_0 \|^2$.

Let $\Cal K \subset L^2(\tilde M)$ be the
right Hilbert $B$-module generated by $1$ and all the vectors
$\xi_{i_1}^{j_1} \cdots  \xi_{i_n}^{j_n} \in \beta$ such that
if $i_{k_0} = 3$ for some $k_0 \leq n$ then
$j_k \leq n_1$, $\forall k < k_0$.
Note that $\Cal K M \subset \Cal K$ and so since
$w_0 \in M$ and $v_0 \in \Cal K$
we have that $v_0w_0 \in \Cal K$.

Let $\zeta, \zeta' \in \tilde S_{m_0}^{R, i}$, $K > m_0$, and $b \in B$,
then since $K > n_1$ we have that
$P_{\Cal K}(\zeta' \xi_i^K b \zeta) = \zeta' \xi_i^K b \zeta$ if
$\zeta \in M$ and $\zeta' \in \Cal K$, and
$P_{\Cal K}(\zeta' \xi_i^K b \zeta) = 0$ otherwise.
Also if $\zeta \in \tilde S_{m_0}$ and $B \in B$ then
$P_{\Cal K}(\zeta b) = \zeta b$ if $\zeta \in \Cal K$
and $P_{\Cal K}(\zeta b) = 0$ otherwise.
Hence $P_{\Cal K}(\Cal H_{n_0}) \subset P_{\Cal K}(\Cal H_{m_0}')
                \subset \Cal H_{m_0}' \subset \Cal H_{m_0}$.

Let us write $v_0^*$ and $w_0$ in $\beta$ as
$v_0^* = \Sigma_{\xi_v \in \beta} \xi_v b_{\xi_v}$ and
$w_0 = \Sigma_{\xi_w \in \beta} \xi_w b_{\xi_w}$.
Take $\xi_v, \xi_w \in \beta$ such that $b_{\xi_v} \not= 0$
and $b_{\xi_w} \not= 0$,
where $\xi_v = \xi_{i_1}^{j_1} \cdots  \xi_{i_n}^{j_n}$,
$\xi_w = \xi_{k_1}^{l_1} \cdots  \xi_{k_m}^{l_m}$.
Thus
$\theta (\xi_w b_{\xi_w}) (\xi_v b_{\xi_v})^*
= (u_{k_1} \xi_{k_1}^{l_1} u_{k_1}^* \cdots u_{k_m}
\xi_{k_m}^{l_m} u_{k_m}^* b_{\xi_w})
    (b_{\xi_v}^* {\xi_{i_n}^{j_n}}^* \cdots {\xi_{i_1}^{j_1}}^*)$.
Let us assume that $n \geq 3m$ by adding on $1$'s at
the end of this word if necessary.
If $k_s \leq n_0$, $\forall 1 \leq s < m$
then since $v_0^* \in {\text {sp}} \tilde S_{n_1}B$, and $m < L$,
by decomposing $u_{k_m}^* b_{\xi_w} b_{\xi_v}^* {\xi_{i_n}^{j_n}}^*$
as it's expectation onto $B$ plus something with terms in
$B \otimes L(\Bbb F_2)$ and zero expectation onto $B$,
we write $\theta (\xi_w b_{\xi_w}) (\xi_v b_{\xi_v})^*$ as something
in $\Cal H_{n_0}$ plus something in $\Cal K^{\perp}$.
Otherwise if $k_s > n_0$ for some $s < m$ then by decomposing
$u_{k_m} (\xi_{k_m}^{l_m} (u_{k_m}^* b_{\xi_w}
    b_{\xi_v}^* {\xi_{i_n}^{j_n}}^*) {\xi_{i_{n-1}}^{j_{n-1}}}^*)
    {\xi_{i_{n-2}}^{j_{n-2}}}^*$
just as above into it's expectation onto $B$ plus something
with terms in $B \otimes L(\Bbb F_2)$ and zero expectation onto $B$,
we write $\theta (\xi_w b_{\xi_w}) (\xi_v b_{\xi_v})^*$ as something
with shorter words plus something in $\Cal K^{\perp}$ .
Hence by induction we have shown that
$\theta (\xi_w b_{\xi_w}) (\xi_v b_{\xi_v})^* \in \Cal H_{n_0}
+ \Cal K^{\perp}$ and hence also
$\theta (w_0)v_0 \in \Cal H_{n_0} + \Cal K^{\perp}$.

As $P_{\Cal K} (\Cal H_{n_0}) \subset \Cal H_{m_0}$ we have that
$P_{\Cal K}(\theta (w_0)v_0) \subset \Cal H_{m_0}$.
Thus
$$| \langle v_0w_0, \theta (w_0)v_0 \rangle |
= | \langle v_0w_0, P_{\Cal K}(\theta (w_0)v_0) \rangle |
= | \langle P_{m_0}(v_0w_0), P_{\Cal K}(\theta (w_0)v_0) \rangle |
$$
$$
\leq \| P_{m_0}(v_0w_0) \|_2 (\| w_0 - w \|_2 \| v_0 \| + \| v_0 \|_2)
$$
$$
\leq (\| P_{m_0}(v_0w) \|_2 + \|v_0\| \|w_0 - w\|_2)
(\| w_0 - w \|_2 \| v_0 \| + \| v_0 \|_2)
$$
$$
< (\varepsilon/4 \| v_0 \| + \varepsilon/6 \| v_0 \|)
(\varepsilon/6 \| v_0 \| + \|v_0\|)
< \varepsilon/4.
$$
Hence we have shown that
$$
\| v \|_2^2 = \| vw \|_2^2
$$
$$
= \langle vw, \theta(w)v \rangle
  \leq 2\| v - v_0 \|_2
    +  | \langle v_0w, \theta(w)v_0 \rangle |
$$
$$
\leq 2\| v - v_0 \|_2 + 2\| v_0 \|^2 \| w - w_0 \|_2
  + | \langle v_0w_0, \theta(w_0)v_0 \rangle |
$$
$$
< \varepsilon/3 + \varepsilon/3 + \varepsilon/4
< \varepsilon,
$$
which contradicts
the assumption that $v$ is non-zero. \hfill $\square$

\heading 5. Rigid subalgebras in AFP factors: General Bass-Serre
type results
\endheading

We have shown in Theorem 4.3 that if $Q$ is a relatively rigid von
Neumann subalgebra of an AFP algebra $M=M_1 *_B M_2$ then there
exists a non-trivial Hilbert bimodule $\Cal H \subset L^2(M)$
intertwining $Q$ into one of the $M_i$'s. We now deduce that a
corner of $Q$ can be conjugated by a unitary element into that
same $M_i$. When $M_1, M_2$ are factors, one can in fact uniquely
partition 1 with projections $q_1, q_2 \in Q'\cap M$ such that
$Qq_i$ is unitary conjugate into $M_i$, $i=1,2$. This general
Bass-Serre type result will be used in the next sections to derive
more specific statements in the cases $B=\Bbb C$, $B=A$ abelian
Cartan and $B= R$ the hyperfinite II$_1$ factor.

\proclaim{5.1. Theorem} Let $(M_i, \tau_i)$, $i=0,1,2$, be finite
von Neumann algebras with a common von Neumann subalgebra $B
\subset M_i$, $i=0,1,2,$ such that
$\tau_{0|B}=\tau_{1|B}=\tau_{2|B}$, and such that the inclusions
$B \subset M_i$ are homogeneous. Denote $M=M_0 *_B M_1*_{B}M_2$.
Let $Q\subset M$ be a relatively rigid diffuse von Neumann
subalgebra. Assume no corner of $Q$ can be embedded into $M_0$
inside $M$.

$1^\circ$. There exist $i \in \{1,2\}$, projections $q\in Q$,
$q''\in Q'\cap M$ with $qq''\neq 0$ and a unitary element $u\in
\Cal U(M)$ such that $uqQqq''u^*\subset M_i$.

$2^\circ$. If $M_1, M_2$ are factors then there exists a unique
pair of projections $q_1, q_2 \in Q'\cap M$ such that $q_1+q_2=1$
and $u_i(Qq_i)u_i^*\subset M_i$ for some unitaries $u_i\in\Cal
U(M)$, $i=1,2$. Moreover, these projections lie in the center of
$Q'\cap M$.
\endproclaim
\noindent {\it Proof}.  $1^\circ$. By Theorem 4.3 and (2.1 in [P6])
there exist $i \in \{0,1,2\}$, $q \in \Cal P(Q)$, $p \in \Cal
P(M_i)$ non-zero projections, an isomorphism $\psi$ of $qQq$ into
$pM_ip$ and a non-zero partial isometry $v \in M$ such that $vv^*
\in (qQq)' \cap qMq$, $v^*v \in \psi(qQq)' \cap pMp$ and $xv =
v\psi(x)$, $\forall x \in qQq$.
By hypothesis, $i$ cannot be equal
to $0$, thus $i\in \{1,2\}$.
Note that by shrinking $q$ if necessary, we may assume
$xv=0$ for $x\in qQq$ implies $x=0$. Also,
if we denote by $q'$ the support projection of $E_{M_1}(v^*v)$,
then by replacing if necessary $\psi$ by $q'\psi(\cdot)q'$
it follows that we may assume $q'=p$.

Now note that if a corner of
$\psi(qQq)$ can be embedded into $pBp$ inside $pM_ip$ then by
Lemma 1.4.5 a corner of $Q$ can
be embedded into $B$ (thus into $M_0\supset
B$ as well) inside $M$, contradicting the hypothesis. Thus,
no corner of $\psi(qQq) \subset
pM_ip$ can be embedded into $B$ inside $M_i$, so we can
apply Theorem 1.2.1 to conclude that $Q_0'\cap pMp \subset pM_ip$.
Hence, $v^*v \in Q_0'\cap pMp \subset M_i$. Taking $q''=vv^*$ and
$u\in M$ a unitary element such that $uqq''=v$, the statement
follows.

$2^\circ$. Let $z=z(q)$ denote the central support of $q$ in $Q$
and note that $zq''$ is then the central support of $qq''$ in
$Qq''$. By the factoriality of $M_i$, $i=1,2$, it follows that
there exists a unitary $u\in \Cal U(M)$ such that $Q q''z \subset
uM_iu^*$. (Indeed, this is because whenever $Q_0\subset N$ is an
inclusion of finite von Neumann algebras,  $q_0 \in \Cal P(Q_0)$
and $N_0\subset N$ is a subfactor with $qQ_0q\subset N_0$, then
there exists $u\in \Cal U(N)$ such that $Q_0z(q)\subset uN_0u^*$.)
Thus, the projection $p_0'=q''z \in Q'\cap M $ together with the
unitary element $u$ satisfy the condition $uQp_0'u^* \subset M_i$.

Let $\Cal F$ be the set of all families of mutually orthogonal
projections $\{p_i'\}_{i\in I}\subset \Cal P(Q'\cap M)$, with the
property that $\forall i\in I, \exists j(i) \in \{1,2\}$ (unique by
Theorem 1.2.1)  and $v_i \in\Cal U(M)$ such that
$v_iQ{p_i}'v_i^*\subset M_{j(i)}$. $\Cal F$ is clearly inductively
ordered with respect to the order given by inclusion. Let
$\{p_i'\}_{i\in I}$ be a maximal element. Let $q_1'=\Sigma_{ j(i)=1}
p_i'$, $q_2'=\Sigma_{ j(i)=2} p_i'$ and set $q'= 1 - q_1' - q_2'$.

Assume $q'\neq 0$. Since $Q$ is diffuse, there exists $q\in \Cal
P(Q)$ such that $\tau(q'q)=1/n$ for some integer $n \geq 1$. Let
$\tilde{Q} \subset M$ be a von Neumann algebra isomorphic to $M_{n
\times n} (qQqq')$ with $qQqq'$ equal to the upper-left corner
$qq'\tilde{Q}qq'$ and $qq'$ having central trace $1/n$ in
$\tilde{Q}$. By (Sec. 4 in [P5]) $\tilde{Q} \subset M$ is a rigid
inclusion. Thus, we can apply the first part of the proof to get
$i\in \{1,2\}$, $ 0 \neq \tilde{q}' \in \tilde{Q}'\cap M$ and a
unitary element $w \in M$ such that $w \tilde{Q} \tilde{q}' w^*
\subset M_i$. Since $qq'$ has scalar central trace in $\tilde{Q}$,
it follows that the projection $p=qq' \tilde{q}'\in (qQq)'\cap
qMq=q(Q'\cap M)q$ is non-zero. Thus $p=qq''$, for some projection
$q'' \in Q'\cap M$ with $q''\leq q'$. Since $p \leq \tilde{q}'$,
we also have $w (qQqq'') w^* \subset M_i$, implying that if we
denote $p'= z(q)q''$ (where $z(q)$ is the central support of $q$
in $Q$) then $p' \in Q'\cap M$, $p'\leq q'$ and there exists a
unitary element $u$ in $M$ such that $u(Qp')u^* \subset M_i$. But
then $\{p_i'\}_{i \in I} \cup \{p'\}$ lies in  $\Cal F$, thus
contradicting the maximality of $\{p'_i\}_{i\in I}$.

We have thus shown that $q_1'+q_2'=1$. On the other hand, by the
factoriality of the $M_k$'s, $k = 1,2$, for each fixed $k$ we can
choose the unitary elements $\{v_i \mid j(i)=k\}$ which satisfy
$v_i (Qp'_i) v_i^* \subset M_k$ so that $v_ip'_iv_i^*$ be mutually
orthogonal projections in $M_k$. Taking $u_k\in \Cal U(M)$ to be a
unitary element extending $\Sigma_{j(i)=k} v_i$, it follows that
$u_k(Qq_k')u_k^*\subset M_k$, $k=1,2$.

This proves the existence part of part $2^\circ$. But the
uniqueness part is then clear, since if $p_1', p_2'$ is another
pair of projections in $Q'\cap M$ satisfying $p_1'+ p_2'=1$,
$v_i(Qp_i')v_i^* \subset M_i$ for some $v_i \in \Cal U(M)$,
$i=1,2$,  and we assume $x=p_1' q_2'\neq 0$, then the partial
isometry $w$ in the polar decomposition of $x$ lies in $Q'\cap M$
and if we denote $p=ww^*$ then $v_1(Qp)v_1^* \subset M_1$ while
$u_2w^*(Qp)wu_2^* \subset M_2$, contradicting Theorem 1.1.

To finish the proof of $2^\circ$ we need to show that $q_1', q_2'$
are in the center of $Q'\cap M$. Since $q_1'+q_2'=1$, this amounts
to showing that their central supports in $Q'\cap M$ are disjoint.

Assume by contradiction that there exist non-zero projections
$q_i''\leq q_i', q_i''\in Q'\cap M$ with $u' q_1'' {u'}^* = q_2''$
for some $u'\in \Cal U(Q'\cap M)$. But then
$u_k(Qq_k'')u_k^*\subset M_k$, $k =1,2$ are diffuse and are
conjugate by the unitary element $u_2 u'u_1^*$, contradicting
Theorem 1.1 again. \hfill $\square$

\proclaim{5.2. Theorem} Let $I$ be a set of indices with $0\in I$
and $(M_i, \tau_i)$, $i\in I$, a family of finite von Neumann
algebras with a common von Neumann subalgebra $B \subset M_i$,
such that $\tau_{0|B}=\tau_{i|B}, \forall i$. Assume $M_i$ are
factors for $i\neq 0$ and the inclusions $B \subset M_i$ are
homogeneous, $\forall i\in I$. Denote
$M=*_{B,i\in I} M_i$ the free
product with amalgamation over $B$ of the algebras $M_i, i\in I$. Let
$t > 0$ and $Q\subset M^t$ be a relatively rigid diffuse von
Neumann subalgebra such that no corner of $Q$ can be embedded into
$M_0$ inside $M$ and such that the normalizer of $Q$ in $M^t$
generates a factor $N$. Then there exists a unique $i\in I
\setminus \{0\}$ and a unitary element $u\in M^t$ such that $uQu^*
\subset M^t_i$. Moreover, such $u$ satisfies $uNu^*\subset M^t_i$,
and in fact $u\tilde{N}u^* \subset M^t_i$, where
$\tilde{N}=\tilde{N}(N,M_1^t;B)$ is as defined in $1.2.2$.
\endproclaim
\noindent {\it Proof}. Note first that $Q\subset M^t$ rigid
implies $Q$ is countably generated (see e.g. [PeP]). Thus, there
exists a countable subset $S\ni 0$ of indices $i\in I$ such that
$Q\subset (*_{B,i\in S} M_i)^t$.
By 3.2, $Q \subset (*_{B,i\in S} M_i)^t$
is rigid and by 1.2.1 all $N$ is contained in $(*_{B,i\in S} M_i)^t$.
This shows that it is sufficient to prove the
statement in the case $M_i, i\geq 0$, is a sequence of algebras.

Moreover, since $Q \subset M^t$ is rigid and since the factors
$\tilde{M}(K,t)=(M_0 *_B (*_{B,k\in K} M_k)^t$,
with $K$ finite
subset of $\{1,2, ...\}$, tend to $M^t$, it follows by ([P5]) that
there exists a non-zero projection $q' \in Q'\cap M^t$, a unitary
element $v\in M^t$ and a finite set $K\subset \{1,2,...\}$ such
that $v(Qq')v^* \subset \tilde{M}(K,t)$. But $q'\in N$ and by (3.4
in [P6]) $Qq'$ is quasi-regular in $q'Nq'$, so by 1.2.1 we have
$v(q'Nq')v^* \subset \tilde{M}(K,t)$. Since $N$ is a factor, $v$
can be modified so that $vNv^* \subset \tilde{M}(K,t)$. In
particular $vQv^* \subset \tilde{M}(K,t)$.

Since $Q$ is diffuse, there exists $q\in \Cal P(Q)$ such that
$\tau(q) \leq t^{-1}$. Thus, we may assume $v(qQq)v^* \subset
\tilde{M}_K \overset \text{\rm def} \to =  M_0
*_B (*_{B,k\in K} M_k)$, and notice that by (4.7 in [P5]) the unital
inclusion $v(qQq)v^* \subset p\tilde{M}_Kp$ is rigid, where
$p=vqv^*$. Since $K$ is finite, 5.1 applies to get $i\in K$,
$0\neq q_i' \in (vqQqv^*)' \cap p\tilde{M}_Kp$ and a unitary $w\in
\tilde{M}_K$ such that $w(vqQqv^*q_i')w^* \subset M_i$. Moreover,
since $\tau(q) \leq t^{-1}$ we can view $w(vqQqv^*q_i')w^*$ as a
(possibly non-unital) subalgebra of $M_i^t$. Since $q_i'\in
(vqQqv^*)'\cap pMp$, it follows that $q_i'\in vqNqv^*$ (recall
that $N$ is generated by the normalizer of $Q$ in $M$). By (3.4 in
[P6]) and 1.2.1 again, it follows that $w(q_i'vqNqvq_i')w^*
\subset M_i$, implying that $wq_i'vq$ can be extended to a unitary
element $u\in M^t$ such that $uNu^* \subset M^t_i$, thus $uQu^*
\subset M^t_i$. Also, by 1.2.1 $i$ is unique with this property
while by Corollary 1.2.3 it follows that $u\tilde{N}u^* \subset
M^t_i$. \hfill $\square$

\heading 6. Amalgamation over $\Bbb C$: Free product factors with
prescribed $\mycal F(M)$.
\endheading

We first apply Theorem 5.1 to plain free product factors, where
the result becomes an analogue of the classical Kurosh theorem for
groups. The first Kurosh-type results in Operator Algebra
framework were obtained by N. Ozawa in ([O2]). He proved that if
$N$ is a non-prime non-hyperfinite II$_1$ subfactor of a free
product $M=M_1 * M_2$ of semiexact finite factors  $M_1, M_2$,
then $N$ can be unitarily conjugated into either $M_1$ or $M_2$
(this is analogue to ``Kurosh subgroup theorem''). As a
consequence, he showed that if two free products $*_i M_i$, $*_j
N_j$ of non-hyperfinite, non-prime, semiexact factors $N_i, M_j$
are isomorphic then the ``length'' of the two free products must
be the same and each $N_i$ is unitary conjugate to $M_i$, after
some permutation of indices (this is analogue to ``Kurosh
isomorphism theorem'').

In turn, our results cover different classes of algebras. Thus,
our analogue of the ``Kurosh subgroup theorem'' allows $M_1, M_2$
to be arbitrary finite von Neumann algebras but only gives
information about relatively rigid subalgebras $Q$ of $M_1*M_2$.
Our corresponding ``isomorphism theorem'', which in fact we obtain
for amplifications of free products,  will require the factors
$N_i, M_j$ to be either w-rigid, i.e. to have diffuse, regular,
relatively rigid subalgebras, or to be group measure space factors
associated to actions of w-rigid ICC groups. In particular, it
holds true for $N_i, M_j$ II$_1$ factors with the property (T) (in
the sense of [CJ]) and more generally for tensor products of
property (T) II$_1$ factors with arbitrary finite factors.
Moreover, since the factors $N=L(\Bbb Z^2 \rtimes \Bbb F_n)$ in
([P5]) are w-rigid and have trivial fundamental group (see Theorem
7.11 for a different proof), this will allow us to obtain large
classes of factors with trivial fundamental group, different from
the ones in ([P5,6]). More generally, using also ([DyR]), we
construct a completely new class of factors with prescribed
fundamental group which, unlike the ones in ([P5,6]), have no
Cartan subalgebras (by [V2], cf. 6.6 below).

\proclaim{6.1. Theorem} Let $(M_i, \tau_i)$, $i=0,1,2$ be finite
von Neumann algebras and denote $M=M_0 * M_1 *M_2$. Assume no
direct summand of $(M_0, \tau_0)$ has relatively rigid diffuse von
Neumann subalgebras $($e.g. $M_0=\Bbb C$, or more generally $M_0$
with Haagerup property$)$. Let $Q\subset M$ be a relatively rigid
diffuse von Neumann subalgebra of $M$.

$1^\circ$. There exist $i \in \{1,2\}$, $q\in \Cal P(Q)$, $q'\in
\Cal P(Q'\cap M)$ and $u\in \Cal U(M)$ such that $qq'\neq 0$ and
$uqQqq'u^*\subset M_i$.

$2^\circ$. If in addition $M_1, M_2$ are factors,  then there
exists a unique pair of projections $q_1', q_2'\in Q'\cap M$ such
that $q_1'+q_2'=1$ and $u_i(Qq_i')u_i^*\subset M_i$ for some
unitaries $u_i\in\Cal U(M)$, $i=1,2$. Moreover, $q_1', q_2' \in
\Cal Z(Q'\cap M)$.

$3^\circ$. If instead of  $M_1, M_2$ we consider a whole family of
finite factors $M_i, i\geq 1$, and we take $Q\subset M^t=(M_0
* M_1 * M_2*...)^t$ a rigid inclusion, for some $t > 0$, and we
assume the normalizer of $Q$ in $M^t$ generates a factor $N$, then
there exists a unique $i \geq 1$ and a unitary element $u\in M^t$
such that $uQu^* \subset M^t_i$. Moreover, such $u$ satisfies
$uNu^*\subset M^t_i$, and in fact $u\tilde{N}u^* \subset M^t_i$,
where $\tilde{N}=\tilde{N}(N,M_i^t;\Bbb C)$ is as defined in
$1.2.2$.
\endproclaim
\noindent {\it Proof.} As in the proof of 5.2, note that $Q$
relatively rigid implies $Q$ countably generated. Thus, there
exist countably generated von Neumann subalgebras $M_i^0\subset
M_i, i=0,1,2,$ such that $Q \subset M^0_0*M^0_1*M^0_2$. Thus, to
prove 1$^\circ$ it is clearly sufficient to prove it in the case
$M_i$ are countable generated, $i=0,1,2$. But then  each $\Bbb C
\subset M_i$ is homogeneous by Lemma 4.2. Let us show that no
corner of $Q$ can be embedded into $M_0$ inside $M$. Assume this
is not true. By (2.1. in [Po6]) it follows that there exist $q\in
Q, p\in M_0$ non-zero projections, a unital isomorphism $\psi$ of
$qQq$ into $pM_0p$ and a non-zero partial isometry $v\in M$ such
that $vv^*\in (qQq)'\cap qMq,v^*v\in \psi(qQq)'\cap pMp$ and
$xv=v\psi(x),\forall x\in qQq$. Denote $q''=vv^* \in Q'\cap M$.

Since $\psi(qQq)\subset pM_0p$ is a diffuse von Neumann
subalgebra, by Theorem 1.2.1, it follows that $\psi(qQq)'\cap pMp
\subset pM_0p$. Thus $v^*v\in pM_0p$. This shows that
$v^*qQqv\subset pM_0p$, which in turn implies that $qQqq''\subset
wM_0w^*$, for some unitary element $w\in M$ extending $v$. Since
$Q \subset M$ is rigid, by (Proposition 4.7 in [Po5])
$qQqq''\subset qq''Mqq''$ is also rigid, which trivially implies
$qQqq'' \oplus (1- qq'')\Bbb C \subset M$ is rigid. But then
$w^*(qQqq'' \oplus (1- qq'')\Bbb C)w \subset M_0$ follows rigid by
Corollary 3.2. By taking a suitable amplification of
$w^*(qQqq'')w$ in $M_0$ and using again (4.7 in [Po5]), this
implies a direct summand of $M_0$ contains a relatively rigid
diffuse von Neumann subalgebra, a contradiction.

Altogether, this shows that the conditions required in  5.1 are
satisfied, so part $1^\circ$ of the statement follows as a
particular case of that theorem.

For part $2^\circ$, simply notice that if $M_1, M_2$ are factors
then the countably generated von Neumann subalgebras $M_0^0,
M_1^0, M_2^0$ with the property that $Q \subset M_0^0 * M_1^0 *
M_2^0$ can be chosen so that $M_1^0, M_2^0$ are factors as well,
so 5.1.2$^\circ$ applies.

Part $3^\circ$ follows then from part $2^\circ$ and 5.2. \hfill
$\square$

\vskip .05in \noindent {\it 6.2. Definition}. A finite von Neumann
algebra $(M,\tau)$ is {\it weakly rigid} (w-rigid) if it contains
a regular, relatively rigid, diffuse von Neumann subalgebra, i.e.
a subalgebra $Q \subset M$ such that $\Cal N_M(Q)''=M$ and $Q
\subset M$ is a rigid inclusion (or $Q$ is a relatively rigid
subalgebra of $M$ [P5]). Note that if $G$ is a w-rigid group as
defined in ([P5-8]), i.e. $G$ contains an infinite normal subgroup
with the relative property (T) of Kazhdan-Margulis, then $L(G)$ is
w-rigid. Also, if $M$ is w-rigid and $P$ is an arbitrary finite
von Neumann algebra then $M \overline{\otimes} P$ is w-rigid.

\proclaim{6.3. Theorem} Let $(M_0, \tau_{M_0}), (N_0, \tau_{N_0})$
be finite von Neumann algebras which have no relatively rigid
diffuse subalgebras $($e.g. $M_0$, $N_0$ have Haagerup's compact
approximation property$)$. Let $M_1, \ldots, M_m$ and $N_1,
\ldots, N_n$ be $\text{\rm II}_1$ factors, where $n, m \geq 1$ are
some cardinals $($finite or infinite$)$ and assume each $M_i, N_j$
is w-rigid. If $\theta$ is an isomorphism of $M=*_{i=0}^{m}M_i$
onto $N^t$, where $N=*_{j=0}^{n} N_j$ and $t > 0$, then $m=n$ and,
after some permutation of indices, $\theta(M_i)$ and $N^t_i$ are
unitary conjugate in $N^t$, $\forall i\geq 1$.
\endproclaim
\vskip .05in \noindent {\it Proof}. For each $m \geq i \geq 1$ let
$Q_i \subset M_i$ be a regular, relatively rigid, diffuse von
Neumann subalgebra. Since $M_i$ are factors and $Q_i\subset M_i$
are regular inclusions, it follows that for any $s >0$ the factor
$M^s_i$ contains a regular, diffuse, relatively rigid subalgebra.
To see this, note that $\Cal N_{M_i}(Q_i)$ acts ergodically on the
center of $Q_i$, so  $\Cal Z(Q_i)$ is either diffuse or atomic. In
both cases we can find projections $q'\in \Cal Z(Q_i)$ and $q\in
Q_iq'$ such that the central trace of $q$ in $Q_i$ is a scalar
mutiple of $q'$ and $\tau(q)=s/k$, for some $k \geq t$. By (Lemma
3.5 in [P6I]) it follows that the inclusion $qQ_iq\subset qM_iq$
is regular and by (4.7 on [P5]) it is rigid as well. But then
$M_{k\times k}(qQ_iq) \subset M_{k \times k}(qM_iq)=M_i^s$ is
regular and rigid ([P5]).

Moreover, note that any diffuse von Neumann subalgebra $B_i
\subset M_i^s$ satisfies $B_i'\cap M^s \subset M_i^s$. To see
this, note first that by taking direct sums of $k$ copies of $B_i$
embedded diagonally into $M_{k \times k}(M_i^s)=M_i^{ks}$, with
$k$ sufficiently large, we may assume $s\geq 1$. Then take $p\in
\Cal P(B_i)$ with $\tau(p)=1/s$ and note that if we assume by
contradiction that $B_i'\cap M^s \neq B_i'\cap M^s_i$, then
$(pB_ip)'\cap pM^sp\neq (pB_ip)'\cap pM^s_ip$. But $(pM_i^sp
\subset pMp)\simeq (M_i\subset M)$ and $M$ splits off $M_i$ as a
free product. Thus, by Theorem 1.2.1 the relative commutant in
$pM^sp=M$ of the diffuse subalgebra $pB_ip\subset pM_i^sp=M_i$
must be contained in $pM^s_ip$, contradicting the assumption.

Taking now $s=1/t$, it follows that $M_i^{1/t}$ has a diffuse,
regular, relatively rigid subalgebra, implying that
$P_i=\theta^{1/t}(M_i^{1/t})$ has such a subalgebra $B_i$ as well.
In particular, the inclusion $B_i\subset N$ is rigid. In addition,
$B_i'\cap N \subset P_i$. Since the inclusion $B_i\subset N$ is
rigid and regular, by 6.1.3$^\circ$ there exists a unique $j(i)\in
\{1,2,..., n\}$ and a unitary $u_i \in N$ such that $B_i \subset
u_i N_{j(i)} u_i^*$ and $P_i \subset u_i N_{j(i)} u_i^*$. Thus,
there exists a unique $j(i)$ such that for some unitary $v_i \in
N^t$ we have $\theta(M_i)=P^t_i \subset v_i N^t_{j(i)} v_i^*$.

Similarly, by applying the above to $\theta^{-1}$ we get for each
$n \geq j \geq 1$ a unique $m \geq k(j)\geq 1$ and a unitary
element in $w_j \in M$ such that $\theta^{-1}(N^t_j) \subset w_j
M_{k(j)} w_j^*$. Altogether, for each $m \geq i \geq 1$ we get
$$
M_i = \theta^{-1}(\theta(M_i)) \subset \theta^{-1} (v_i N^t_{j(i)}
v_i^*) \subset u_i M_{k(j(i))} u_i^*, \tag 6.3.1
$$
where $u_i = w_{j(i)} \theta^{-1}(v_i)$. By Theorem 1.2.1 it
follows that $k(j(i))=i$, i.e. $k \circ j =id$. Similarly, $j
\circ k=id$. Thus $m=n$, $j$ and $k$ are onto isomorphisms and the
inclusions $(6.3.1)$ are in fact equalities. \hfill $\square$

\vskip .05in

In the next statement, for a finite permutation $\pi \in S_m$ and
$m \geq i \geq 1$ we denote by $m(\pi,i)$ the cardinality of the
set $\{\pi^k(i) \mid k \geq 1 \}$.

\proclaim{6.4. Corollary} Let $m \in \Bbb N$ and let $M_1, \ldots,
M_m$ be w-rigid $\text{\rm II}_1$ factors. Let $(M_0, \tau_{M_0})$
be a finite von Neuman algebra which contains no rigid diffuse von
Neumann subalgebras. Denote $M = *_{i = 0}^m M_i$. Then $\mycal
F(M) \subset \cup_{\pi \in S_m} \cap_{i = 1}^m \mycal
F(M_i)^{m(\pi,i)^{-1}} \subset \cap_{i = 1}^m \Cal F(M_i)^{1/m!}$.
In particular if one of the factors $M_i, m\geq i \geq 1,$ has
trivial fundamental group then so does $M$.
\endproclaim
\noindent {\it Proof}. For $t\in\mycal F(M)$, let $\theta:
M\rightarrow M^t$ be an isomorphism. Applying the previous
theorem, we get that there exists $\pi\in S_m$ such that
$\theta(M_i)$ and $M^{t}_{\pi(i)}$ are unitary conjugate in $M^t$.
In particular, we have that $M_i\cong M^{t}_{\pi(i)},\forall 1\leq
i\leq m$ which by induction implies that $M_i\cong
M^{t^{k}}_{\pi^{k}(i)},\forall 1\leq i\leq m,k\in \Bbb N $.

Fixing $i$ and letting $k=m(\pi,i)$ we obtain that
$t^{m(\pi,i)}\in\mycal F(M_i)$,
 or equivalently $t\in\mycal F(M_i)^{m(\pi,i)^{-1}}$.
By intersecting over all values of $i$, taking the union over all
possible permutations and noticing that $m(\pi,i)$ divides $m$, the
result follows. \hfill $\square$

\vskip .05in

For the next corollary we denote by $\Cal X_{fin}$ the set of all
finite tuples of positive numbers $(t_i)_{1 \leq i \leq n} \subset
\Bbb R_+^*$, $n \geq 2$, and by $\Cal X_\infty$ the set of all
infinite sequences $(t_j)_{j \geq 1} \subset \Bbb R_+^*$. Also, we
let $\Cal X = \Cal X_{fin} \cup \Cal X_\infty$. If $X=(t_i)_i$,
$Y=(s_j)_j$ are in $\Cal X$ we then write $X \sim Y$ if both have
the same ``length'' and there exists a permutation (bijection) $\pi$
of the (common) set of indices $\{1,2, ... \}$ such that $s_{\pi(i)}
= t_i$, $\forall i $.

Given a II$_1$ factor $M$ and $X=(t_i)_i\in \Cal X$, we denote
$M^X=*_i M^{t_i}$. Note that if $X, Y\in \Cal X$ and $X \sim_\pi
Y$ then $\pi$ implements a natural isomorphism $\theta_\pi : M^X
\simeq M^Y$, in the obvious way. For $t > 0$ and $X=(t_i)_i \in
\Cal X$ we denote $tX=(tt_i)_i \in \Cal X$.

\proclaim{6.5. Corollary} Let $M$ be a w-rigid II$_1$ factor with
trivial fundamental group, e.g. $M=L(\Bbb Z^2 \rtimes SL(2,\Bbb
Z))$ $($cf. $\text{\rm [P5]})$.

$1^\circ$. If $X, Y \in \Cal X$ then $M^X \simeq M^Y$ iff $X \sim
Y$ and iff $M^X * L(\Bbb F_k) \simeq M^Y * L(\Bbb F_k)$, for some
$1 \leq k \leq \infty$.

$2^{\circ}$.  $\mycal F(M^X)=\{1\}$ for all $X \in \Cal X_{fin}$.
Moreover, if we denote by $\Cal X_0$ the set of elements $X$ in
$\Cal X_{fin}$ with $\min X = 1$, then $\{M^X \mid X \in \Cal
X_0\}$ is a continuous family of mutually non stably isomorphic
$\text{\rm II}_1$ factors.

$3^\circ$. For each $X \in \Cal X_\infty$ denote $\Cal S_X = \{ t
\in \Bbb R_+^* \mid tX \sim X \}$. Then $\mycal F(M^X) = \Cal
S_X$. In particular, if $S \subset \Bbb R_+^*$ is an infinite,
countable subgroup and $X \in \Cal X_\infty$ has the elements of
$S$ as entries, each one repeated with the same $($possibly
infinite$)$ multiplicity, then $\mycal F(M^{tX})=S$, $\forall t >
0$. Moreover, $M^{t_1X}, M^{t_2X}$ are stably isomorphic iff $t_1,
t_2$ are in the same class in $\Bbb R_+^*/S$. Thus, $\{ M^{tX}
\mid t \in \Bbb R_+^*/S \}$ is a continuous family of mutually non
stably isomorphic $\text{\rm II}_1$ factors all with fundamental
group equal to $S$.

$4^\circ$. If $S \subset \Bbb R_+^*$ is an arbitrary infinite
$($possibly uncountable$)$ subgroup then the $\text{\rm II}_1$
factor $M^S=*_{s\in S} M^s$ has fundamental group equal to $S$.
\endproclaim
\vskip .05in \noindent {\it Proof}. $1^\circ$. If $M^X \cong M^Y$
then by Theorem 6.2 $X,Y$ have the same length and there exists a
bijection $\pi$ of the corresponding (finite or infinite) set of
indices $\{1,2,...\}$ such that $M^{t_i} = M^{s_{\pi(i)}}$,
$\forall i$. Since $M$ has trivial fundamental group, this implies
$t_i = s_{\pi(i)}, \forall i,$ and thus $X \sim_\pi Y$. The second
equivalence has exactly the same proof, using 6.3 with $M_0 =
L(\Bbb F_k)$.

$2^{\circ}$. If $X=(t_i)_{1\leq i \leq n}$, $Y = (s_j)_{1\leq j
\leq m}$ are in $\Cal X_{fin}$ and $M^X \simeq (M^Y)^t$ then $M^X
* L(\Bbb F_\infty)\simeq (M^{Y})^t * L(\Bbb F_\infty)$. But by
([DyR]) the last factor is isomorphic to $M^{tY} * L(\Bbb
F_\infty)$. Thus, $M^X
* L(\Bbb F_\infty) \simeq M^{tY} * L(\Bbb F_\infty)$, which by
part $1^\circ$ implies $X \sim tY$. Thus, if $X=Y$ (respectively
$X,Y\in \Cal X_0$) then $t=1$ and we get $\mycal F(M^X)=\{1\}$
(resp. $X \sim Y$). This implies both statements.

$3^\circ$. By ([DyR]), if $X \in \Cal X_\infty$ then $(M^X)^t
\cong M^{(tX)}$. Thus $M^X \cong (M^X)^t$ if and only if $X \sim
tX$, which readily implies all statements.

$4^\circ$. It is easy to see that the proof of the amplification
formula $(M^S)^t\simeq M^{tS}$ in ([DyR]) doesn't in fact depend
on the fact that the infinite set $S$ is countable. Thus, since
for $t\in S$ we have $tS = S$ as sets, it follows that
$(M^S)^t=M^{tS}=M^S$, thus $S \subset \mycal F(M^S)$. Conversely,
if $s\in \Bbb R_+^*$ satisfies $(M^S)^s\simeq M^S$, then
$M^{sS}\simeq M^S$, which by 6.3 implies $tS=S$, so that $s\in S$.
\hfill $\square$

\vskip .05in

\noindent {\bf 6.6. Remarks}. Dima Shlyakhtenko pointed out to us
that by combining Voiculescu's initial argument for showing that
$L(\Bbb F_n)$ has no Cartan subalgebras, with Kenley Jung's
``monotonicity'' ([Ju]) it follows that any free product of type
II$_1$ factors is ``Cartan-less'' (see [Sh] for a detailed
argument). Thus, unlike the examples of factors with prescribed
fundametal group in ([P6]), which are group measure space factors
associated with equivalence relations coming from Connes-St\o rmer
Bernoulli actions, the examples of factors $M$ that we produce
here (in Corollary 6.5) have no Cartan subalgebras, and altogether
no diffuse hyperfinite ``core''. In particular, they cannot be
written as crossed products of the form $M=R \rtimes \Gamma$ with
$R$ the hyperfinite factor.

It is interesting to note that Theorem 5.1 can be used to give a
completely new proof of the by now classical result of Connes and
Jones showing that property (T) factors cannot be embedded into
the free group factor ([CJ]). Thus, rather than using Haagerup's
property (i.e. ``deformation by compact c.p. maps''), as the
original proof does, this new proof uses a ``deformation by
automorphisms'' of  the free group factors.

\proclaim{6.7. Corollary} $(\text{\rm [CJ]})$ For every n, $2\leq
n \leq\infty$ the free group von Neumann algebra $L(\Bbb F_n)$
contains no relatively rigid diffuse subalgebra.
\endproclaim
\vskip .05in \noindent {\it Proof}. If we write $L(\Bbb F_n)$ as
$L(\Bbb Z) * L(\Bbb Z) * ... * L(\Bbb Z)$ and then  apply
recursively the first part of Theorem 5.1 and Corollary 3.2, it
follows that a corner of $L(\Bbb Z)$  contains a rigid diffuse von
Neumann subalgebra, a contradiction. \hfill $\square$

\heading 7. Amalgamation over Cartan subalgebras: vNE/OE rigidity
results
\endheading

In this section we apply Theorem 5.1 to study group measure space
factors of the form $A \rtimes_\sigma \Gamma$ where $\Gamma$ is a
free product of groups $\Gamma=\Gamma_0 * \Gamma_1 *...$ and
$\sigma$ is a free ergodic m.p. action of $\Gamma$ on $A =
L^\infty(X,\mu)$, for $(X,\mu)$ a probability space. Such a factor
can alternatively be viewed as a free product with amalgamation
$M=M_0 *_A M_1 *_A ...$, where $M_i =A \rtimes_{\sigma_i}
\Gamma_i$, $\sigma_i=\sigma_{|\Gamma_i}$, with $A \subset M$, with
the algebra of coefficients $A$ of the crossed product $A \rtimes
\Gamma$ now becoming the  ``core'' of the amalgamated free
product. It is this form that will allow us to use 5.1.

Following ([P6]), we regard an isomorhism of such group measure
space factors as a {\it von Neumann equivalence} (vNE) of the
corresponding actions $(\sigma, *_i \Gamma_i)$. Thus, the main
result we prove in this section is a rigidity result showing that
vNE of actions of free products of groups $\Gamma_i$ satisfying
some weak rigidity conditions (of property (T) type) entails the
orbit equivalence (OE) of the actions $\sigma$, with component by
component OE of the actions $(\sigma_i, \Gamma_i)$. Due to its
analogy to similar statements on (amalgameted) free products of
groups in Bass-Serre theory, we refer to this as vNE {\it
Bass-Serre rigidity}. We note that when applied to isomorphisms of
group measure space factors that come from OE of the actions, they
give OE {\it Bass-Serre rigidity} results.

Since we study the group measure space factors $M= A \rtimes
(\Gamma_0 * \Gamma_1 *...)$ as AFP factors $(A \rtimes \Gamma_0)
*_A (A \rtimes \Gamma_1) *_A...$, it is worth noticing that if
$M=M_0 *_A M_1 *_A...$ is an AFP factor coming from Cartan
subalgebra inclusions $A=L^\infty(X, \mu) \subset M_i$, $i\geq 0$,
then the AFP ``core'' $A$ follows regular in $M$, but in general
it may not be maximal abelian (thus not Cartan) in $M$. For
instance, if $A$ is a Cartan subalgebra of a II$_1$ factor $N$
then $A$ is not maximal abelian in $M=N*_A N$, because for any
$u\in \Cal N_N(A)$ with $E_A(u)=0$ the element $u * u^{-1}$ is
still perpendicular to $A$ yet acts trivially on it. For more on
general properties of AFP factors arising from Cartan inclusions
we send the reader to ([Ko], [U1,2]).

In case $M_i = A \rtimes_{\sigma_i} \Gamma_i$, with each $\sigma_i$
a free m.p. action, then there exists a unique m.p. action $\sigma$
of $\Gamma= \Gamma_0 * \Gamma_1 * ...$ on $A$ such that
$\sigma_{|\Gamma_i}=\sigma_i$, $\forall i$, and we still have the
natural identification $M=A \rtimes_\sigma \Gamma = M_0 *_A M_1*_A
...$, as in the case $\sigma$ a free action mentioned above. Then
$A$ is Cartan in $M$ iff $\sigma$ is a free action, i.e. iff
$\sigma_i$ are ``freely independent'' actions (in the obvious
sense). For general equivalence relations (or Cartan subalgebras),
the definition of ``free independence'' was formulated by Gaboriau
([G1]) and is reminded below. Recall that if $B$ is a finite von
Neumann algebra, $p,q$ are non-zero projections in $B$ and
$\theta:pBp \rightarrow qBq$ is a $^*$-morphism, then $\theta$ is
called {\it properly outer} if $b\in B$, $\theta(x)b=bx$, $\forall
x\in B$, implies $b=0$.

\vskip .05in \noindent {\bf 7.1. Definition} ([G1]). Let $\{\Cal
R_i\}_{i\in I}$ be a family of countable, measurable, measure
preserving equivalence relations on the same standard non atomic
probability space $(X,\mu)$ (see e.g. [FM]). We alternatively view
each $\Cal R_i$ as a pseudogroup of {\it local m.p. isomorphisms}
$\phi: Y_1 \simeq Y_2$ with $Y_1, Y_2 \subset X$ measurable and the
graph of $\phi$ contained in $\Cal R_i$ ([Dy], [FM]). We say that
$\{\Cal R_i\}_i$ are {\it freely independent} if for any $n$ and any
properly outer local isomorphisms $\phi_j \in \Cal R_{i_j}$, $1 \leq
j \leq n$, $i_j \in I$, with $i_j \neq i_{j+1}, 1\leq j \leq n-1$,
the product $\phi_1 \phi_2 ... \phi_n$ is properly outer.

\vskip .05in

In the case when each of the equivalence relations $\Cal R_i$ is
generated by properly outer automorphisms, Definition 7.1 can be
viewed as a particular case of the following:

\vskip .05in \noindent {\bf 7.1'. Definition}. Let $(B, \tau)$ be
a finite von Neumann algebra and $S_i \subset \text{\rm
Aut}(B,\tau), i\in I,$ a family of sets of $\tau$-preserving
automorphisms, with each $\theta \in S_i$ either properly outer or
equal to $id_B$. We say that $\{S_i\}_i$ are {\it freely
independent} if for any $n$ and any $\theta_j \in S_{i_j}\setminus
\{id_B\}$, $1 \leq j \leq n$, $i_j \in I$, with $i_1 \neq i_2\neq
... \neq i_n,$, the product $\theta_1 \theta_2 ... \theta_n$ is
properly outer.

The next lemma translates the freeness conditions in definition
7.1 into the framework of operator algebras (see [U2]):

\proclaim{7.2. Lemma} Let $(M_n, \tau_n), n \geq 1,$ be finite von
Neumann algebras with a common Cartan subalgebra $A \subset M_n$,
such that $\tau_{n|A}=\tau_{m|A}, \forall n,m$. Then $A \subset
M_1 *_A M_2 *_A...$ is a Cartan subalgebra if and only if the
equivalence relations $\Cal R_n=\Cal R_{A \subset M_n}$, $n \geq
1$, are freely independent.
\endproclaim
\noindent {\it Proof}. Since $A$ is clearly regular in $M=M_1 *_A
M_2 *_A...$, all we need to prove is that $A \subset M$ is maximal
abelian iff $\{\Cal R_n\}_n$ are freely independent. But this is
trivial by the definitions of freeness and respectively of the
amalgamated free product over $A$. \hfill $\square$

\vskip .05in The next result, essentially due to Tornquist ([To]),
shows that a sequence of actions of countable groups (or merely
countable m.p. equivalence relations) can be made ``freely
independent'' by conjugating each one of them with a suitable m.p.
automorphism. We include a proof, based on Lemma A.1 in the
Appendix, for the reader's convenience.

\proclaim{7.3. Proposition} $1^\circ$. Let $(X, \mu)$ be a standard
non atomic probability space and $\sigma_n : G_n \rightarrow
\text{\rm Aut}(X,\mu)$ be free m.p. actions of discrete countable
groups $G_n$, $n\geq 1$. Then there exists a free m.p. action
$\sigma$ of $G=*_n G_n$ on $(X,\mu)$ such that $\sigma_{|G_n}$ is
conjugate to $\sigma_n$, $\forall n\geq 1$. More generally, if
$\{\Cal R_n\}_n$ are standard equivalence relations on $(X,\mu)$
then there exists an equivalence relation $\Cal R$ on $(X,\mu)$
generated by a family of freely independent sub equivalence
relations $\Cal R'_n \subset \Cal R, n \geq 1,$ such that $\Cal
R_n \simeq \Cal R'_n, \forall n$.

$2^\circ$. Let $(M_n, \tau_n)$ be countably generated finite von
Neumann algebras with a common diffuse Cartan subalgebra $A
\subset M_n$, $n \geq 0$, such that $\tau_{n|A}=\tau_{m|A}$,
$\forall n,m$. Then there exist Cartan subalgebra inclusions $\{A
\subset N_n \}_n$ such that $(A \subset N_n) \simeq (A \subset
M_n)$, $\forall n,$ and such that $A$ is a Cartan subalgebra in
$N_0 *_A N_1 *_A N_2 ...$.
\endproclaim
\noindent {\it Proof}. $1^\circ$. This is an immediate application
by induction of Lemma
A.1, once we notice that any $\Cal R_n$ can be extended to a
countable m.p. equivalence relation $\Cal S_n$ on $(X,\mu)$ which
is generated by countably many properly outer m.p. automorphisms.

$2^\circ$. By Lemma 7.2, we only need to make the equivalence
relations $\Cal R_{A \subset M_n}$ freely independent, so the
first part applies. \hfill $\square$

\vskip .05in We also notice the following general ``compression
formula'' for restrictions of ``free products of equivalence
relations'' $\Cal R = *_{i=1}^n \Cal R_i$, i.e. for relations
$\Cal R$ that are generated by freely independent sub-equivalence
relations $\Cal R_i\subset \Cal R, 1\leq i \leq n$.

\proclaim{7.4. Proposition} $1^\circ$. Let $M_i, 1 \leq i \leq
n$, be $\text{\rm II}_1$ factors, for some $2\leq n\leq \infty$,
with a common Cartan
subalgebra $A$ and assume $A \subset M = M_1 *_A M_2 .. *_A M_n$
is Cartan. If $p\in A$ is a projection of trace $1/m$ for some
integer $m \geq 1$ then the Cartan subalgebra inclusion $Ap\subset
pMp$ is naturally isomorphic to $Ap \subset M_0 *_{Ap} pM_1p
*_{Ap} ...*_{Ap} pM_np $, where $(Ap \subset M_0)=(Ap \rtimes \Bbb
F_{(n-1)(m-1)})$ for some free action of the free group $\Bbb
F_{(n-1)(m-1)}$ on $Ap$.

$2^\circ$. Let $\Cal R_1, ..., \Cal R_n$ be freely independent,
countable, ergodic, m.p. equivalence relations on the same standard
probability space $(X,\mu)$ and denote by $\Cal R$ the equivalence
relation they generate. If $Y \subset X$ is a subset of measure
$1/m$ then the restriction $\Cal R^Y$ of $\Cal R$ to $Y$ is
generated by the freely independent ergodic sub-equivalence
relations $\Cal R^Y_i, 1 \leq i \leq n$, and $\Cal R_0$, where $\Cal
R_i^Y$ is the restriction of $\Cal R_i$ to $Y$ and $\Cal R_0$ is
generated by a free m.p. action of a free group with $(n-1)(m-1)$
generators $\Bbb F_{(n-1)(m-1)}$ on $Y$.
\endproclaim
\noindent {\it Proof}. It is clearly sufficient to prove part
$1^\circ$. By the representation of AFP algebras $(1.1.1)$, the
von Neumann algebra generated by $pM_ip, 1 \leq i \leq n$, in
$pMp$ is isomorphic to the AFP algebra $pM_1p *_{Ap} ... *_{Ap}
pM_np$. On the other hand,  since $\tau(p)=1/m$ and each $M_i$ is
a factor, by Dye's Theorem ([D]) there exist matrix units
$\{e_i^{jk}\}_{1\leq j,k \leq m}$ in the normalizing groupoid
$\Cal G\Cal N_{M_i}(A)$ of $A \in M_i$ such that $e^{11}_i=p$ and
$e^{jj}_i=e^{jj}_{i'}$, $\forall 1\leq i, i' \leq n$ and $\forall
1\leq j \leq m$. Denote by $u^j_i = e^{1j}_1 e^{j1}_i\in pMp$,
$2\leq i \leq n,$ $2 \leq j \leq m$, and notice that there are
$(n-1)(m-1)$ many such unitary elements. The expansion $(1.1.1)$
of $M_1 *_A ... *_A M_n$ implies that $\{u^j_i\}_{i,j}$ are the
generators of a free group $\Bbb F_{(n-1)(m-1)}$, all of whose
elements $\neq 1$ are perpendicular to $Ap$. Since $A \subset M$
is Cartan, this implies the action implemented by $\Bbb
F_{(n-1)(m-1)}$ on $Ap$ is free. Moreover, if we denote $M_0 = Ap
\vee \{u^j_i\}_{i,j}'' \simeq A_p \rtimes \Bbb F_{(n-1)(m-1)}$
then it is immediate to check that if $v=v_{i_1} v
_{i_2} ... v_{i_k}\in pMp$ is an ``alteranting word'',
with $v_{i_l}$ in the normalizing groupoid of
$Ap$ in $pM_{i_l}p, \forall l,$ and $i_1 \neq i_2 \neq ... \neq
i_k$ in $\{0,1, ... ,n\}$, $E_{Ap}(v_{i_l})=0, \forall l$,
then $v$ has expectation $0$ on $Ap$ as well, $E_{Ap}(v)=0$. Thus, if
we denote by $N\subset pMp$ the von Neumann algebra generated by
$pM_ip, 0 \leq i \leq n,$ then $(Ap \subset N) = (Ap\subset M_0
*_{Ap} pM_1p *_{Ap} ... *_{Ap} pM_np)$.

Finally, since $e^{j1}_1 u^j_i=e^{j1}_i$, we have that $M$ is
generated by $N=pNp$ and the matrix unit $\{e_1^{jk}\}_{j,k}$.
This also implies that $pMp$ is generated by $\vee_{i=0}^n pM_ip$.
Altogether, this shows that $(Ap \subset pMp)=(Ap\subset
N)=(Ap\subset M_0 *_{Ap} pM_1p *_{Ap} ... *_{Ap} pM_np)$. \hfill
$\square$

\vskip .05in The above result shows in particular that if $\Bbb
F_n \curvearrowright X$ is a free m.p. action on the probability
space with restrictions to each of the generators of   $\Bbb F_n$
acting ergodically on $X$ then the amplification of the
corresponding orbit equivalence relation $\Cal R_{\Bbb F_n}$ by
$1/m$ is an equivalence relation that can be implemented by a free
ergodic m.p. action of a free group with $m(n-1)+1$ generators,
thus recovering the ``compression formula'' $(\Cal R_{\Bbb
F_n})^{1/m}=\Cal R_{\Bbb F_{m(n-1) +1 }}$ in ([Hj]). On the other
hand, Proposition 7.4 can be viewed as an AFP version of the
``compression formula'' for plain free product factors in ([DyR]).

The statement $7.5$ below will require the following:

\vskip .05in \noindent {\it 7.5.0. Notation}. Let
$\{\Gamma_{i,j}\}_{j=0}^{n_i}$ be discrete countable groups, $1
\leq n_i \leq \infty$, $i=1,2$. Denote $G_i = \Gamma_{i,0} *
\Gamma_{i,1}*... * \Gamma_{i,n_i}$, $i=1,2$. Let $\sigma_i:G_i
\rightarrow \text{\rm Aut}(X_i,\mu_i)$ be a free, m.p. action on a
standard probability space $(X_i,\mu_i)$. Denote
$A_i=L^\infty(X_i,\mu_i)$, $M_i=A_i \rtimes_{\sigma_i} G_i$ and
$M_{ij}=A_i \rtimes_{\sigma_{ij}}\Gamma_{ij}$, where
$\sigma_{ij}=(\sigma_i)_{|\Gamma_{ij}}$, $\forall n_i \geq j \geq
0, i=1,2$.

\vskip .05in

The general result we prove shows that, under suitable weak
rigidity conditions on the groups $\Gamma_{ij}$, an isomorphism
between the factors $M_1, M_2$ must take each of the ``component
inclusions'' $(A_i \subset M_{ij})$ onto each other, modulo some
permutation of indices and unitary conjugacy. Since the weak
rigidity assumption on the $\Gamma_{ij}$'s is somewhat technical,
we display the conditions separately and give right away a list of
examples when they are satisfied:

\vskip .1in \noindent {\it 7.5.1. Assumption}. $\Gamma_{10},
\Gamma_{20}$ have Haagerup property, $i=1,2$, and if both
$\Gamma_{10}, \Gamma_{20}$ are finite then at least one of the
$n_i$ must be $\geq 2$. For each $j \geq 1$, $i=1,2$,
$\Gamma_{ij}$ contains a subgroup $H_{ij}$ with the following
properties:

\vskip .05in $(a)$. $H_{ij}$ is non virtually abelian and the pair
$(\Gamma_{ij}, H_{ij})$ has the relative property (T) ([M]).

$(b)$. The normalizer $N_{ij}$ of $H_{ij}$ in $\Gamma_{ij}$
is ICC in $\Gamma_{ij}$ (i.e. $|\{hgh^{-1} \mid
h\in N_{ij} \}| = \infty$,
$\forall g\in \Gamma_{ij} \setminus \{e\}$) and $\sigma_{ij}$
is ergodic on $N_{ij}$.

$(c)$. For any proper intermediate subgroup $N_{ij} \subset
N_{ij}' \subset \Gamma_{ij}$ there exists $g \in \Gamma_{ij}
\setminus N_{ij}'$ such that $g(N_{ij}')g^{-1} \cap N_{ij}'$ is
non virtually abelian. \vskip .1in

Note that condition $(c)$ above on the inclusion $N_{ij} \subset
\Gamma_{ij}$ is similar to the {\it wq-normal} condition in
([P6,8]). It is equivalent to the existence of a well ordered,
strictly increasing family of intermediate subgroups $\{G_{l} \mid
0 \leq l \leq L \}$ such that $G_0 = N_{ij}$, $G_{L}=\Gamma_{ij}$
and $G_{k+1} = \{g \in \Gamma_{ij} \mid  gG_kg^{-1} \cap G_k$ non
virtually abelian$\}$, $\forall k$.

If $\Gamma_{ij}, j\geq 1,$ is $ICC$ and has a normal, non
virtually abelian subgroup with the relative property $(\text{\rm
T})$ then conditions $(a), (b), (c)$ are trivially satisfied.

Related to conditions $(a)$ and $(b)$, note that if a non virtually
abelian group $H$ is normal in an ICC group $G$ then $L(H)$ has no type I
summand, i.e. it is of type II$_1$. Indeed, because $G$ acts
ergodically on the center of $L(H)$, so if $L(H)$ has non-zero
type I part then it is homogeneous of type I$_n$ for some $2 \leq
n < \infty$, contradicting ([Ka]).

\proclaim{7.5. Theorem (vNE Bass-Serre rigidity)} With the
notations $7.5.0$ and assumptions $7.5.1$, if $\theta : M_1 \simeq
M_2^t$, for some $t > 0$, then $n_1 = n_2$ and there exists a
permutation $\pi$ of indices $j \geq 1$ and unitaries $u_j\in
M^t_2$ such that $\text{\rm
Ad}(u_j)(\theta(M_{1,j}))=M^t_{2,\pi(j)}$, $\text{\rm
Ad}(u_j)(\theta(A_1))=A_2^t$, $\forall j\geq 1$. In particular,
$\Cal R_{\sigma_1} \simeq \Cal R_{\sigma_2}^t$ and $\Cal
R_{\sigma_{1,j}} \simeq \Cal R^t_{\sigma_{2,\pi(j)}}$, $\forall
j\geq 1$. Also, if $\Gamma_{10}=\Gamma_{20}=1$, $2\leq n_1 <
\infty$ then the existence of such an isomorphism
forces $t=1$.
\endproclaim
\noindent {\it Proof}. We denote by $Q_{ij} \subset M_{ij}$ the
``rigid part'' of $M_{ij}$, i.e. $Q_{ij}=L(H_{ij})$ where
$H_{ij}\subset \Gamma_{ij}$ is a subgroup satisfying properties
$(a), (b), (c)$ in $7.5.1$. Also, we denote by $P_{ij}$ the von
Neumann algebra generated by the normalizer of $Q_{ij}$ in
$M_{ij}$. Thus, $P_{ij} \supset L(N_{ij})$, where $N_{ij}$ is the
normalizer of $H_{ij}$ in $\Gamma_{ij}$ in $(b)$. Notice that
conditions $7.5.1$  imply that $P_{ij}'\cap M_{ij} =\Bbb C$ so
that by 1.2.1 we also have $P_{ij}'\cap M_i = \Bbb C$, $\forall j
\geq 1$, $i=1,2$. Also, note that $Q_{ij}$ is relatively rigid in
$M_{ij}$, thus in $M_i$.

Assume first that $t \leq 1$. For simplicity, denote $A_2=A$,
$M_2=M$, $G_2=G$, $\Gamma_{2j}=\Gamma_j$ and let $\{u_g\}_{g\in G}
\subset M$ denote the canonical unitaries. Let $q\in A$ be so that
$\tau(q)=t$ and $\theta(1)=q$. Fix $i\geq 1$. Since $Q_{1i}\subset
M_1$ is a rigid inclusion, $Q=\theta(Q_{1i})\subset qMq$ is a
rigid inclusion ([P5]).

Let us show that no corner of $Q$ can be embedded into $M_{20}$
inside $M$. Assume by contradiction that there exist non-zero
$q_0\in \Cal P(Q)$, $p_0 \in \Cal P(M_{20})$, a unital isomorphism
$\psi$ of $q_0Qq_0$ into $p_0M_{20}p_0$ and a non-zero partial
isometry $v \in M_{20}$ such that $v^*v \in (q_0Qq_0)'\cap
q_0Mq_0$, $vv^* \in \psi(q_0Qq_0)'\cap p_0Mp_0$ and $vy =
\psi(y)v$, $\forall y\in q_0Qq_0$. Since $Q$ is type II,
$\psi(q_0Qq_0)\subset M_{20}$ is type II, so no corner of
$\psi(q_0Qq_0)$ can be embedded into $A$. By 1.2.1, this implies
$\psi(q_0Qq_0)'\cap M \subset M_{20}$. Thus $vQv^*=v(q_0Qq_0)v^*
\subset M_{20}$. But $\Gamma_{20}$ has Haagerup property, so by
([P5]) there exist unital trace preserving $A$-bimodular c.p. maps
$\phi_n$ on $M_{20}$ such that $\phi_n \rightarrow id_{M_{20}}$
and $\phi_n$ compact relative to $A$. Then $\phi_n * id
\rightarrow id_{M}$. But by ([P5]), $vQv^* \subset vv^* M vv^*$ is
a rigid inclusion. By ([P5]) again this implies $\lim_n \|(\phi_n
* id) (x) - x\|_2=0$ uniformly for $x\in (vQv^*)_1$. Since
$(\phi_n * id)(x)=\phi_n(x)$ for $x\in M_{20} \supset vQv^*$, this
implies the maps $\phi_n$, which are $A$-bimodular and compact
relative to $A$, tend uniformly to the identity on the unit ball
of the type II$_1$ algebra $vQv^*$. By ([P5]), this implies a
corner of $vQv^*$ can be embedded into $A$ inside $M_{20}$, a
contradiction.

Since no corner of $Q$ can be embedded into $M_{20}$ inside $M$ we
can apply Theorem 5.1. Thus, there exist $j=j(i)\geq 1$, a
non-zero projection $q'\in Q'\cap qMq$ and a unitary element $u\in
M$ such that $uQq'u^* \subset qM_{2,j}q=q(A \rtimes
\Gamma_{2,j})q$. Since $P=\theta(P_{1i})$ is generated by the
normalizer of $Q$, we have $q'\in P$ and by (3.5.1$^\circ$ in
[P6]), $Qq' \subset q'Pq'$ is quasi-regular. By Theorem 1.2.1 it
follows that $uq'Pq'u^*\subset qM_{2,j}q$. Since $P$ is a factor
and a corner of it is contained in $qM_{2j}q$, it follows that $u$
can be suitably modified so that to satisfy $uPu^*\subset
qM_{2j}q$. If $v\in qMq$ is such that $P_v = vPv^* \cap P$ is type
II$_1$, then $P_v\subset qM_{2,j}q$ and no corner of $P_v$ can be
embedded into the abelian algebra $A$. Since $(P_v)v = vP \subset
vM_{2j}$, by Theorem 1.2.1 it follows that $v\in M_{2,j}$. By
applying this recursively, by 7.5.1 $(c)$ it follows that
$\theta(L(\Gamma_{1i})) \subset qM_{2,j}q$.

Thus, if $\{u_{1,h}\mid h\in G_1 \}$ denotes the canonical
unitaries in $M_1 = A_1 \rtimes G_1$, then $v_h=\theta(u_{1,h})\in
qMq, h\in G_1,$ are in the normalizer of $\theta(A_1)$ in $qMq$.
Thus, the unitary elements $\{uv_hu^* \mid h \in \Gamma_{1i}\}
\subset qM_{2,j}q$ normalize $u\theta(A_1)u^*$ and they generate a
II$_1$ von Neumann algebra $N\subset qM_{2,j}q$. By Proposition
1.4.1, this implies $u\theta(A_1)u^* \subset M_{2,j}$ and
$w(u\theta(A_1)u^*)w^* = Aq$ for some $w\in \Cal U(qM_{ij}q)$.
Taking $v=wu$, it follows that $v(\theta(L(\Gamma_{1,i})))v^*
\subset M_{2,j}$ and $v\theta(A_1)v^* \subset M_{2,j}$, thus
$v\theta(M_{1,i})v^*\subset M_{2,j}$.

This shows that if $t \leq 1$ then $\forall i \geq 1$ there exists
$j=j(i) \geq 1$ (unique by Theorem 1.2.1) and a unitary element
$v\in M_2^t$ such that $\text{\rm Ad}(v)\theta(M_{1,i}) =
M_{2,j}^t$ and Ad$(v)(\theta(A_1))=A_2^t$.

Let us now consider the case $t > 1$. Let $n\geq t$ be an integer
and note that if we let $B=M_{n \times n}(A_2) \subset  M_{n
\times n} (M_2) = M$, $G=G_2$ and extend $\sigma_2$ to the action
$\sigma$ which acts trivially on $M_{n \times n}(\Bbb C)\subset
M_{n \times n} (M_2) = M$, then $M= B \rtimes_\sigma G_2$. Let $q
\in A=D_n \otimes A_2$ be a projection of trace $t/n$.

The hypothesis then states that $\theta: M_1 \simeq qMq$ is an
onto isomorphism. Fix $i \geq 1$ and let
$Q=\theta(Q_{1i})\subset qMq$, $P=\theta(P_{1i}) \subset qMq$. As
in the case $t \leq 1$, it follows that there exist $j=j(i)\geq 1$
and a unitary element $u$ in $qMq=M_2^t$ such that $uPu^* \subset
q(B \rtimes \Gamma_{2,j})q=M^t_{2,j}$. In particular, $\{
\theta(u_{1,h}) \mid h \in \Gamma_{1,i}\}\subset qM_{2,j}q$ and
they normalize $u\theta(A_1)u^*$. By applying Theorem 1.4.1 to
$A_0 = u\theta(A_1)u^* \subset qBq$, it follows that
$u(\theta(A_1))u^* \subset q(B \rtimes \Gamma_{2,j})q$ as well, so
Ad$(u)(\theta(M_{1,i}) \subset M_{2,j}^t$. Since
Ad$(u)(\theta(A_1))$ is regular in $M_2^t=qMq$, by Corollary 1.2.4
and (A.1 in [P5]) it follows that there exists a unitary in
$M_{2,j}$ that conjugates $\theta(A_1)$ onto $A_2^t$.

Since we have dealt with both cases $t \geq 1$ and $t \leq 1$, we
can apply the above equally well to $\theta$ and $\theta^{-1}$, to
obtain the following: $\forall n_1 \geq i \geq 1, n_2 \geq j \geq
1$, there exist unitaries $u_i \in M_2^t$, $v_j\in M_1$ and
indices $j(i) \in \{1,2,..., n_2\}, i(j) \in\{1, 2, ..., n_1\}$
such that
$$
\theta(A_1 \subset M_{1,i})\subset u_i(A_2\subset M_{2,j(i)})^t
u_i^*, 1\leq i \leq n_1,
$$
$$
\theta^{-1}((A_2 \subset M_{2,j})^t) \subset v_j(A_1\subset
M_{1,i(j)})v_j^*, 1\leq j \leq n_2,
$$
so altogether
$$
\theta( M_{1,k}) \subset \text{\rm Ad}(u_k \theta(v_{j(k)}))
(M_{1,i(j(k))}), 1\leq k \leq n_1.
$$
which by Theorem 1.1 entails $i(j(k))=k$, $\forall k$. Similarly
$j(i(k))=k$, $\forall k$. Thus, $n_1 = n_2 = n$ and $j$ defines a
permutation $\pi$ of the set of indices $1 \leq i \leq n$.

To prove the last part note that the equivalence $\Cal R_{\sigma_1}
\simeq \Cal R^t_{\sigma_2}$ and ([G1]) imply
$$
\Sigma_{j=1}^{n} \beta_1(\Gamma_{1j}) + (n-1) = \beta_1(\Cal
R_{\sigma_1}) = \beta_1(\Cal R_{\sigma_2})/t = (\Sigma_{j=1}^{n}
\beta_1(\Gamma_{2j}) +(n-1))/t.
$$
On the other hand, by the equivalence
$\Cal R_{\sigma_{1j}}  \simeq \Cal
R_{\sigma_{2,\pi(j)}}^t$ we get
$\beta_1(\Gamma_{1j})=\beta_1(\Gamma_{2,\pi(j)})/t$, $\forall j
\geq 1$ while by  ([BeVa]) and assumptions $7.5.1$ all
$\Gamma_{1j}$'s have $0$ first $\ell^2$-Betti number,
$\beta_1(\Gamma_{1j})$. (Indeed,
this follows easily from Corollary 4 in [BeVa] and
the argument in the proof of 2.4 in [P8].) Altogether we get
$n-1=(n-1)/t$, implying that $t=1$.

\hfill $\square$.

Before stating specific OE applications, recall from ([Fu1]) that
two groups $\Gamma, \Lambda$ are said to be {\it measure
equivalent} (ME) with {\it dilation constant} $t>0$ if there
exists free, m.p. actions $(\sigma, \Gamma), (\theta,\Lambda)$
such that $\Cal R_{\sigma}\simeq \Cal R_\theta^t$. We'll use the
notation $\Gamma \sim_{OE_t} \Lambda$ to denote this property. It
was recently proved in ([G2]) that if $\Gamma_i \sim_{OE_1}
\Lambda_i$, $\forall i\geq 0$, then $*_j\Gamma_j \sim_{OE_1} *_j
\Lambda_j$. Note that an alternative proof of this fact follows
from 7.3.1$^\circ$ (see also the comment 2.27 in [MS]).

The OE rigidity result below, of Bass-Serre type, can alternatively
be viewed as a converse to Gaboriau's above ME result, for free
products of w-rigid ICC groups. Note however that we need the
actions involved  to be ``separately ergodic'' (not assumed in
[G2]).

\proclaim{7.6. Corollary (OE Bass-Serre rigidity)} Let $\Gamma_0,
\Lambda_0$ be Haagerup groups and $\Gamma_i, \Lambda_j$, $1\leq i
\leq n \leq \infty$, $1\leq j \leq m \leq \infty$, be ICC groups
having normal, non virtually abelian subgroups with the relative
property $(\text{\rm T})$. Let $\sigma$ $($resp. $\theta)$ be a
free ergodic m.p. action of $\Gamma=\Gamma_0
* \Gamma_1 * ...$ $($resp.
$\Lambda=\Lambda_0 * \Lambda_1 *...)$ on the probability space
such that $\sigma_j=\sigma_{|\Gamma_j}$ $($resp.
$\theta_j=\theta_{|\Lambda_j})$ is egodic $\forall j \geq 1$. If
$\Cal R_{\sigma, \Gamma} \simeq \Cal R_{\theta, \Lambda}^t$ then
$n=m$ and there exists a permutation $\pi$ of the set of indices
$\geq 1$ such that $\Cal R_{\sigma_i,\Gamma_i} \simeq \Cal
R_{\theta_{\pi(i)}, \Lambda_{\pi(i)}}^t$, $\forall i\geq 1$.
\endproclaim

The condition on the groups $\Gamma_i, \Lambda_j$, $i,j\geq 1$, in
7.6 can be weakened by using the full generality of Theorem 7.5:

\proclaim{7.6'. Corollary} The same statement as in $7.6$ holds
true if we assume each $(\sigma_{|\Gamma_i}, \Gamma_i)$ $($resp.
$(\theta_{|\Lambda_i}, \Lambda_i))$, $i \geq 1$, satisfies
conditions $7.5.1 (a), (b), (c)$.
\endproclaim

\proclaim{7.7. Corollary} Let $\Gamma_i, 0\leq i \leq n$,
$\Gamma=\Gamma_0 * \Gamma_1 *...$, $\sigma$ be as in $7.6$ or
$7.6'$. Assume $\text{\rm Out}(\Cal R_{\sigma_1, \Gamma_1})=\{1\}$
and $(\sigma_1, \Gamma_1)$ is not orbit equivalent to $(\sigma_i,
\Gamma_i)$ for any $i \neq 1$. Then $\text{\rm Out}(\Cal
R_{\sigma, \Gamma})=\{1\}$ and $\text{\rm Out}(A \rtimes_\sigma
\Gamma)=\text{\rm H}^1(\sigma,\Gamma)$. Also, if $n$ is finite and
either all $\Gamma_i$ are finitely generated or there exists
$i\geq 1$ with $\beta_n(\Gamma_i)\neq 0, \infty$, then $\mycal F(A
\rtimes_\sigma \Gamma)=\{1\}$.
\endproclaim
\noindent {\it Proof}.  The first part is trivial by Theorem 7.5
while the last part follows from ([G1]). \hfill $\square$

Outer automorphism groups of equivalence relations are usually
hard to calculate and there are only a few special families of
group-actions $(\sigma_1, \Gamma_1)$ for which one knows that
Out$(\Cal R_{\sigma_1})=\{1\}$ (cf. [Ge2], [Fu2],[MS]). Similarly
for the 1-cohomology group H$^1$, where the only known
calculations are in ([Ge1], [PSa], [P8]). We remind below some
examples from ([Ge2], [Fu2], [MS]) where both calculations can be
made. We add a new construction of examples, in 7.8.3 below, which
uses the Monod-Shalom OE rigidity theorem to calculate Out and
([P8]) to calculate H$^1$.

\vskip .05in \noindent {\it 7.8.1. Example} (Gefter [Ge2], Furman
[Fu2]). Take $\Gamma_1$ a lattice in $SO(p,q)$, with $p > q \geq
2$ and notice that $\Gamma_1$ has rk$_{\Bbb R} \geq 2$ (thus has
property T) and admits a dense embedding into the compact Lie
group $SO(n)$, where $n=p+q$. Let $\sigma_1$ be the action by left
translation of $\Gamma_1$ on the homogeneous space
$SO(n)/SO(n-1)$. Then Out$(\Cal R_{\sigma_1})=\{1\}$.

\vskip .05in \noindent {\it 7.8.2. Example} (Monod-Shalom [MS]).
Let $\Cal G=SO(n), n\geq 5,$ and $\Lambda_0$ an ICC,  torsion
free, Kazhdan group which admits a dense embedding into $\Cal G$
and has no outer automorphisms (such groups exists for any $n \geq
5$, for instance lattices in $SO(p,q)$ with $p>q\geq 2$ and
$n=p+q$ as in 7.8.1). Let $K$ be any torsion free group embeddable
into $\Cal G$ and $K_0 \subset K$ a non-trivial subgroup such that
$K_0$ is not isomorphic to $K$ (for instance, $K_0=\Bbb F_r
\subset \Bbb F_s=K$, for some $s > r \geq 1$). Let $\Gamma_1 =
(\Lambda_0 * K) \times (\Lambda_0 * K_0)$ and note that any
automorphism of the group  $\Gamma_1$ is inner on $\Lambda_0
\times \Lambda_0$ (by Kurosh or Bass-Serre). Let $\sigma_1$ be the
action of $\Gamma_1=(\Lambda_0
* K) \times (\Lambda_0 * K_0)$ on $\Cal G$ by left-right
translation. Notice that this action is free ergodic on $\Lambda_0
\times \Lambda_0$ and that by ([MS]) one can choose the embedding
$K \subset \Cal G$ such that $\sigma_1$ is free on $\Gamma_1$ (by
considering all embeddings $gKg^{-1}$, $g \in \Cal G,$ and using a
Baire category argument). Then $\Gamma_1$ is in the class $\Cal
C_{geom}$ of Monod-Shalom and Out$(\Cal R_{\sigma_1})=\{1\}$.

\vskip .05in \noindent {\it 7.8.3. Example.} Let this time
$\Lambda_0$ be any torsion free ICC group with only inner
automorphisms and which cannot be decomposed as $\Bbb Z *
\Lambda_0'$ (note that if $\Lambda_0$ is w-rigid then it does have
this latter ``free in-decomposability'' property). Let $K$ be any
torsion free group with $K_0 \subset K$ a non-trivial subgroup
such that $K_0$ is not isomorphic to $K$. Let $\Lambda=\Lambda_0 *
K$ and $\Gamma_1 = (\Lambda_0 * K) \times (\Lambda_0 * K_0)\subset
\Lambda \times \Lambda$. Denote by $\sigma_1$ the action of
$\Gamma_1$ on the product probability space $(X, \mu)=\Pi_{g\in
\Lambda} (X_0, \mu_0)_g$ by left-right (double) Bernoulli shifts.
Then Out$(\Cal R_{\sigma_1})=\{1\}$. Indeed, noticing that
$\Lambda_0 * K, \Lambda_0 * K_0 \in \Cal C_{geom}$ and that
$\sigma_1$ is separately ergodic (on $\Lambda_0 * K$ and
$\Lambda_0 * K_0$) it follows from ([MS]) that any automorphism
$\theta \in$Aut$(\Cal R_{\sigma_1})$ is an inner perturbation of a
conjugacy $\sigma_1 \sim \sigma_1 \circ \gamma$ with respect to
some $\gamma \in $ Aut$(\Gamma_1)$. But $\Gamma_1$ has only inner
automorphisms (again by Kurosh or Bass-Serre). Thus, any such
$\theta$ is an inner perturbation of an automorphism of the
probability space that commutes with $\sigma_1(\Gamma_1)$. But
this commutant is trivial if for instance $(X_0,\mu_0)$ is atomic
with un-equal weights, as shown by the following:

\proclaim{7.9. Proposition} Let $\Lambda$ be a countable discrete
group with two subgroups $\Lambda_1, \Lambda_2 \subset \Lambda$
such that $\Lambda_0 = \Lambda_1 \cap \Lambda_2$ satisfies $|
\{hgh^{-1} \mid h \in \Lambda_0\} |=\infty$, $\forall g \in
\Lambda, g\neq e$. Let $(B_0,\tau_0)$ be a finite von Neumann
algebra and denote $(B,\tau) = \Pi_{g\in \Lambda} (B_0,\tau_0)_g$.
Let $\sigma$ be the action of $\Lambda_1 \times \Lambda_2$ on
$(B,\tau)$ given by $\sigma_{h_1,h_2}((x_g)_g)=(x'_g)_g$, where
$x'_g = x_{h_1^{-1}gh_2}$, $\forall g\in \Lambda, h_i\in
\Lambda_i$, $i=1,2$. Then $\sigma$ is a free, separately mixing
action of $\Lambda_1 \times \Lambda_2$ on $(B,\tau)$ and we have:

$1^\circ$. If $\theta \in \text{\rm Aut}(B,\tau)$ commutes with
$\sigma(\Lambda_1 \times \Lambda_2)$ then there exists a unique
$\theta_0 \in \text{\rm Aut}(B_0,\tau_0)$ such that $\theta$ is
the product action given by $\theta_0$, i.e. $\theta=\Pi_g
(\theta_0)_g$.

$2^\circ$. Any $\delta \in \text{\rm Aut}(\Lambda)$ satisfying
$\delta(\Lambda_i) = \Lambda_i$, $i=1,2$, implements an
automorphism $\Delta=\Delta(\delta) \in \text{\rm Aut}(B,\tau)$,
by $\Delta((b_g)_g)=(\delta(b_g))_g$, which satisfies $\Delta
\sigma \Delta^{-1} =\sigma \circ \delta$.

$3^\circ$. If $(B_0, \tau_0)=(L^\infty(X_0, \mu_0), \int \cdot
d\nu_0)$ for some atomic probability space  $(X_0, \mu_0)$, then
the commutant of $\sigma(\Lambda_1 \times \Lambda_2)$ in
$\text{\rm Aut}(X, \mu)$ is equal to $\text{\rm Aut}(X_0,
\mu_0)=\text{\rm Aut}(B_0, \tau_0)$. Moreover, if $\sigma$ is
conjugate to another double Bernoulli shift $\sigma'$ with atomic
base space $(X_0',\mu_0')$, then $(X_0, \mu_0) \simeq
(X_0',\mu_0')$.
\endproclaim
\vskip .05in \noindent {\it Proof}. Part $2^\circ$ is trivial. To
prove part $1^\circ$, it is sufficient to show that any $\theta$
commuting with $\sigma$ must take the subalgebra $B_0^e=... 1
\otimes (B_0)_e \otimes 1...$ of $B$ into itself. Indeed, because
if we denote $\theta_0=\theta_{|B_0^e}$ and regard it as an
automorphism of $B_0$ then $\theta \circ ( \Pi_g
(\theta_0)_g)^{-1}$ still commutes with $\sigma(\Lambda_1 \times
\Lambda_2)$ and it acts as the identity on $B^e_0$, thus on
$\sigma(g_1,g_2)(B^e_0)$, $\forall g_1, g_2 \in \Lambda_0$. Since
$\Lambda_1\Lambda_2 = \Lambda$, the latter generate all $B$, thus
$\theta=\Pi_g (\theta_0)_g$.

To show that $\theta$ leaves $B_0^e$ globally invariant, it is
sufficient to show that that the fixed point algebra $\{b \in B
\mid \sigma(g,g)(b)=b, \forall g\in \Lambda_0\}$ coincides with
$B_0^e$. This in turn follows trivially from the fact that for any
finite subset $F \subset \Lambda\setminus \{e\}$ there exists $g
\in \Lambda_0$ such that $gFg^{-1} \cap F =\emptyset$. To see that
this latter property holds true, note that if some finite set
$\emptyset \neq F\subset \Lambda \setminus \{e\}$ would satisfy
$|gFg^{-1} \cap F| \geq 1$, $\forall g\in \Lambda_0$, then the
``left-right'' representation $\pi(g)(f) = \lambda(g)\rho(g)(f)$
on $\ell^2(\Lambda)$ would satisfy $\langle \pi(g)(\chi_F), \chi_F
\rangle \geq |F|^{-1}$ for all $g \in \Lambda_0$. Taking the
element $f$ of minimal Hilbert norm in $\overline{\text{\rm co}}^w
\{ \pi(g)(\chi_F) \mid g\in \Lambda_0\}\subset \ell^2(\Lambda)$,
it follows that $f\geq 0$, $f \neq 0$ (because $\langle f,
\chi_F\rangle \geq |F|^{-1}$) and $\pi(g)(f)=f, \forall g\in
\Lambda$. But then any appropriate ``level set'' $K$ for $f$ will
be finite, non-empty and will satisfy $\pi(g)(\chi_K)=\chi_K$,
$\forall g\in \Lambda_0$, i.e. $gKg^{-1}=K$, $\forall g\in
\Lambda_0$, implying that $\Lambda$ has elements with finite
conjugacy class, a contradiction.

The first part of $3^\circ$ is trivial by $2^\circ$. Then to see
that conjugacy of double Bernoulli shifts entails isomorphism of
the base spaces, note that $\sigma$ conjugate to $\sigma'$ implies
that all ``diagonal'' actions  $\sigma \otimes \sigma''$,
$\sigma'\otimes \sigma''$ must also be conjugate, $\forall
\sigma''$, thus having isomorphic commutants in Aut. Taking
$\sigma''$ to be itself a double Bernoulli $\Gamma_1$-action of
base $(X_0'', \mu_0'')$, it follows that Aut$((X_0, \mu_0) \times
(X_0'', \mu_0'')) \simeq \text{\rm Aut}((X_0', \mu_0') \times
(X_0'', \mu_0''))$, $\forall (X_0'', \mu_0'')$, which easily
implies the result. \hfill $\square$

\vskip .05in \noindent {\it 7.10.0. Notation}. Denote by $w\Cal
T_2$ the class of groups $G$ that have a non virtually abelian
subgroup $H_0 \subset G$ such that: $(G, H_0)$ is a property (T)
pair; the normalizer $H$ of $H_0$ in $G$ satisfies $|\{hgh^{-1}
\mid h\in H\}|=\infty$, $\forall g \in G\setminus \{e\}$; the
wq-normalizer of $H$ in $G$ generates $G$. Note that any group in
$w\Cal T_2$ is ICC and that if $G$ is ICC and has a normal, non
abelian, relatively rigid subgroup then $G \in w\Cal T_2$. Thus,
any group of the form $G=H_0 \times K$ with $H_0$ ICC Kazhdan and
$K$ either ICC or equal to 1 is w-rigid and thus in $w\Cal T_2$.
Also, if $G\in w\Cal T_2$ then $(G
* K_0)\times K \in w\Cal T_2$ for any ICC group $K$ and any
arbitrary group $K_0$.

\proclaim{7.10. Corollary} Let $\Gamma_0$ be a Haagerup group and
$\Gamma_i \in w\Cal T_2$, for $1 \leq i \leq n$, where $ 1\leq n
\leq \infty$. Assume $\Gamma_1$ is as in $7.8.1, 7.8.2$, or
$7.8.3$. Then $\Gamma=\Gamma_0
* \Gamma_1 *... $ has a free ergodic m.p. action $\sigma$ with
$\text{\rm Out}(\Cal R_{\sigma})=\{1\}$. Moreover, we have:

$1^\circ$. If $\Gamma_1$ is as in $7.8.2$ or $7.8.3$, then there
exist uncountably many non stably orbit equivalent actions
$\sigma$ of $\Gamma$ with $\text{\rm Out}(\Cal R_{\sigma})=\{1\}$
and $\mycal F(\Cal R_\sigma)=\{1\}$.

$2^\circ$. If in addition $\Gamma_0$ is a product of amenable
groups then given any discrete countable abelian group $K$, the
uncountable family of actions $\sigma$ in $1^\circ$ can be taken
to satisfy $\text{\rm H}^1(\sigma, \Gamma)= \Bbb G_0^{n-1} \times
\Bbb G \times \Pi_{j \geq 1} \text{\rm Char}(\Gamma_j) \times
K^{n-1}$, where $\Bbb G$ is the polish group $\Cal U(L^\infty(\Bbb
T, \lambda))$ and $\Bbb G_0=\Bbb G/\Bbb T$.
\endproclaim
\noindent {\it Proof}. By 7.3, we can take the free m.p. action
$\sigma$ of $\Gamma$ on $A = L^{\infty}(X, \mu)$ so that for each
$i\neq 1$,  $\sigma_i=\sigma_{|\Gamma_i}$ is a (left) Bernoulli
$\Gamma_i$-action, or a quotient of it as in ([P8]), and of the
form $7.8.1-7.8.3$ for $i=1$. Notice that in example $7.7.1$ the
group $\Gamma_1$ has property (T) and is ICC, thus $\Gamma_1\in
w\Cal T_2$. Then in both $7.8.2$ and $7.8.3$ , $\Gamma_1$ has
$\Lambda_0 \times \Lambda_0$ as a relatively rigid subgroup, which
is wq-normal in $\Gamma_1$, with $\Lambda\times \Lambda$ ICC in
$\Gamma_1$. Furthermore, since Out$(\Cal R_{\sigma_i})$ is huge
for $i \geq 2$ and trivial for $i=1$, $\sigma_1$ cannot be stably
orbit equivalent to $\sigma_i, i\geq 2$. The existence of ``many
actions'' $\sigma$ in the case 7.8.2 and 7.8.3 follows from the
existence of uncountably many non-stably OE relations $(\sigma_1,
\Gamma_1)$ of the form $7.8.2$ (cf [MS]) and respectively of the
form $7.8.3$ (by Proposition 7.9.3$^\circ$ and [MS]).

Finally, the calculation of the 1-cohomology groups follows from
(2.12, 3.1, 3.2 in [P8]) and the fact that in all cases 7.8.1,
7.8.2, 7.8.3 one has H$^1(\sigma_1,\Gamma_1)=\text{\rm
Char}(\Gamma_1)$. Indeed, in case $\sigma_1$ is as in 7.8.1 or
7.8.2 then this calculation follows from ([Ge1] and 2.12 in [P8]),
while in case 7.8.3 the calculation  is in ([P8]). \hfill
$\square$

\proclaim{7.11. Corollary} With $\Gamma= *_{i\geq 0} \Gamma_i$,
$K$, $\sigma$ as in $7.10.2^\circ$,  denote $A=L^\infty(X,\mu)$
and $M=A \rtimes_\sigma \Gamma$. Then $\mycal F(M)=\{1\}$ and
$$
\text{\rm Out}(M)=\text{\rm H}^1(\sigma,\Gamma)=\Bbb G_0^{n-1}
\times \Bbb G \times \Pi_{j \geq 1} \text{\rm Char}(\Gamma_j)
\times K^{n-1}.
$$
\endproclaim
\noindent {\it Proof.} Trivial by 7.7 and 7.10. \hfill $\square$
\vskip .05in

Note that Out$(M)$ is abelian and non-locally compact in all
examples in 7.11 above, but if we denote $\tilde{\text{\rm
Out}}(M)$ the quotient of $\text{\rm Out}(M)$ by the connected
component of $id_M$ (which is closed in Out$(M)$, with the latter
a polish group in all examples considered), then $\tilde{\text{\rm
Out}}(M)$ is the quotient of $\Pi_{j \geq 1} \text{\rm
Char}(\Gamma_j) \times K^{n-1}$ by the connected component of 1,
which for $n < \infty$ is totally disconnected separable, locally
compact group.

We end this section by mentioning another rigidity result, which
from an isomorphism of group measure space factors corresponding
to relatively rigid actions ([P5]) of free products of groups
derives the orbit equivalence of the actions. This type of results
were first obtained in ([P5]) for HT group-actions, and in ([P6]),
for Bernoulli shift actions of groups containing infinite
subgroups (not necessarily normal) with the relative property (T).

\proclaim{7.12. Theorem (vNE/OE rigidity)} Let $(M_i, \tau_i)$ be
type $\text{\rm II}_1$ von Neumann algebras with a common Cartan
subalgebra $A\subset M_i$, $i=1,2$, such that
$\tau_{1|A}=\tau_{2|A}$. Assume $M=M_1 *_A M_2$ is a factor and
$A$ is Cartan in $M$ $($N.B. By $7.2$ this is the same as
requiring that $\Cal R_{A\subset M_i}$, $i=1,2$, are freely
independent$)$. If $A_0\subset M^t$ is a rigid Cartan subalgebra,
for some $t > 0$, then there exists a unitary element $u \in M^t$
such that $uA_0u^* = A^t$.
\endproclaim
\noindent {\it Proof.}  It is clearly sufficient to prove this in
the case $t=1$. If some corner of $A_0$ can be embedded into $A$
inside $M$, then the statement follows by (A.1 in [P5]). If we
assume this is not the case, then we can apply Theorem 7.5 to get
a non-zero $p\in \Cal P(A_0)$ such that $vA_0pv^*\subset M_i$ for
some $i\in\{1,2\}$ and $v\in \Cal U(M)$. Moreover, since $M$ is a
factor and $A_0$ is Cartan in $M$, we may assume $vpv^*$ is
central in $M_i$, so in particular $p_1=vpv^*$ lies in $A$ (the
latter being maximal abelian in $M_i$). Then $vA_0pv^*$ is Cartan
in $p_1Mp_1$, so by Corollary 1.2.4 we have $p_1M_2p_1=Ap_1$,
contradicting the fact that $M_2$ is type II. \hfill $\square$

\vskip .05in

\proclaim{7.13. Corollary} Let $\Cal F$ be the class of groups
that can be written as a free product of two $($or more$)$
infinite groups. Let $\sigma : G \rightarrow (X, \mu)$ be a free
ergodic m.p. action of a group $G \in \Cal F$ and denote
$M=L^\infty(X,\mu) \rtimes_\sigma G$ the corresponding group
measure space $\text{\rm II}_1$ factor, with $A=L^\infty(X,\mu)
\subset M$ the corresponding Cartan subalgebra. Let $M_0$ be a
$\text{\rm II}_1$ factor with a relatively rigid Cartan subalgebra
$A_0 \subset M_0$ and denote $\Cal R_0=\Cal R_{A_0\subset M_0}$.

$1^\circ$. If $\theta: M_0 \simeq M^t$ for some $t > 0$ then
$\theta$ can be perturbed by an inner automorphism so that to take
$A_0$ onto $A^t$. In particular, $\Cal R_{0} \simeq \Cal
R_\sigma^t$ and thus $\beta_n(\Cal R_0)=\beta_n(G)/t$, $\forall
n$.

$2^\circ$. If $G \in \Cal F$ satisfies $\beta_n(G) \neq 0,\infty$
for some $n$ $($for instance, if $G=\Gamma_1 * \Gamma_2$, with
$\Gamma_{1,2}$ finitely generated infinite groups, in which case
$\beta_1(G)\neq 0,\infty)$ and the action $\sigma$ is relatively
rigid, then $\mycal F(M)=\{1\}$.
\endproclaim
\noindent {\it Proof.} All statements are trivial by 7.12 and
([G1]). \hfill $\square$

\vskip .05in Note that 7.13.$2^\circ$ above shows in particular
that $\mycal F(L(\Bbb Z^2 \rtimes \Bbb F_n))=\{1\}$ for any $\Bbb
F_n \subset SL(2,\Bbb Z)$, with $2 \leq n < \infty$, thus giving a
new proof (but still using Gaboriau's work [G1]) of one of the
main result in ([P5]).

\vskip .05in \noindent {\it 7.14. Definition}. A countable,
measurable, standard m.p. equivalence relations $\Cal R$ is a $FT$
{\it equivalence relation} if it is of the form $\Cal R=\Cal
R^t_\sigma$, where $t > 0$ and $(\sigma,\Gamma)$ are free ergodic
m.p. actions on the probability space $(X,\mu)$ with the following
properties: $(a)$. The group $\Gamma$ is a free product of two (or
more) infinite groups; $(b)$. $\sigma$ is relatively rigid, in the
sense of (Definition 5.10.1 in [P5]), i.e. $L^{\infty}(X,\mu)
\subset L^\infty(X,\mu) \rtimes_\sigma G$ is a rigid inclusion
(4.2.1 in [P5]). The above Corollary 7.13 thus shows that all OE
invariants for $FT$ equivalence relations $\Cal R$ are in fact vNE
invariants for $\Cal R$, i.e. are isomorphism invariants of the
associated group measure space II$_1$ factors $M=L(\Cal R,w)$,
where $w\in \text{\rm H}^2(\Cal R)$. We denote by $\Cal F\Cal T$
the class of all such II$_1$ factors $M$.

Note that if $\Cal R=\Cal R^t_{\sigma,\Gamma}$ for some free
ergodic m.p. action $(\sigma,\Gamma)$ with $\Gamma$ a free product
of two infinite groups, then $\Cal R$ is $HT_{_{s}}$ in the sense
of ([P5]) if and only if it is $FT$ and $\Gamma$ has Haagerup
property. Thus, all equivalence relations coming from
amplifications of actions $\sigma$ of non-amenable subgroups
$\Gamma \subset SL(2,\Bbb Z)$ on $L^\infty(\Bbb T^2,\lambda)$ are
$FT$ actions. However, actions $\sigma$ of groups such as
$SL(n,\Bbb Z)* H$, with $H$ infinite group and $\sigma_{|SL(n,\Bbb
Z)}$ isomorphic to the canonical action of $SL(n,\Bbb Z)$ on $(\Bbb
T^n,\lambda)$, give $FT$ equivalence relations which are not
$HT_{_{s}}$. Thus, the class $\Cal F\Cal T$ provides additional
group measure space II$_1$ factors for which orbit equivalence
invariants of the actions, such as Gaboriau's $\ell^2$-Betti
numbers, become isomorphism invariants of the factors.

We end by mentioning an application of 7.3.1$^\circ$ which brings
some light to (Problems 5.10.2, 6.12.1 in [P5]) and to the problem
of existence of ``many'' non-OE actions for non-amenable groups,
as a consequence of (Corollary 7 in [GP]):

\proclaim{7.15. Corollary} $1^\circ$. The class of groups $\Gamma$
which admit free ergodic  relatively rigid m.p. actions on the
probability space is closed to free products with arbitrary groups
$\Gamma'$. Also, if $\Gamma$ is a $H_{T_{_{s}}}$ group $($i.e.
$\Gamma$ has Haagerup's property and admits relatively rigid
actions, see 6.11 in ${\text{\rm [P5]}})$, then $\Gamma * \Gamma'$
is $H_{T_{_{s}}}$, for any $\Gamma'$ with Haagerup property.

$2^\circ.$ If $\Gamma$ admits a relatively rigid action $($e.g. if
$\Gamma\subset SL(2,\Bbb Z)$ non-amenable, or $\Gamma$ an
arithmetic  lattice in an absolutely simple, non-compact Lie group
with trivial center, cf $\text{\rm [P5], [Va]})$ and $\Gamma'$ is
an arbitrary infinite amenable group then $\Gamma
* \Gamma'$ has uncountably many non-stably OE free ergodic m.p. actions
on the probability space. Also, if $\Gamma_0$ is an arbitrary
group and $\Gamma_1,\Gamma_2$ are non-trivial amenable groups, at
least one having more than $2$ elements, then $\Gamma_0 * \Gamma_1
*\Gamma_2$ has uncountably many non-stably OE free ergodic m.p.
actions on the probability space.
\endproclaim
\noindent {\it Proof.} Part $1^\circ$ is a trivial consequence of
7.3.1$^\circ$ and of the (trivial) property of relatively rigid
equivalence relations $\Cal R$ that any $\Cal R_0$ that contains
$\Cal R$ is also relatively rigid (see e.g. 4.6.2$^\circ$ in
[P5]). Part 2$^\circ$ is just the combination of (Corollary 7 in
[GP]) and 7.3.
\hfill $\square$

\heading 8. Amalgamation over $R$: Factors with no outer
automorphisms
\endheading

We prove here another rigidity result for AFP factors, this time
in the case $M=M_0*_R M_1 *_R * ...$, where $R$ is the hyperfinite
II$_1$ factor. As an application, we obtain factors $M$ with
Out$(M)=\{1\}$, thus answering a well known problem posed by A.
Connes in 1973.

Like in the group measure space case in Section 7, we only
consider crossed product inclusions $(R \subset M_i)=(R\subset R
\rtimes_{\sigma_i} \Gamma_i)$, with the $\sigma_i$ freely
independent, i.e. implementing a free action $\sigma$ of
$\Gamma=\Gamma_0 * \Gamma_1
* ...$ on $R$. Thus, $M$ will be viewed alternatively as a crossed
product factor $M=R \rtimes_\sigma \Gamma$, with the algebra of
coefficients $R$ having trivial relative commutant in $M$.

The key assumption is that the action $(\sigma, *_i \Gamma_i)$ has
the relative property (T), i.e. $R \subset M$ be a rigid inclusion
in the sense of ([P5]). The rigidity result shows the uniqueness,
modulo unitary conjugacy, of the ``core'' $R$ of such factors.
Since the normalizer of $R$ in $M$ completely encodes the group
$\Gamma$, we can recover completely the isomorphism class of the
groups $\Gamma_i$, by classical Bass-Serre theory. The result is
similar to the vNE/OE Rigidity 7.11 (where however only the orbit
equivalence class of $\Gamma$ could be recovered) and to the
unique crossed-product decomposition result in ([P9]). But since
we also get the component by component unitary conjugacy of the
factors $M_i$, it is again a Bass-Serre type rigidity result.

Through this theorem, the calculation of Out$(M)$ and $\mycal
F(M)$ reduces to the calculation of the commutant of
$\sigma(\Gamma)$ in Out$(R)$, like in  [P9], where however no such
commutant could be calculated! This time, due to Bass-Serre
arguments and the possibility of choosing the actions $(\sigma_i,
\Gamma_i)$ with prescribed properties (cf. 8.2 below), we can
control such commutants and calculate Out$(M)$ completely for
large classes of factors.

\proclaim{8.1. Lemma} Let $M_n, n \geq 0,$ be $\text{\rm II}_1$
factors with a common subfactor $N \subset M_n$. Then $N \subset
M=M_0 *_N M_1 *_N M_2 *_N...$ is irreducible and regular if and
only if $N\subset M_n$ is regular, irreducible $\forall n \geq 1$
and the groups of outer automorphisms $\Gamma_n = \{\text{\rm
Ad}(u) \mid u \in \Cal N_{M_n}(N) \}/\Cal U(N)$, $n \geq 0$, on
$N$ are freely independent.
\endproclaim
\noindent {\it Proof}. Trivial by the definitions of freeness and
respectively of the amalgamated free product over $N$. \hfill
$\square$

\vskip .05in

The result below is the analogue for actions on the hyperfinite
II$_1$ factor $R$ of the result on the existence of freely
independent actions on the probability space in Proposition 7.3. It
shows the existence of free actions $\sigma$ of groups $\Gamma=
\Gamma_0
* \Gamma_1 * \Gamma_2 *...$ on $R$ such that the restriction of
$\sigma$ to each individual group $\Gamma_j$ is conjugate to a
prescribed free action of $\Gamma_j$ on $R$. It will be frequently
used in this Section. The proof relies on Lemma A.2 in the Appendix.

\proclaim{8.2. Proposition} Let $\sigma_n : \Gamma_n \rightarrow
\text{\rm Aut}(R)$ be free actions of countable discrete groups
$G_n$, $n \geq 0$. Then there exists a free action $\sigma$ of the
group $G=*_n \Gamma_n$ on $R$ such that $\sigma_{|\Gamma_n}$ is
conjugate to $\sigma_n$, $\forall n\geq 0$.
\endproclaim
\noindent
{\it Proof}. For each $n \geq 0$ denote $\tilde{G}_n =G_0 *
G_1 * ...* G_n$. Assume we have constructed a map
$\tilde{\sigma}_n$ of $\tilde{G}_n$ into $\text{\rm Aut}(R)$ such
that the quotient map $\tilde{\sigma}'_n$ of $G$ into $\text{\rm
Aut}(R)/\text{\rm Int}(R)$ is a faithful group morphism with
$\tilde{\sigma}_{n|G_i}$ conjugate to $\sigma_i$, $\forall 0 \leq
i \leq n$. We then apply A.2 to $\{\tilde{\sigma}_n(g) \mid g \in
\tilde{G}_n\} \cup \{\sigma_{n+1}(h) \mid h \in G_{n+1}\}$ to get
an automorphism $\theta_{n+1}$ of $R$ such that $\tilde{\sigma}_n
(\tilde{G}_n)$ and $\theta_{n+1} \sigma_{n+1}(G_{n+1})
\theta_{n+1}^{-1}$ are freely independent. Denoting by
$\tilde{\sigma}_{n+1}$ the map of $\tilde{G}_{n+1} = \tilde{G}_n
* G_{n+1}$ into $\text{\rm Aut}(R)$ which restricted to $\tilde{G}_n$ equals
$\tilde{\sigma}_n$ and restricted to $G_{n+1}$ equals
$\theta_{n+1} \sigma_{n+1}(G_{n+1}) \theta_{n+1}^{-1}$, the
statement follows by induction.
\hfill $\square$

\proclaim{8.3. Theorem} Let $G_i = \Gamma_{i,0} * \Gamma_{i,1}
* ... * \Gamma_{i,n_i}$ with $\Gamma_{ij}, 0\leq j \leq n_i$,
non-trivial groups, for some $1\leq n_i \leq \infty$, $i=1,2$. For
each $i=1,2$ let $\sigma_i:G_i \rightarrow \text{\rm Aut}(N_i)$ be
a free, ergodic m.p. action on a $\text{\rm II}_1$ factor $N_i$.
Denote $M_i=N_i \rtimes_{\sigma_i} G_i$, $i=1,2$ and assume $N_i
\subset M_i$ are rigid inclusions $i=1,2$. Let $\theta : M_1
\simeq M^t_2$, for some $t > 0$. Then we have:

$1^\circ$. There exists $u \in \Cal U(M_2^t)$ such that $\text{\rm
Ad}(u)(\theta(N_1))=N_2^t$. Thus, $G_1 \simeq G_2$ and $\sigma_1,
\sigma^t_2$ are cocycle conjugate actions with respect to the
identification $G_1\simeq G_2$.

$2^\circ$. If in addition $\Gamma_{i,0}$ are free groups and
$\Gamma_{ij}$ are free indecomposable, not equal to the infinite
cyclic group, $\forall 1\leq j \leq n_i$, $i=1,2$, then
$\Gamma_{10} \simeq \Gamma_{20}$, $n_1 = n_2$ and there exists a
permutation $\pi$ of the indices $j \geq 1$ and unitaries $u_j\in
M^t_2$ such that $\text{\rm
Ad}(u_j)(\theta(M_{1,j}))=M^t_{2,\pi(j)}$, $\text{\rm
Ad}(u_j)(\theta(N_1))=N^t_2$, $\forall j\geq 1$. In particular,
$\Gamma_{1,j}\simeq \Gamma_{2,\pi(j)}$ and $\sigma_{1,j},
\sigma^t_{2,\pi(j)}$ are cocycle conjugate with respect to this
identification of groups, $\forall j\geq 1$.
\endproclaim
\vskip .05in \noindent {\it Proof}. We first prove that a corner
of $\theta^{-1}(N_2^t)$ can be embedded into $N_1$ inside $M_1$.
Assume this is not the case. By applying recursively
5.1.3$^\circ$, it follows that there exist a unitary element $u\in
\Cal U(M_1)$ and some $n_1 \geq j \geq 1$ such that
$u\theta^{-1}(N_2^t)u^*\subset M_{1,j}$. Since
$u\theta^{-1}(N_2^t)u^*$ is regular in $M_1$, using again the
assumption by contradiction, Corollary 1.2.4 implies that a corner
of $u\theta^{-1}(N_2^t)u^*$ can be embedded into $N_1$ inside
$M_{1,j}$ (thus inside $M_1$ as well). Altogether, this shows that
a corner of $\theta^{-1}(N_2^t)$ can be embedded into $N_1$ inside
$M_1$.

Similarly, a corner of $\theta(N_1)^{1/t}=\theta^{1/t}(N_1^{1/t})$
can be embedded into $N_2$ inside $M_2$. Thus, a corner of $N_1$
can be embedded into $\theta^{-1}(N_2^t)$ inside $M_1$. Since both
$N_1$ and $\theta^{-1}(N_2^t)$ are regular in $M_1$, with $\Cal
N_{M_1}(N_1)/\Cal U(N_1) \simeq G_1$ and with the other similar
quotient isomorphic to $G_2$, and since both $G_1, G_2$ are ICC
(being free products of non-trivial groups), the unitary conjugacy
of $N_1, \theta^{-1}(N_2^t)$ in $M_1$ (equivalently of
$\theta(N_1), N_2^t$ in $M_2^t$) follows from the following
general:

\proclaim{8.4. Lemma} Let $M$ be a $\text{\rm II}_1$ factor and
$P, Q \subset M$ be irreducible, regular subfactors. Assume
$\Gamma=\Cal N_M(P)/\Cal U(P)$, $\Lambda = \Cal N_M(Q)/\Cal U(Q)$
are ICC groups. Also, assume each one of the inclusions $P \subset
M$, $Q \subset M$ is an amplification of a genuine crossed product
inclusion. If $L^2(M)$ contains non-zero $P-Q$ Hilbert bimodules
$\Cal H, \Cal K \subset L^2(M)$ such that $\text{\rm dim}
_{P}(\Cal H)< \infty$, $\text{\rm dim} (\Cal K)_{Q}< \infty$ then
$P, Q$ are unitary conjugate in $M$.
\endproclaim
\noindent {\it Proof}. We first prove that $L^2(M)$ is generated
by irreducible $P-Q$ Hilbert bimodules that are finite dimensional
both as left $P$ modules and as right $Q$ modules. We'll actually
prove this by only using the fact that $P,Q$ are quasi-regular in
$M$. Note that by (1.4 in [P5]) $\Cal H^0=\Cal H \cap M$ is dense
in $\Cal H$ and contains an orthonormal basis$/Q$. Similarly,
since $Q$ is quasi-regular in $M$ and it is a factor, $L^2(M)$ is
generated by Hilbert $Q-Q$ bimodules $\Cal H_\beta$ such that
$\Cal H^0_\beta = \Cal H_\beta \cap M$ is dense in $\Cal H_\beta$
and contains both left and right orthonormal basis$/Q$. But then
$\Cal H^0 \cdot \Cal H_\beta^0$ span all $L^2(M)$ and are finite
dimensional$/Q$. Equivalently, $P'\cap J_MQ'J_M$ is generated by
projections that have finite trace in $J_MQ'J_M$. Similarly, $P'
\cap JQ'J$ is generated by projections that have finite trace in
$P'$. Thus, $\Cal A=P'\cap JQ'J$ is generated by projections that
are finite with respect to both traces, thus corresponding to
Hilbert $P-Q$ bimodules which are finite dimensional both from
right and left. Since $P,Q$ are factors, by ([J]) each such
bimodule is a direct sum of irreducible bimodules.

Let now $\Cal H\subset L^2(M)$ be an irreducible $P-Q$ bimodule.
By ([J]) we have dim$(_P\Cal H) \cdot \text{\rm dim}(\Cal H_Q)
\geq 1$ and the equality means the orthogonal projection $p_\Cal
H$ of $L^2(M)$ onto $\Cal H$ satisfies: $p_\Cal H\in P'\cap
\langle M, e_Q \rangle$, $Tr_{\langle M, \rangle Q}(p_\Cal H)=1$
and $p_\Cal H\langle M, e_Q\rangle p_\Cal H=Pp_\Cal H$. Thus, by
the proof of (Lemma 1 in [P7]) $up_\Cal H u^*=e_Q$ for some $u\in
\Cal U(M)$, which also satisfies $uPu^*=Q$.

Assume now that $Tr(p_\Cal H)>1$, by (2.1 in [P6]) there exists a
projection $p\in P$, a unital isomorphism $\psi: pPp \rightarrow Q$
and a partial isometry $v\in M$ such that $vv^*=p$, $q'=v^*v \in
\psi(pPp)'\cap Q$ and $xv = v\psi(x), \forall x\in pPp$. Moreover,
the finite dimensionality plus irreducibility of $\Cal H$ as a $P-Q$
bimodule implies that $Q_1=\psi(pPp)$ has finite index in $Q$, has
trivial relative commutant in $Q$ and $q'$ is minimal in $Q_1'\cap
M$.

By appropriately amplifying $Q\subset M$ we may assume this
inclusion is a genuine crossed product inclusion $Q \subset
Q\rtimes_\sigma \Gamma$. Denote by $\{u_g\}_g\subset M=Q
\rtimes_\sigma \Gamma$ the canonical unitaries implementing
$\sigma$ on $Q$. Let $q'=\Sigma_g x_g u_g$, with $x_g \in Q$. By
identification of Fourier series, it follows that $x_gu_g \in
Q_1'\cap M$, $\forall g$. Thus $x_gx_g^*=x_gu_gu_g^*x_g^* \in
Q_1'\cap Q=\Bbb C$, so that all $x_g$ are scalar multiples of
unitaries in $Q$. Let $K_0\subset \Gamma$ be the support of this
Fourier expansion of $q'$. Let also $K \subset \Gamma$ be the set
of all $k\in \Gamma$ such that $u_k$ can be perturbed a unitary in
$Q$ so that to fix $Q_1$ pointwise. Since $[Q:Q_1] < \infty$ and
$Q_1'\cap Q=\Bbb C$, $K$ is a finite subgroup of $\Gamma$ and $K_0
\subset K$.

Since $Q_1'\cap Q = \Bbb C$, by Connes' vanishing $1$-cocycle for
finite groups, the unitaries $w_k \in Q$ satisfying
Ad$(w_ku_k)_{|Q_1}=id_{Q_1}$ can be chosen of the form
$\sigma_k(w)w^*$, $k\in K$, for some unitary element $w\in Q_1$.
Thus, by perturbing all $\{u_g\}_{g\in \Gamma}$ by a 1-cocycle we
may assume Ad$(u_k)$ act trivially on $Q_1$. Let $\Gamma_0\subset
\Gamma$ be the subgroup of all $g \in \Gamma$ such that $u_g$ can
be perturbed by a unitary in $Q$ so that to normalize $Q_1$.
Clearly $K \subset \Gamma_0$ and $K$ is normal in $\Gamma_0$. We
will prove that $\Gamma_0$ has finite index in $\Gamma$, thus
contradicting the hypothesis.

By the minimality of $q'$ in $Q_1'\cap M$, it follows that $q'$ is
minimal in the group algebra $L(K)=\text{\rm sp}\{u_k \mid k\in
K\}$. Identify $pPp \subset pMp$ with $Q_1q' \subset q'Mq'$ via
Ad$(v)$. Let $\{v_h\mid h \in \Lambda\}$ be a choice of canonical
unitaries in $M=P \rtimes \Lambda$, which we assume commute with
$p\in P$ (we can do that for each $h$ by perturbing if necessary
with unitaries in the factor $P$). For each $h \in \Lambda$,
$h\neq e$, let $v_h=\Sigma_g x^h_g u_g\in Q \rtimes_\sigma \Gamma$
be the Fourier expansion of $v_h$ and denote $\theta_h$ the action
implemented by Ad$(v_h)$ on $Q_1\simeq Q_1q'=pPp$. Thus, $v_h y =
\theta_h(y) v_h, \forall y\in Q_1$. Identifying the Fourier series
in $\{u_g\}_g$, this implies $\theta_h(y)(x^h_gu_g) = (x^h_gu_g)
y, \forall y\in Q_1$. As before, this implies each $x^h_g$ is a
scalar multiple of a unitary in $Q$ and that Ad$(x^h_gu_g)$
normalizes $Q_1$. Thus the support $K_h$ of the Fourier series for
$v_h$ is contained in $\Gamma_0$. Since $Q_1q'\subset q'Mq'$ is
the closure of the span of elements in $Q_1v_hq', h\in \Lambda$
and each $v_h$ is supported on $\Gamma_0$, as a Fourier expansion
in $\{u_g\}_g$ with coefficients in $Q$, it follows that
$q'Mq'\subset \oplus_{g\in \Gamma_0} L^2(Q)u_g$. In particular,
since $q'\in L(K)\subset L(\Gamma_0)$ we get
$q'L(\Gamma)q'=q'L(\Gamma_0)q'$, which clearly implies $\Gamma_0$
has finite index in $\Gamma$. But this implies $\Gamma_0$ is also ICC,
so in particular it cannot have a non-trivial
normal subgroup $K$. This contradiction finishes the
proof. \hfill $\square$

\vskip .05in \noindent {\it End of proof of} 8.3. By Lemma 8.4,
$\theta(N_1), N^t_2$ are conjugate by a unitary element, so we may
assume $\theta(N_1)=N^t_2$. Thus, $\theta$ implements an
isomorphism between the groups $\Gamma_1=\Cal N_{M_1}(N_1)/\Cal
U(N_1)$ and $\Gamma_2=\Cal N_{M_2}(N_2)/\Cal U(N_2)$. But then by
the classical Kurosh Theorem and the condition on ``free
indecomposability'' of the groups $\Gamma_{ij}$, it follows that
$n_1=n_2=n$ and that there exists a permutation $\pi$ of the
indices $ 1 \leq j \leq n$ such that $g_j \Gamma_{1,j}g_j^{-1} =
\Gamma_{2,\pi(j)}$, for some elements $g_i\in G$. Thus,
$u_j=u_{g_j}$ normalizes $N_2^t$ and Ad$(u_j)$ takes
$\theta(M_{1,j})$ onto $M^t_{2,\pi(j)}$. \hfill $\square$

\vskip .05in \noindent {\it 8.5. Notation}. We denote by $f\Cal T_R$
the class of free actions $\sigma: \Gamma_0 * \Gamma_1 \rightarrow
\text{\rm Aut}(R)$ on the hyperfinite II$_1$ factor $R$, with the
properties:

\vskip 0.05in \noindent $(8.5.1)$. $\Gamma_0$ is free
indecomposable; $\Gamma_1$ is w-rigid (in particular free
indecomposable).

\vskip 0.05in \noindent $(8.5.2)$. $\sigma_0=\sigma_{|\Gamma_0}$
has the relative property (T), i.e. $R \subset R
\rtimes_{\sigma_0} \Gamma_0$ is a rigid inclusion;
$\sigma_1=\sigma_{|\Gamma_1}$ is a non-commutative Bernoulli shift
action of $\Gamma_1$ on $R = \overline{\otimes}_g (N_0,
\tau_0)_g$, where $N_0=R$ or $N_0 = M_{n\times n}(\Bbb C)$ for
some $n \geq 2$.

\vskip 0.05in \noindent $(8.5.3)$. $\sigma(\Gamma_1)$ and the
normalizer of $\sigma(\Gamma_0)$ in Out$(R)$ (which is countable
by [P5]) are free independent.

\vskip .1in

\proclaim{8.6. Lemma} Let $\Gamma_1$ be an arbitrary  w-rigid
group and $\Gamma_0=SL(n, \Bbb Z)$, $n \geq 2$, more generally
$\Gamma_0$ a free indecomposable, arithmetic lattice in an
absolutely simple, non-compact Lie group with trivial center. Then
$\Gamma_0 * \Gamma_1$ has $f\Cal T_R$ actions on $R$.
\endproclaim
\noindent {\it Proof}. By ([Va]) any such $\Gamma_0$ has a free
ergodic action on some $\Bbb Z^m$ such that the pair $(\Bbb Z^m
\rtimes \Gamma_0, \Bbb Z^m)$ has the relative property (T) of
Kazhdan-Margulis ([Ma]). By ([Ch], [NPS]) it follows that $\Gamma_0$
admits a free action $\sigma_0$ on the hyperfinite II$_1$ factor $R$
such that $R \subset R \rtimes_{\sigma_0} \Gamma_0$ has the relative
property (T). By ([P5]) it follows that the normalizer $\Cal N_0$ of
$\sigma_0(\Gamma_0)$ in Out$(R)$ is countable. Let $\sigma'$ be a
fixed copy Bernoulli shift action of $\Gamma_1$. By 8.2 it follows
that there exists an automorphisms $\theta$ of $R$ such that
$\theta(\sigma'(\Gamma_1))\theta^{-1}$ is freely independent from
$\Cal N_0$.

Thus, if we denote by $\sigma$ the unique action of
$\Gamma=\Gamma_0 * \Gamma_1$ on $R$ given by $\sigma_{|\Gamma_0} =
\sigma_0$ and $\sigma_{|\Gamma_1}= \theta \sigma' \theta^{-1}$,
then all conditions $(8.5.1)-(8.5.3)$ are satisfied. \hfill
$\square$

\proclaim{8.7. Theorem} Let $\sigma_i: \Gamma_{i0} * \Gamma_{i1}
\rightarrow \text{\rm Aut}(R_i)$, $i=1,2$, be $f\Cal T_R$ actions
and denote $G_i=\Gamma_{i0} * \Gamma_{i1}$, $M_i = R_i
\rtimes_{\sigma_i} G_i$, $i=1,2$. If $\theta: M_1 \simeq M_2^t$ is
an isomorphism, for some projection $t > 0$ then $t=1$ and there
exist a unitary $u\in M_2$, a character $\gamma$ of $\Gamma_{20}$ and
isomorphisms $\delta : G_1 \simeq G_2$, $\Delta : R_1 \simeq R_2$
such that $\theta_0 = {\text{\rm Ad}}(u) \circ
\theta^\gamma \circ \theta$ satisfies
. Moreover, any other
isomorphism $\theta': M_1 \simeq M_2$ is a perturbation of
$\theta$ by an automorphism of $M_2$ of the form $\text{\rm Ad}(v)
\circ \theta^{\gamma'}$ for some $v \in \Cal U(M_2)$ and
$\gamma'\in \text{\rm Char}(G_2)$.
\endproclaim
\noindent {\it Proof}. Since $R_i \subset R_i \rtimes_{\sigma_i}
\Gamma_{i0}$ are rigid inclusions, $R_i \subset R_i
\rtimes_{\sigma_i} G_i = M_i$ are rigid as well. We can thus apply
Theorem 8.3 to get a unitary element $v\in M_2^t$ such that
$v(\theta(R_1))v^*=R_2^t$ and $v(\theta(M_{11}))v^*=(M_{2j})^t$ for
some $j\in \{0,1\}$, where $M_{ij}=R_i \rtimes \Gamma_{ij}$, $i=1,2,
j=0,1$. But by ([P5] or [P6]) $R_1 \subset M_{11}$ is not rigid
while $(R_2 \subset M_{20})^t$ is rigid, so the only possibility is
that $j=1$. Thus, $(R_1\subset M_{11})\simeq (R_2 \subset M_{21})^t$
and both inclusions come from crossed products associated with
non-commutative Bernoulli shift actions of w-rigid groups. By ([P3])
this implies $t=1$.

On the other hand, Ad$v \circ \theta$ implements an isomorphism
$\delta : \Gamma_{10} * \Gamma_{11} \simeq \Gamma_{20}
* \Gamma_{21}$, which takes $\Gamma_{11}$ onto $\Gamma_{21}$. By
Kurosh Theorem, $\delta(\Gamma_{10})=g\Gamma_{20}g^{-1}$ for some
$g\in G_2$. But by ([GoS]), the groups $g\Gamma_{20}g^{-1}$ and
$\Gamma_{21}$ can generate $\Gamma_{20} * \Gamma_{21}$ only if
$g=g_1g_2$ for some $g_i \in \Gamma_{2i}$, $i=0,1$. By conjugating
with $g$, we may thus assume the unitary $v$ is so that
$\theta_0=\text{\rm Ad}(v)\circ \theta$ implements an isomorphism
$\delta: G_1\simeq G_2$ that takes $\Gamma_{1j}$ onto
$\Gamma_{2j}$, $j=0,1$. Thus, after identifying $R_1$ with $R_2$
via $\Delta=(\theta_0)_{|R_1}$ and $G_1\simeq G_2$ via $\delta$,
we are left with finding all automorphisms $\alpha$ of $M_2$ that
take $R_2$ onto itself and take the canonical unitaries $u_g$ into
unitaries $w_gu_g$, $g\in G_2$, for some $w:G_2 \rightarrow \Cal
U(R_2)$ a 1-cocycle for $\sigma_2$.

By ([P3]) this implies $w$ is co-boundary when restricted to
$\Gamma_{21}$, i.e. $w_g=\sigma_g(w)w^*$, $\forall g\in
\Gamma_{21}$, for some unitary $w\in R_2$. Thus, by replacing
$\alpha$ by Ad$(v^*)\circ\alpha$, we may assume $\alpha (u_g)=u_g,
\forall g\in \Gamma_{21}$. Thus $\alpha_{|R_2}\in \text{\rm
Aut}(R_2)$ commutes with $\sigma_2(\Gamma_{21})$, while still
normalizing $\sigma_2(\Gamma_{20})$. But by condition $(8.5.3)$
the latter condition implies $\alpha_{|R_2}$ is freely independent
from $\sigma_2(\Gamma_{21})$. This contradicts the commutation
condition with $\sigma_2(\Gamma_{21})$, unless $\alpha_{|R_2}$ is
inner. By perturbing $\alpha$ by Ad$(w_0)$ for an appropriate
$w_0\in \Cal U(R_2)$ we may thus assume $\alpha_{R_2}=id_{R_2}$.
Thus $\alpha$ is given by a character of $G_2$. \hfill $\square$

\vskip .05in The above theorem shows in particular that the
fundamental group of any $f\Cal T_R$ factor $M= R \rtimes_\sigma
(\Gamma_0 * \Gamma_1)$, corresponding to an $f\Cal T_R$ action
$(\sigma, \Gamma_0*\Gamma_1)$, is trivial, while its Out-group is
equal to Char$(\Gamma_0)\times\text{\rm Char}(\Gamma_1)$. By Lemma
8.6 one can take $\Gamma_0=SL(n,\Bbb Z)$, which has only trivial
characters, thus making Aut$(M)=\text{\rm Char}(\Gamma_1)$, with
$\Gamma_1$ an arbitrary w-rigid group. For instance, one can take
$\Gamma_1=SL(3, \Bbb Z)\times H$, where $H$ is an arbitrary
discrete abelian group, thus getting Aut$(M)=\hat{H}$. We thus
obtain:

\proclaim{8.8. Corollary} Let $\sigma: \Gamma_0 * \Gamma_1
\rightarrow \text{\rm Aut}(R)$ be a $f\Cal T_R$ action and denote
$M=R \rtimes_\sigma (\Gamma_0 * \Gamma_1)$. Then we have:

$1^\circ$. $\mycal F(M)=\{1\}$ and $\text{\rm Out}(M) = \text{\rm
Char}(\Gamma_0) \times \text{\rm Char}(\Gamma_1)=\text{\rm
Out}(M^\infty)$.

$2^\circ$. Given any compact abelian group $K$, there exists
$(\sigma, \Gamma_0 * \Gamma_1)$ such that the the corresponding
$f\Cal T_R$ factor $M$ satisfies $\text{\rm Out}(M)=K=\text{\rm
Out}(M^\infty)$. For instance, if $\Gamma_0=SL(n,\Bbb Z)$,
$\Gamma_1 = SL(m,\Bbb Z) \times \hat{K}$ for some $n,m \geq 3$,
then $\text{\rm Out}(M) = K$.
\endproclaim

\vskip .1in \noindent {\bf 8.9. Remark}. One can use $A.3.2^\circ$
in lieu of A.2 in all the above proofs, to construct more II$_1$
factors with small, calculable symmetry groups. Thus, let $f\Cal
T'_R$ be the class of free actions $\sigma: \Gamma_0 * \Gamma_1
\rightarrow \text{\rm Aut}(R)$ on the hyperfinite II$_1$ factor $R$,
satisfying the properties:

$(a)$. $\Gamma_0, \Gamma_1$ free indecomposable, $\neq \Bbb Z$;

$(b)$. $R \subset R \rtimes_{\sigma_0} \Gamma_0$ is a rigid
inclusion and $\sigma_1=\sigma_{|\Gamma_1}$ is non-cocylce conjugate
to $\sigma_0=\sigma_{|\Gamma_0}$ (note that this is indeed the case
if $\sigma_1$ is a non-commutative Bernoulli shift);

$(c)$. $\sigma(\Gamma_1)$ is free independent with respect to the
set $\Cal N(\sigma_0(\Gamma_0))\cup \Cal
N^{op}(\sigma_0(\Gamma_0))$, consisting of all automorphisms and
anti-automorphisms of $R^\infty$ that normalize
$\sigma_0(\Gamma_0)$.

By $A.3.2^\circ$ and 8.6 it follows that there exist such actions
$\sigma$ for any linear group $\Gamma_0$ as in 8.6 and for any
$\Gamma_1$ free indecomposable. It then follows as in 8.8 that the
corresponding crossed product II$_1$ factors $M=R \rtimes_\sigma
(\Gamma_0 * \Gamma_1)$ satisfy Out$(M)=\text{\rm Char}(\Gamma_0
* \Gamma_1)$. Moreover, if $t \in \mycal F(M)$ and $\theta: M \simeq M^t$
then $\theta$ must normalize $\sigma(\Gamma_0)$, contradicting
condition $(c)$ above. Thus, $\mycal F(M)=\{1\}$. Notice that to
get this calculation we no longer have to use the results in
([P3]) on the fundamental group of $R \subset R\rtimes_{\sigma_1}
\Gamma_1$. Similarly, if $\alpha$ is an anti-automorphism of $M$
then the same argument shows that it must normalize
$\sigma_0(\Gamma_0)$, in contradiction with the choice $(c)$.
Altogether, this shows that in addition to the properties $\mycal
F(M)=\{1\}$, Out$(M)=\text{\rm Char}(\Gamma_0 * \Gamma_1)$, the
factors $M$ in the class $f\Cal T'_R$ have no anti-automorphisms
either. This provides a fairly large new family of factors with
this latter property, after Connes first examples in ([C4]). Thus,
if we choose the groups $\Gamma_0, \Gamma_1$ without characters,
e.g. $\Gamma_i = SL(n_i,\Bbb Z), n_i\geq 3$, then the resulting
factors $M$ have no  outer symmetries at all. Moreover, noticing
that $f\Cal T'_R$ factors are w-rigid, it follows that any factor
of the form $N
* M$, with $N$ a property (T) II$_1$ factor (e.g. $N=L(PSL(n,\Bbb Z))$, $n
\geq 3$) and $M \in f\Cal T'_R$ has all the above properties (as a
consequence of Theorem 6.3 and of the properties of $M$) and in
addition has no Cartan subalgebras (by [V2], cf. Remark 6.6)!

\heading Appendix: Constructing freely independent actions.
\endheading

We now prove the technical results needed in the proofs of
Propositions 7.3 and 8.2, which established the existence of free
actions of groups $\Gamma=\Gamma_0
* \Gamma_1 * ...$ on $A=L^\infty(X,\mu)$
and on $R$ with restrictions to $\Gamma_i$ isomorphic to
prescribed actions $(\sigma_i, \Gamma_i)$, $\forall i$.

More precisely, we prove that given any countable set
$\{\theta_n\}_n$ of properly outer automorphisms of $(X, \mu)$
(resp. of $R$) the set $\Cal V$ of $\theta \in \text{\rm
Aut}(X,\mu)$ (resp $\theta \in \text{\rm Aut}(R)$) with the
property that all alternating words $\theta_{i_0} \theta
\theta_{i_1} \theta^{-1} \theta_{i_2} \theta \theta_{i_3}
\theta^{-1} ... $ are properly outer is $G_\delta$-dense in
Aut$(X,\mu)$ (resp. Aut$(R)$). Writing $\Cal V$ as an intersection
of open sets $\Cal V_n$ is obvious, and the non-trivial part is to
show that each $\Cal V_n$ is dense. To prove the density, in the
commutative case (A.1) we use a maximality argument inspired from
([P7]), while in the hyperfinite case (A.2) we use directly a
result from ([P7]), not having to re-do such maximality argument.

The idea of using Baire category, in both the proof of A.1 and
A.2, was triggered by Remark 2.27 in ([MS]) and the Category Lemma
in ([To]). In fact, the commutative case A.1 below is essentially
contained in ([To]). We have included the complete proof, with a
different treatment of the ``density'', for the reader's
convenience.

\proclaim{A.1. Lemma} Let $(X, \mu)$ be a standard non atomic
probability space and $\{\theta_n\}_{n \geq 1}\in \text{\rm
Aut}(X,\mu)$ a sequence of properly outer m.p. automorphisms of
$(X, \mu)$. Denote by $\Cal V\subset \text{\rm Aut}(X,\mu)$ the
set of all $\theta\in \text{\rm Aut}(X,\mu)$ with the property
that $\theta_{i_0}\Pi^n_{j=1} (\theta \theta_{i_{2j-1}}
\theta^{-1} \theta_{i_{2j}})$ is properly outer, $\forall n \geq
1$, $i_1, i_2, ..., i_{2n-1} \in \{1,2,3,...\}$ and $i_0, i_{2n}
\in \{0,1,2,...\}$, where $\theta_0 =id_X$. Then $\Cal V$ is a
$G_\delta$ dense subset of $\text{\rm Aut}(X,\mu)$. In particular,
$\Cal V \neq \emptyset$.
\endproclaim
\noindent {\it Proof}. Let $A=L^\infty(X,\mu)$, $\tau= \int \cdot
\text{\rm d}\mu$. Denote by $\Cal F$ the set of all finite
partitions of the identity $\{p_i\}_i \subset \Cal P(A)$. If $\rho
\in \text{\rm Aut}(X,\mu)$ then we still denote by $\rho$ the
automorphism it implements on $(A,\tau)$. As usual, $\text{\rm
Aut}(A,\tau)$ is endowed with the topology given by pointwise $\|
\cdot \|_2$ convergence, with respect to which it is metrizable
and complete.

For each $\rho \in \text{\rm Aut}(A,\tau)$ denote $k(\rho) = \inf
\{\|\Sigma_i \rho(p_i)p_i\|_2 \mid \{p_i\}_i \in \Cal F\}$. Note
that $\rho$ is properly outer iff $k(\rho)=0$. Also, if $\Cal D_n$
denotes the set of $\rho \in \text{\rm Aut}(A,\tau)$ with $k(\rho)
< 1/n$ then $\Cal D_n$ is clearly open in Aut$(A,\tau)$. Given an
$n$-tuple $(\theta_1, ..., \theta_n) \subset \text{\rm
Aut}(A,\tau)$, we denote by $\Cal V_n = \Cal V(\theta_1,...,
\theta_n)$ the set of all $\rho \in \text{\rm Aut}(A,\tau)$ with
$\theta_{i_0}\Pi^l_{j=1} (\rho \theta_{i_{2j-1}} \rho^{-1}
\theta_{i_{2j}}) \in \Cal D_n$, for all $1 \leq l \leq n$ and all
choices $i_1, i_2, .., i_{2l-1} \in \{1,2,..., n\}$, $i_0, i_{2l}
\in \{0, 1, 2, ..., n\}$.

It is immediate to see that $\Cal V_n$ is open in Aut$(A,\tau)$
and that $\cap_n \Cal V_n = \Cal V$. We have to prove that each
$\Cal V_n$ is also dense in Aut$(A,\tau)$, i.e., that any fixed
$\rho'\in \text{\rm Aut}(A, \tau)$ can be approximated arbitrarily
well (in the point norm-$\|\cdot \|_2$ convergence on $A$) by some
$\rho \in \Cal V_n$.

By replacing if necessary $\{\theta_j\}_j$ by the properly outer
automorphisms $\{\theta_j\}_j \cup \{\rho' \theta_k
{\rho'}^{-1}\}_k$, it follows that in order to prove the density
of $\Cal V_n$ it is sufficient to prove that $id_A$ is in the
closure of $\Cal V_n=\Cal V(\theta_1,..., \theta_n)$, for any
$n$-tuple of properly outer automorphisms.

To this end we will use the ultrapower II$_1$ factor $R^\omega$
([McD]) as framework. Thus, we choose $\omega$ a free ultrafilter
on $\Bbb N$ and let $R^\omega=\ell^\infty(\Bbb N, R)/\Cal
I_\omega$, where $\Cal I_\omega$ is the ideal associated with the
trace $\tau_\omega((x_n)_n) = \underset n \rightarrow \omega \to
\lim \tau(x_n)$, i.e. $\Cal I_\omega = \{x=(x_n)_n \mid
\tau_\omega(x^*x)=0\}$.

We regard the abelian von Neumann algebra $(A,\tau)$ as a Cartan
subalgebra of $R$. By ([Dy]) given any $\rho \in \text{\rm Aut}(A,
\tau)$, any finite set $F \subset A$ and $\varepsilon > 0$, there
exists $v\in \Cal N_R(A)$ such that $\|\rho(a)-\text{\rm
Ad}(v)(a)\|_2 \leq \varepsilon$, $\forall a\in F$. Thus, there
exist $u_n \in \Cal N_R(A)$ such that $u=(u_n)_n \in R^\omega$
satisfies $uau^*=\rho(a), \forall a\in A\subset A^\omega \subset
R^\omega$. Note that $\rho$ is properly outer iff $E_{A'\cap
R^\omega}(u)=0$.

Write $A=\overline{\cup_m D_m}^w$, for some increasing sequence of
finite dimensional subalgebras of $A$. Let $\Cal N_m \overset
\text{\rm def} \to = D_m'\cap \Cal N_R(A)$, i.e., the part of the
normalizer of $A$ in $R$ that leaves $D_m$ pointwise fixed. To
prove that $id_A$ is in the closure of $\Cal V_n$, it is
sufficient to prove:

\vskip .1in \noindent (A.1.1). Let $U_0=1, U_1, U_2, ..., U_n \in
\Cal N_{R^\omega}(A)$ with $E_{A'\cap R^\omega}(U_i)=0, \forall
i\neq 0$. For all $m \geq 1$, there exists $u\in \Cal N_m$ such that
$E_{A'\cap R^\omega}(u_{i_0}\Pi_{j=1}^l
(uU_{i_{2j-1}}u^*U_{i_{2j}})) =0$, for all $1 \leq l \leq n$ and all
choices $i_1, i_2, .., i_{2l-1} \in \{1,2,..., n\}$, $i_0, i_{2l}
\in \{0, 1, 2, ..., n\}$.

\vskip .1in

We construct $u$ by a maximality argument, ``patching together''
partial isometries in $\Cal G_m \overset \text{\rm def} \to = \{
vp \mid v\in \Cal N_m, p\in \Cal P(A)\}=D_m'\cap \Cal G\Cal
N_R(A)$.

Thus, for each $v\in \Cal G_m$ with $vv^*=v^*v$ and each $1 \leq k
\leq 2n$ we denote $\Cal V_k^v$ the set of all products of the form
$\Pi_{j=0}^{k-1} (U_{i_j}v^{\alpha_j})U_{i_k}$, where $i_j \in
\{1,2,..., n\}$, $\forall 1\leq j \leq k-1$, $i_0, i_k \in
\{0,1,...,n\}$, $\alpha_j \in \{\pm 1\}$, $\alpha_1 \neq \alpha_2
\neq ... \neq \alpha_k$. We also put $\Cal V^v_0=\{U_i \mid 1\leq
i\leq n\}$. Note that if $i_0, i_{2l} \in \{0,1,...,n\}$, $i_1, i_2,
..., i_{2l-1} \in \{1,2,..., n\}$, then $U_{i_0} \Pi_{j=1}^l
(vU_{2j-1} v^* U_{2j}) \in \Cal V_{2l}^v$. We denote
$$
\Cal W = \{v \in \Cal G_m \mid vv^*=v^*v, E_{A'\cap
R^\omega}(x)=0, \forall x\in \Cal V^v_k, 1\leq k \leq 2n\}.
$$
We endow $\Cal W$ with the order given by: $v\leq v'$ if
$v=v'v^*v$. $(\Cal W, \leq)$ is clearly inductively ordered. Let
$v_0\in \Cal W$ be a maximal element. Assume $0\neq p=1-v_0v_0^*
\in \Cal P(A)$. Since $E_{A'}(x)=0$ for all $x\in \Cal V^{v_0}_k$,
$1\leq k \leq 2n$, it follows that $E_{A'}(pxp)=0$ as well. Since
all such elements $w=pxp$ satisfy $ww^*, w^*w\in A$ and
$wAw^*=Aww^*$, by ([Dy]) it follows that

\vskip .1in \noindent (A.1.2). $\exists$ $0\neq q_1 \in \Cal P(Ap)$
such that: $q_1D_m=\Bbb Cq_1$, $q_1 w q_1=E_{A'}(w)q_1=0$, $\forall
w\in \cup_{k=0}^{2n} \Cal V_k^{v_0}$. \vskip .1in

Let $v_1\in \Cal N_{q_1Rq_1}(Aq_1)$. Note that $v_1 \in \Cal G_m$,
$v_1v_1^*=v_1^*v_1=q_1$. Denote $u=v_0+v_1$. We will show that if
$v_1\in \Cal N_{q_1Rq_1}(Aq_1)$ is chosen appropriately then $u\in
\Cal W$, thus contradicting the maximality of $v_0$. Write
$$x=\Pi_{j=0}^{k-1} (U_{i_j}u^{\alpha_j})U_{i_k} =\sum_\beta
\Pi_{j=0}^{k-1} (U_{i_j}v^{\alpha_j}_{\beta_j})U_{i_k}=\sum_\beta
y_\beta,
$$
where the sum is taken over all choices
$\beta=(\beta_j)_{j=0}^{k-1} \in \{0,1\}^k$. We will show that
$v_1$ can be taken such that $E_{A'}(y_\beta)=0$, $\forall \beta$,
$\forall x\in \cup_{k=0}^{2n} \Cal V^u_{k}$. For
$\beta=(0,0,...,0)$ we have $y_\beta=\Pi_{j=0}^{k-1}
(U_{i_j}v^{\alpha_j}_0)U_{i_k}\in \Cal V_k^{v_0}$, so that
$E_{A'}(y_\beta)= 0$, by the fact that $v_0\in \Cal W$.

The $k$ terms $y_\beta$ corresponding to just one occurrence of
$v_1$ (i.e. $\beta=(\beta_1,..., \beta_k)$ with all $\beta_i=0$
except one), are of the form $w_0v_1w_1$, with $w_0, w_1\in
\cup_{j=0}^{k-1} \Cal V^{v_0}_j$. Thus, each $w_i$ satisfies
$w_iw_i^*, w_i^*w_i \in A$, $w_iAw_i^*=Aw_iw_i^*$, $i=0,1$. By
shrinking $q_1$ recursively, we may assume
$(w_0^*w_0)q_1(w_1w_1^*)$ is either equal to $0$ or to $q_1$, for
all such $y_\beta$, $\forall 0 \leq k \leq 2n$. For the $y_\beta$
for which $(w_0^*w_0)q_1(w_1w_1^*)=0$ we have $y_\beta=0$ and
there is nothing to prove. For the $y_\beta$ with
$(w_0^*w_0)q_1(w_1w_1^*)=q_1$, take $u_0, u_1 \in \Cal
N_{R^\omega}(A)$ such that $u_0q_1=w_0q_1, q_1u_1=q_1w_1$. Then
$E_{A'}(y_\beta)=E_{A'}(u_0v_1u_1)= u_0E_{A'}(v_1u_1u_0)u_0^*$, so
that $\|E_{A'}(y_\beta)\|_1 = \|E_{A'}(v_1u_1u_0)\|_1$. Shrinking
$q_1$ recursively again, all conditions so far are still satisfied
while we can assume $q_1u_1u_0q_1$ is either $0$ or an element in
$A'\cap R^\omega$, for all $u_0, u_1$, coming from all $w_0, w_1$
arising as above. Thus, if we take $v_1 \in \Cal
N_{q_1Rq_1}(Aq_1)$ to be properly outer, then in both cases
$E_{A'}(v_1u_1u_0)=0$, $\forall u_0, u_1$. We have thus shown that
$q_1, v_1$ can be chosen so that for all
$\beta=(\beta_1,...,\beta_k)$ having just one occurrence of $1$ we
have $E_{A'} (y_\beta)=0$ as well.

Finally, the $y_\beta$ with at least two occurrences of $v_1$ can
be written as $y_\beta=x_0 (v^{\alpha}_1w v^{\alpha'}_1)x_1$, with
$w\in \Cal V^{v_0}_l$, for some $0\leq l \leq k-2$, $\alpha,
\alpha'\in \{\pm 1\}$ and $x_i$ partial isometries. Thus
$$
\|E_{A'}(y_\beta)\|_1 \leq \|y_\beta\|_1 \leq
\|v^{\alpha}_1wv^{\alpha'}\|_1=\|q_1wq_1\|_1 =\|E_{A'}(w)q_1\|_1=0
$$
the latter equality by (A.1.2).

Altogether, $E_{A'}(x)=0$, $\forall x\in \cup_{k=0}^{2n} \Cal
V^u_{k}$, showing that $u\in \Cal W$. But this contradicts the
maximality of $v_0$. Thus $v_0$ must be a unitary element,
finishing the proof that $\Cal V_n$ is dense in Aut$(A,\tau)$ and
thus the proof of the statement. \hfill $\square$

\proclaim{A.2. Lemma} Let $\{\theta_n\}_{n \geq 1}$ be a sequence of
properly outer automorphisms of the hyperfinite $\text{\rm II}_1$
factor $R$. Denote by $\Cal V$ the set of all automorphisms $\theta$
of $R$ such that any automorphism of $R$ of the form
$\theta_{i_0}\Pi^n_{j=1} (\theta \theta_{i_{2j-1}} \theta^{-1}
\theta_{i_{2j}})$ is outer, $\forall n \geq 1$, $i_1, i_2, ...,
i_{2n-1} \in \{1,2,3,...\}$ and $i_0, i_{2n} \in \{0,1,2,...\}$,
where $\theta_0 =id_R$. Then $\Cal V$ is a $G_\delta$ dense subset
of $\text{\rm Aut}(R)$. In particular, $\Cal V\neq \emptyset$.
\endproclaim
\noindent {\it Proof}. Let $\{u_n\}_n$ be a sequence of unitariy
elements in $R$, dense in $\Cal U(R)$ in the norm $\|\cdot \|_2$ and
with each element repeated infinitely many times. For each $x\in R$
and $\rho \in {\text{\rm Aut}}(R)$ denote by $k(\rho,x)$ the unique
element of minimal norm $\| \cdot \|_2$ in $K(\rho,x)=
\overline{\text{\rm co}}^w \{\rho(v) x v^* \mid v\in \Cal U(R)\}$.
Let $\Cal D_n$ be the set of automorphisms $\rho$ of $R$ with the
property that $\|k(\rho, u_i)\|_2 < 1/n, \forall 1\leq i \leq n$.

We claim that $\Cal D_n$ is open in ${\text{\rm Aut}}(R)$. To see
this, let $\rho \in \Cal D_n$ and for each $ 1\leq i \leq n$ choose
$v^i_1, v^i_2, ..., v^i_{m_i} \in \Cal U(R)$ such that $\|m_i^{-1}
\Sigma_{j=1}^{m_i} \rho(v^i_j) u_i {v_j}^{i*} \|_2 < 1/n$. If
$\delta > 0$ is sufficiently small, then any $\rho'\in \text{\rm
Aut}(R)$ satisfying $\| \rho'(v^i_j)-\rho(v^i_j)\|_2 < \delta$,
$\forall 1\leq i\leq n$, $\forall 1\leq j \leq m_i$, will satisfy
$\|m_i^{-1} \Sigma_{j=1}^{m_i} \rho'(v^i_j) u_i {v_j}^{i*} \|_2 <
1/n$, implying that $\|k(\rho',u_i)\|_2 < 1/n$, $\forall 1\leq i
\leq n$, thus $\rho'\in \Cal D_n$.

Denote by $\Cal V_n=\Cal V_n(\theta_1, ..., \theta_n)$ the set of
all $\rho \in \text{\rm Aut}(R)$ with the property that \newline
$\sigma_{i_0}\Pi^l_{j=1} (\rho \theta_{i_{2j-1}} \rho^{-1}
\theta_{i_{2j}}) \in \Cal D_n$, for all $1 \leq l \leq n$ and all
choices $i_1, i_2, .., i_{2l-1} \in \{1,2,..., n\}$, $i_0, i_{2l}
\in \{0, 1, 2, ..., n\}$. Since $\Cal D_n$ is open and
$$\text{\rm Aut}(R) \ni \rho \mapsto
\theta_{i_0}\Pi^l_{j=1} (\rho \theta_{i_{2j-1}} \rho^{-1}
\theta_{i_{2j}})
$$
is a continuous map, for each $m$ and $i_j$ as before, it follows
that $\Cal V_n$ is open in Aut$(R)$.

It is immediate to see that $\cap_n \Cal V_n = \Cal V$. Thus, in
order to show that $\Cal V$ is a $G_\delta$-dense subset of Aut$(R)$
we have to prove that each $\Cal V_n$ is dense in Aut$(R)$.
Moreover, arguing as in the proof of A.1 we see that, by replacing
if necessary $\{\theta_n\}_{n \geq 1}$ with the sequence
$\{\theta'_n\}_n = \{\theta_n\}_n \cup \{\rho \theta_m
\rho^{-1}\}_m$, it is enough to show that $id_R$ is in the closure
of $\Cal V_n$. We will prove that in fact $id_R$ is in the closure
of $\Cal V_n\cap \text{\rm Int}(R)$, i.e., given any finite
dimensional subfactor $R_0\subset R$, there exists $u\in \Cal
U(R_0'\cap R)$ such that Ad$(u) \in \Cal V_n$. In turn, this will be
an easy consequence of Theorem 2.1 in [P7].

To make the ideas more transparent, let us consider first the case
when $\{\theta_n\mid n\geq 0\}$ is an enumeration of the
automorphisms of a free cocycle action of a countable group $\Gamma$
on $R$, $\theta:\Gamma \rightarrow \text{\rm Aut}(R)$. Let
$M=R\rtimes \Gamma$, $U_0=1, U_1, U_2, ... \in \Cal U(M)$ be the
canonical unitaries implementing $\theta_0, \theta_1, ...$. Since
the action is free (i.e. $\theta_g$ non-inner $\forall g\neq e$), we
have $R'\cap M=\Bbb C$. Fix a free ultrafilter $\omega$ on $\Bbb N$.
We view $R\subset M$ as subalgebras of constant sequences in the
ultrapower II$_1$ factor $M^\omega$.

Since $(R_0'\cap R)'\cap M=R_0$, by Theorem 2.1 in [P7] there exists
a unitary element $V \in (R_0'\cap R)^\omega\subset M^\omega$ such
that $VRV^* \vee M \simeq R *_{R_0} M$. In particular, if $w \in
\Cal U(R_0'\cap R)$ is a Haar unitary and we put $U=VwV^*\in VRV^*$,
then for any choice of $1\leq l \leq n$, $i_1, i_2, ..., i_{2l-1}
\in \{1,2,...,n\}$, $i_0, i_{2l}\in \{0,1,2,..., n\}$ and $1\leq r
\leq n$ the unitary elements
$$
\text{\rm Ad}(U_{i_0} \Pi_{j=1}^l ((U U_{i_{2j-1}} U^{-1})
U_{i_{2j}})) (U^k) u_r U^{-k}, k=1,2, ... \tag A.2.1
$$
are mutually orthogonal with respect to the scalar product given by
the trace. To see this, we need to show that $\tau(\text{\rm
Ad}(U_{i_0} \Pi_{j=1}^l ((U U_{i_{2j-1}} U^{-1}) U_{i_{2j}})) (U^k)
u_r U^{-k})=0$, $\forall k \neq 0$. This amounts to
$\tau(U^{-k}\Pi_{j=1}^l (U U_{i_{2j-1}} U^{-1} U_{i_{2j}}) U^k
\Pi_{j=1}^l (U^{-1}_{i_{2l-2j}}U U^{-1}_{i_{2l-2j+1}}
U^{-1})u_r)=0$, which does indeed hold true, because after some
appropriate word-reduction we are left with a word of alternating
``letters'' $U^k \in VRV^* \ominus R_0$, $\forall k\neq 0$, and
$U_{i_j}, U^*_{i_j} \in M\ominus R$ $\subset M \ominus R_0$, for all
$i_j \neq 0$.

By the orthogonality of the elements in $(A.2.1)$, it follows that
for large enough $N$ we have
$$
\| N^{-1} \Sigma_{k=1}^N \text{\rm Ad}(U_{i_0} \Pi_{j=1}^l ((U
U_{i_{2j-1}} U^{-1}) U_{i_{2j}})) (U^k) u_r U^{-k}\|_2 < 1/n, \tag
A.2.2
$$
for all choices of $1\leq l \leq n$, $i_1, i_2, ..., i_{2l-1} \in
\{1,2,...,n\}$, $i_0, i_{2l}\in \{0,1,2,..., n\}$ and $1\leq r \leq
n$. Writing $U$ as a sequence of unitary elements in $R_0'\cap R$,
$U=(v_m)_m$, it follows that for $m$ large enough $u=v_m$ satisfies
$$
\| N^{-1} \Sigma_{k=1}^N \theta_{i_0} \Pi_{j=1}^l (\text{\rm Ad}(u)
\theta_{i_{2j-1}} \text{\rm Ad}(u^*) \theta_{i_{2j}})(u^k) u_r
u^{-k}\|_2 < 1/n, \tag A.2.3
$$
for all $1\leq l \leq n$, $i_1, i_2, ..., i_{2l-1} \in
\{1,2,...,n\}$, $i_0, i_{2l}\in \{0,1,2,..., n\}$, $1\leq r \leq n$.
Thus Ad$(u)\in \Cal V_n$, finishing the proof of this particular
case.

Now, in the general case we can take $\theta_0=id_R, \theta_1,
\theta_2, ...$ to be a lifting in Aut$(R)$ of an injective group
morphism $\Gamma\rightarrow \text{\rm Out}(R)$, with $\Gamma$
generated by $n$ elements. Notice that the automorphisms $\theta_j
\otimes \theta_j^{op}$ on $R\overline{\otimes} R$ implement a
cocycle action $\tilde{\theta}:\Gamma \rightarrow \text{\rm
Aut}(R\overline{\otimes} R^{op})$ (see e.g. Section 3 in [P10]), so
we can consider the crossed product factor
$\tilde{M}=R\overline{\otimes} R^{op}\rtimes \Gamma$. Denote by
$U_n\in \tilde{M}$ the canonical unitaries implementing
$\theta_n\otimes \theta_n^{op}$. By [P10] again, we can view
$R\overline{\otimes} R \subset \tilde{M}$ as the symmetric
enveloping inclusion associated with a ``diagonal subfactor''
$N\subset R\simeq \Bbb M_{n+1}(N)$, with the embedding of $N$ given
by $x \oplus \theta_1(x)...\oplus \theta_n(x)$. Moreover, the
associated Jones tower $N\subset R \subset N_1 \subset ... \nearrow
N_\infty$ can be viewed as a sequence of subalgebras of $\tilde{M}$,
making a non-degenerate commuting square:

$$
\CD
\ \ R\overline{\otimes} R^{op} \ \ \ \subset\ \ \tilde{M}\\
\noalign{\vskip-6pt}
\ \ \ \cup\ \ \ \ \ \ \ \ \ \ \ \ \ \ \cup \\
\noalign{\vskip-6pt}
R\vee R'\cap N_\infty \subset\ N_\infty\\
\endCD
$$

As before, we view $\tilde{M}$ as the algebra of constant sequences
in $\tilde{M}^\omega$. Since $(R_0'\cap R)'\cap N_\infty=R_0\vee
R'\cap N_\infty$ and each $R'\cap N_k$ is finite dimensional, we can
apply Theorem 2.1 in [P7] to get a unitary element $V \in (R_0'\cap
R)^\omega\subset N_\infty^\omega$ such that $VN_\infty V^* \vee
N_\infty \simeq  N_\infty *_{R_0\vee R'\cap N_\infty} N_\infty$. By
the above commuting square, we then also have $V \tilde{M}V^* \vee
\tilde{M} \simeq \tilde{M} *_{R_0\vee R^{op}} \tilde{M}$.

Like before, take $w \in \Cal U(R_0'\cap R)$ to be a Haar unitary
and denote $U=VwV^*\in VRV^*$. For any choice of $1\leq l \leq n$,
$i_1, i_2, ..., i_{2l-1} \in \{1,2,...,n\}$, $i_0, i_{2l}\in
\{0,1,2,..., n\}$ and $1\leq r \leq n$, the unitary elements
$$
\text{\rm Ad}(U_{i_0} \Pi_{j=1}^l ((U U_{i_{2j-1}} U^{-1})
U_{i_{2j}})) (U^k) u_r U^{-k}, k=1,2, ... \tag A.2.1'
$$
are then mutually orthogonal with respect to the scalar product
given by the trace. Indeed, since $U^k \in VRV^*\ominus R_0\subset
VN_\infty V^* \ominus R_0\vee N'\cap N_\infty$, $\forall k\neq 0$,
and $U_{i_j}, U^*_{i_j} \in \tilde{M}\ominus R\vee R^{op}$, $\forall
i_j \neq 0,$ it follows that $\tau(\text{\rm Ad}(U_{i_0} \Pi_{j=1}^l
((U U_{i_{2j-1}} U^{-1}) U_{i_{2j}})) (U^k) u_r U^{-k})=0$, for $k
\neq 0$, showing the orthogonality in $(A.2.1')$.

Now, by the orthogonality of the elements in $(A.2.1')$, it follows
that for large enough $N$ we have
$$
\| N^{-1} \Sigma_{k=1}^N \text{\rm Ad}(U_{i_0} \Pi_{j=1}^l ((U
U_{i_{2j-1}} U^{-1}) U_{i_{2j}})) (U^k) u_r U^{-k}\|_2 < 1/n, \tag
A.2.2'
$$
for all choices of $1\leq l \leq n$, $i_1, i_2, ..., i_{2l-1} \in
\{1,2,...,n\}$, $i_0, i_{2l}\in \{0,1,2,..., n\}$ and $1\leq r \leq
n$. Writing $U$ as a sequence of unitary elements in $R_0'\cap R$,
$U=(v_m)_m$, it follows again that for large enough $m$ the unitary
element $u=v_m$ satisfies
$$
\| N^{-1} \Sigma_{k=1}^N \theta_{i_0} \Pi_{j=1}^l (\text{\rm Ad}(u)
\theta_{i_{2j-1}} \text{\rm Ad}(u^*) \theta_{i_{2j}})(u^k) u_r
u^{-k}\|_2 < 1/n, \tag A.2.3'
$$
for all $1\leq l \leq n$, $i_1, i_2, ..., i_{2l-1} \in
\{1,2,...,n\}$, $i_0, i_{2l}\in \{0,1,2,..., n\}$, $1\leq r \leq n$.
Thus Ad$(u)\in \Cal V_n$.

\hfill $\square$

\vskip .05in \noindent {\bf A.3. Remarks}. The proofs of both
Theorems A.1 and A.2 can of course be carried out without using
the hyperfinite II$_1$ factor $R$ and its ultrapower $R^\omega$ as
framework, working exclusively  in the space Aut$(X,\mu)$,
respectively Aut$(R)$. But while it is straightforward to re-write
the proof of A.1 this way, the proof of A.2 then becomes much more
tedious, as one can no longer use results from ([P7]). Instead,
one has to go through a similar maximality argument as in the
proof of A.1 but with the non-commutativity requiring some
complicated estimates, similar to (pp 189-192 in [P7]).

When written in the ``Aut$(X,\mu)$ framework'', a suitable
adaptation of the proof of A.1 shows the following result: \vskip
.05in \noindent {\bf A.3.1}$^\circ$. Let $(X, \Cal X, \nu)$ be a
standard probability space and denote by $\Cal A$ the set of all
measurable isomorphisms $\rho$ of $X$ into $X$ such that $\nu
\circ \rho$ is non-singular with respect to $\nu$. Let $\mu$ be
another measure (not necessarily finite) on $(X, \Cal X)$,
equivalent to $\nu$. Denote by Aut$(X, \mu)$ the $\mu$-preserving
automorphisms in $\Cal A$. If $\{\theta_n\}_{n \geq 1}\in \Cal A$
are properly outer then the set $\Cal V\subset \text{\rm
Aut}(X,\mu)$ of all $\theta\in \text{\rm Aut}(X,\mu)$ with the
property that $\theta_{i_0}\Pi^n_{j=1} (\theta \theta_{i_{2j-1}}
\theta^{-1} \theta_{i_{2j}})$ is properly outer, $\forall n \geq
1$, $i_1, i_2, ..., i_{2n-1} \in \{1,2,3,...\}$ and $i_0, i_{2n}
\in \{0,1,2,...\}$, where $\theta_0 =id_X$, is a $G_\delta$ dense
subset of $\text{\rm Aut}(X,\mu)$.

\vskip .05in In turn, the proof of the non-commutative case in the
``Aut$(R)$ framework'' can be adapted to show the following more
general result: \vskip .05in \noindent {\bf A.3.2}$^\circ$. For each
$n \geq 1$, let $\theta_n : R^\infty \rightarrow R^\infty$ be either
endomorphism or anti-endomorphism of the hyperfinite II$_\infty$
factor $R^\infty$ such that $\theta_n$ is outer, $Tr \circ \theta_n$
is a finite multiple of the trace $Tr$ and $\theta_n(R^\infty)'\cap
R^\infty$ is atomic, $\forall n$. Denote by $\Cal V$ the set of all
$Tr$-preserving automorphisms $\theta$ of $R^\infty$ such that any
product $\theta_{i_0}\Pi^n_{j=1} (\theta \theta_{i_{2j-1}}
\theta^{-1} \theta_{i_{2j}})$ is outer, $\forall n \geq 1$, $i_1,
i_2, ..., i_{2n-1} \in \{1,2,3,...\}$ and $i_0, i_{2n} \in
\{0,1,2,...\}$, where $\theta_0 =id_R$. Then $\Cal V$ is a
$G_\delta$ dense subset of $\text{\rm Aut}(R^\infty)$. \vskip .05in

The case $[R^\infty: \theta_n(R^\infty)]<\infty$, $\forall n$, of
A.3.2$^\circ$ follows easily from Theorem 2.1 in [P7], by arguing
exactly as in the proof of A.2 above: The only change in that
argument is the definition of the finite index subfactor $N\subset
R$, which will again be a ``diagonal'' inclusion, but with $R\simeq
N^t$ for some appropriate amplification of $N$ and $N$ embedded into
it by taking a partition $p_0, ..., p_{n}\in \Cal P(R)$ and defining
$N\simeq p_0Rp$, $t=\tau(p_0)^{-1}$, $N\hookrightarrow N^t=R$ by $x
+ \Sigma_{i=1}^n \theta_i(x)p_i$, $x\in p_0Rp_0$, where
$\tau(p_i)/\tau(p_0)=\text{\rm d}Tr\circ \theta_i/\text{\rm d}Tr$,
$\theta_i:p_0Rp_0 \rightarrow p_iRp_i$ being ``corners'' of the
endomorphisms $\theta_i:R^\infty \rightarrow R^\infty$. When the
resulting subfactor is extremal, the rest of the argument is
identical, while in case it is not extremal, then one replaces the
symmetric enveloping algebra of $N\subset R$ be an appropriately
defined enveloping algebra $\tilde{M}$, containing $N_\infty$ and
satisfying appropriate commuting square properties.

Note that this result shows in particular that given any two
subfactors of finite Jones index of the hyperfinite II$_1$ factor,
$P\subset R$ and $Q\subset R$, there exists a subfactor $N\subset R$
having standard invariant the ``free product'' of the standard
invariants $\Cal G_{P,R}, \Cal G_{Q,R}$, as considered in ([BJ]).

\heading  References\endheading

\item{[BeVa]} B. Bekka, A. Valette: {\it Group cohomology,
harmonic functions and the first $L^2$-Betti numbers}, Potential
Analysis, {\bf 6} (1997), 313-326.

\item{[BJ]} D. Bisch, V.F.R. Jones: {\it  Algebras associated to
intermediate subfactors}, Invent. Math. {\bf 128} (1997), 89-157.

\item{[Bo]} F. Boca: {\it Completely positive maps on amalgamated
product $C^*$-algebras.} Math. Scand., {\bf 72} (1993), no. 2,
212-222.

\item{[Bu]} M. Burger: {\it Kazhdan constants for} SL$(3,\Bbb Z)$,
J. reine angew. Math., {\bf 413} (1991), 36-67.

\item{[Ch]} M. Choda: {\it A continuum of non-conjugate property}
(T) {\it actions of} $SL(n,\Bbb Z)$ {\it on the hyperfinite}
II$_1$ {\it factor}, Math. Japon., {\bf 30} (1985), 133-150

\item{[C1]} A. Connes: {\it Classification of injective factors},
Ann. of Math., {\bf104} (1976), 73-115.

\item{[C2]} A. Connes: {\it Outer conjugacy classes of
automorphisms of factors}, Ann. Ecole Norm. Sup., {\bf 8} (1975),
383-419.

\item{[C3]} A. Connes: {\it A type II$_1$ factor with countable
fundamental group}, J. Operator Theory {\bf 4} (1980), 151-153.

\item{[C4]} A. Connes: {\it Sur la classification des facteurs de
type} II, C. R. Acad. Sci. Paris {\bf 281} (1975), 13-15.

\item{[CJ]} A. Connes, V. Jones: {\it Property $T$ for von Neumann
algebras.} Bull. London Math. Soc., {\bf 17} (1985), no. 1, 57-62.

\item{[D]} H. Dye: {\it On groups of measure preserving
transformations}, II, Amer. J. Math, {\bf 85} (1963), 551-576.

\item{[DyR]} K. Dykema, F R\u adulescu: {\it Compressions of free
products of von Neumann algebras.} Math. Ann., {\bf 316} (2000),
no. 1, 61-82.

\item{[FM]} J. Feldman, C.C. Moore: {\it Ergodic equivalence
relations, cohomology, and von Neumann algebras I, II}, Trans.
Amer. Math. Soc. {\bf 234} (1977), 289-324, 325-359.

\item{[Fe]} T. Fernos: {\it Kazhdan's Relative Property} (T): {\it
Some New Examples}, \newline math.GR/0411527.

\item{[Fu1]} A. Furman: {\it Orbit equivalence rigidity}, Ann. of
Math. {\bf 150} (1999), 1083-1108.

\item{[Fu2]} A. Furman: {\it Outer automorphism groups of some
ergodic equivalence relations}, preprint 2003.

\item{[G1]} D. Gaboriau: {\it Invariants $\ell^2$ de r\'elations
d'\'equivalence et de groupes},  Publ. Math. I.H.\'E.S. {\bf 95}
(2002), 93-150.

\item{[G2]} D. Gaboriau: {\it Examples of groups that are measure
equivalent to free groups}, math.DS/0503181.

\item{[GP]} D. Gaboriau, S. Popa: {\it An Uncountable Family of
Non Orbit Equivalent Actions of $\Bbb F_n$}, to appear in Journal
of the AMS, math.GR/0306011.

\item{[Ge1]} S.L. Gefter: {\it On cohomologies of ergodic actions
of a T-group on homogeneous spaces of a compact Lie group}
(Russian), in ``Operators in functional spaces and questions of
function theory'', Collect. Sci. Works, Kiev, 1987, pp 77-83.

\item{[Ge2]} S.L. Gefter: {\it Outer automorphism group of the
ergodic equivalence relation generated by translations of dense
subgroups of compact group on its homogeneous space}, Publ. RIMS,
Kyoto Univ. {\bf 32} (1996), 517-538.

\item{[GoS]} O. N. Golowin, L. E. Syadowsky: {\it \"{U}ber die
Automorphismengruppen der freien Produkte}, Rec. math. Moscou,
N.S. {\bf 4} (1938), 505-514.

\item{[H]} U. Haagerup: {\it An example of a
nonnuclear $C\sp{*}$-algebra, which has
the metric approximation property.}  Invent. Math. {\bf 50} (1979),
279-293.

\item{[dHV]} P. de la Harpe, A. Valette: ``La propri\'et\'e T de
Kazhdan pour les groupes localement compacts'', Ast\'erisque {\bf
175} (1989).

\item{[Hj]} G. Hjorth: {\it A lemma for cost attained}, preprint
UCLA 2002.

\item{[J]} V.F.R. Jones : {\it Index for subfactors}, Invent.
Math. {\bf 72} (1983), 1-25.

\item{[Jol]} P. Jolissaint: {\it Haagerup approximation property
for finite von Neumann algebras.} J. Operator Theory, {\bf 48}
(2002), no. 3, 549-571.

\item{[Ka]} E. Kaniuth: {\it Eberhard Der Typ der regularen
Darstellung diskreter Gruppen}, Math. Ann. {\bf 182} (1969)
334-339.

\item{[K]} D. Kazhdan: {\it Connection of the dual space of a
group with the structure of its closed subgroups}, Funct. Anal.
and its Appl., {\bf1} (1967), 63-65.

\item{[Ko]} H. Kosaki: {\it Free products of measured equivalence
relations}, Journal of Functional Analysis {\bf 207} (2004),
264-299.

\item{[LS]} R. Lyndon, P. Schupp: ``Combinatorial Group Theory'',
Springer-Verlag, \newline Berlin/Heidelberg/New York, 1977.

\item{[McD]} D. McDuff: {\it Central sequences and the hyperfinite
factor}, Proc. London Math. Soc. {\bf 21} (1970), 443-461.

\item{[M]} G. Margulis: {\it Finitely-additive invariant measures
on Euclidian spaces}, Ergodic. Th. and Dynam. Sys. {\bf 2} (1982),
383-396.

\item{[MoS]} N. Monod, Y. Shalom: {\it Orbit equivalence rigidity
and bounded cohomology}, \newline Preprint 2002, to appear in Ann.
of Math.

\item{[NPSa]} R. Nicoara, S. Popa, R. Sasyk: {\it On irrational
rotation} HT {\it factors},

\item{[O1]} N. Ozawa: {\it Solid von Neumann algebras.} Preprint
2003, to Appear in Acta Math., math.OA/0302082.

\item{[O2]} N. Ozawa: {\it A Kurosh type theorem for type II$\sb
1$ factors,} math.OA/0401121.

\item{[OP]} N. Ozawa, S. Popa: {\it Some prime factorization
results for} II$_1$ {\it factors}, Invent. Math. {\bf 156} (2004),
223-234.

\item{[Pe]} J. Peterson: {\it A} 1-{\it cohomology
characterization of property} (T) {\it in von Neumann algebras},
math.OA/0409527.

\item{[PeP]} J. Peterson, S. Popa: {\it On the notion of relative
property (T) for inclusions of von Neumann algebras}, J. Funct.
Anal., {\bf 219} (2005), 469-483.

\item{[P1]} S. Popa: {\it Orthogonal pairs of *-subalgebras in
finite von Neumann algebras.} J. Operator Theory., (1983), no. 9,
253-268.

\item{[P2]} S. Popa: {\it Correspondences}, INCREST preprint 1986,
unpublished,
http://www.math.ucla.edu/~popa/popa-correspondences.pdf.

\item{[P3]} S. Popa: {\it Some rigidity results for
non-commutative Bernoulli shifts}, J. Fnal. Analysis {\bf 230}
(2006), 273-328.

\item{[P4]} S. Popa: {\it Markov traces on  universal Jones
algebras and subfactors of finite index}, Invent. Math. {\bf 111}
(1993), 375-405.

\item{[P5]} S. Popa: {\it On a class of type II$_1$ factors with
Betti numbers invariants.} Ann. of Math {\bf 163} (2006), 809-899.

\item{[P6]} S. Popa: {\it Strong rigidity of II$_1$ factors
arising from malleable actions of  $w$-rigid groups} I, II, Invent.
Math., {\bf 165} (2006), 369-408 and 409-453.

\item{[P7]} S. Popa: {\it Free independent sequences in type
II$_1$ factors and related problems}, Asterisque {\bf 232} (1995),
187-202.

\item{[P8]} S. Popa: {\it Some computations of $1$-cohomology groups
and construction of non orbit equivalent actions}, Journal of the
Inst. of Math. Jussieu {\bf 5} (2006), 309-332.

\item{[P9]} S.Popa: {\it A unique decomposition result for} HT
{\it factors with torsion free core}, J. Fnal. Analysis {\bf 242}
(2007), 519-525 (math.OA/0401138).

\item{[P10]} S. Popa: {\it Some properties of the symmetric enveloping
algebras with applications to amenability and property T}, Documenta
Mathematica, {\bf 4} (1999), 665-744.

\item{[S]} Y. Shalom: {\it Measurable group theory}, in
``Proceedings of the European Congress of Mathematocs'' (Stockholm
2004), European Math Soc, Zurich 2005, 391-424.

\item{[Sh]} D. Shlyakhtenko: {\it On the classification of Full
Factors of Type} III, Trans. AMS {\bf 356} (2004) 4143-4159.

\item{[T]} E. Thoma: {\it Eine Charakterisierung diskreter Gruppen
vom Typ} I, Invent. Math. 6 {\bf 1968} 190-196.

\item{[To]} A. Tornquist: {\it Orbit equivalence and actions of}
$\Bbb F_n$, UCLA preprint 2004.

\item{[U1]} Ueda: {\it Amalgameted free product over a Cartan
subalgebra}, Pacific J. Math. {\bf 191} (1999), 359-392.

\item{[U2]} Ueda: {\it Notes on treeability and cost in operator
algebra framework}, preprint 2005, \newline
math.OA/0504262.

\item{[Va]} A. Valette: {\it Group pairs with relative property}
(T) {\it from arithmetic lattices}, preprint 2004 (preliminary
version 2001).

\item{[V1]} D. Voiculescu: {\it Symmetries of some reduced free
product C$^*$ algebras}, in ``Operator algebras and their
connections with topology and ergodci theory'', Springer Lecture
Notes in Math {\bf 1132} (1985), pp. 566-588.

\item{[V2]} D. Voiculescu: {\it The analogues of entropy and of
Fisher's information theory in free probability} II, Invent. Math.
{\bf 118} (1994), 411-440 .

\item{[Z1]} R. Zimmer: {\it Strong rigidity for
ergodic actions of seimisimple Lie groups},
Ann. of Math. {\bf 112} (1980), 511-529.

\item{[Z2]} R. Zimmer: ``Ergodic theory and semisimple groups'',
Birkha\"user-Verlag, Boston 1984.

\enddocument